\documentclass[11pt,a4paper]{article} 

\usepackage{booktabs} 
\usepackage{array} 
\usepackage{paralist} 
\usepackage{subfig} 

\usepackage{bm} 

\usepackage{amssymb}
\usepackage{amsmath}
\usepackage{amsthm}
\usepackage[utf8]{inputenc}
\usepackage[T1]{fontenc}
\usepackage{lmodern}
\usepackage{ae}
\usepackage{enumerate}
\usepackage{enumitem}
\usepackage{thmtools}
\usepackage{tikz}
\usetikzlibrary{arrows.meta}

\usepackage[all]{xy} 

\usepackage{svg} 
\usepackage{pdfpages} 

\usepackage{geometry} 
\usepackage{graphicx} 

\usepackage{fancyhdr} 
\usepackage{sectsty}
\allsectionsfont{\sffamily\mdseries\upshape} 

\usepackage[nottoc,notlof,notlot]{tocbibind} 
\usepackage[titles,subfigure]{tocloft} 

\usepackage{eso-pic} 

\usepackage{stmaryrd} 

\usepackage{hyperref} 
\usepackage{cleveref} 

\newcommand{\N}{\mathbb{N}}
\newcommand{\Z}{\mathbb{Z}}

\newcommand{\R}{\mathbb{R}}
\newcommand{\D}{\mathbb{D}} 
\newcommand{\C}{\mathbb{C}}
\newcommand{\Hh}{\mathbb{H}} 

\newcommand{\CZ}{\C / \Z}

\newcommand{\piZ}{\pi_\Z}

\newcommand{\id}{\operatorname{id}}

\newcommand{\Ss}{\mathcal{S}} 

\newcommand{\Con}{\mathcal{C}}

\newcommand{\h}{\mathbf{h}} 
\newcommand{\Dd}{\mathcal{D}} 

\newcommand{\bi}{\mathbf{i}}

\newcommand{\W}[1]{W({#1})}
\newcommand{\Wt}[1]{\wt W({#1})}

\newcommand{\Dom}{\on{Dom}}

\newcommand{\setof}[2]{\{#1\,;\,#2\}}

\numberwithin{equation}{section}

\newlength{\speclen}
\declaretheoremstyle[
  thmbox={
    style=S,
    bodystyle=\normalfont\noindent,
    nounderline,
    hskip=\speclen,
  },
]{def}
\declaretheorem[style=def,heading=Definition,
numberlike=equation]{definition}

\newtheorem{theorem}[equation]{Theorem}

\newtheorem{proposition}[equation]{Proposition}
\newtheorem{PropositionDefinition}[equation]{Proposition-Definition}
\newtheorem{lemma}[equation]{Lemma}
\newtheorem{corollary}[equation]{Corollary}
\newtheorem{compl}[equation]{Complement}

  \theoremstyle{definition}

\newtheorem{remark}[equation]{Remark}
\newtheorem{example}[equation]{Example}
\newtheorem{question}[equation]{Question}
\newtheorem{notation}[equation]{Notation}
\newtheorem{hypothesis}[equation]{Hypothesis}

\newcommand{\on}[1]{\operatorname{#1}}
\newcommand{\wt}[1]{\widetilde{#1}}
\newcommand{\ov}[1]{\overline{#1}}
\newcommand{\wh}[1]{\widehat{#1}}
\renewcommand{\Re}{\operatorname{Re}}
\renewcommand{\Im}{\operatorname{Im}}
\newcommand{\eps}{\epsilon}
\newcommand{\tends}{\longrightarrow}
\newcommand{\ext}{\mathrm{ext}}
\newcommand{\wDd}{\widetilde{\mathcal D}}

\usepackage{regexpatch}

\makeatletter
\xpatchcmd*\thmbox@start
{\vsize}
{\pagegoal}
{}{\fail}
\makeatother

\title{Horn maps of holomorphic functions locally pseudo-conjugate on their parabolic basins}
\author{Arnaud Chéritat, Dimitri Le Meur}
\date{} 

\begin{document}

\maketitle

\section*{Abstract}

The lifted horn map of a holomorphic function with a simple parabolic point is well known to be a complete local conjugacy invariant; this is a classical result proved independently by Écalle \cite{b:Ec}, Voronin \cite{b:Vor}, Martinet and Ramis \cite{b:MarRam}.
Lanford and Yampolski have shown in \cite{b:LY} that, if two functions $f_1, f_2$ with simple parabolic points at $z_1, z_2$ are globally conjugate on their immediate parabolic basins, with the conjugacy and its inverse continuous at $z_1$, resp.\ $z_2$, then their horn maps must be cover-equivalent: there are isomorphisms $\psi^+ : \mathcal{D}_1^+\to \mathcal{D}_2^+$ and $\psi^- : \mathcal{D}_1^-\to \mathcal{D}_2^-$ between the top and bottom connected components of their domains, and a translation $T$ on the cylinder, such that $\h_2\circ\psi^+ = T\circ \h_1$ and $\h_2\circ\psi^- = T\circ \h_1$ holds on these domains.
In this article, we introduce a notion of (semi) local conjugacy on immediate parabolic basins, which we call \emph{local pseudo-conjugacy} and which in particular does not make any continuity assumption, and show that the horn maps $\h_1$ and $\h_2$ satisfy the condition above if and only if the two functions $f_1, f_2$ are locally pseudo-conjugate.
This result is a first step to better understand invariant classes by parabolic renormalization.

\newpage

\tableofcontents

\newpage

\section{Introduction}

\subsection{Classical results}

Let $f$ be a holomorphic function with a parabolic point at $0$ of multiplier $1$ and one attracting axis.
This is amounts to the fact that $f$ admits a Taylor expansion at $0$ of the form $f(z) = z + a z^2 + o(z^2)$ where $a \in \C^*$. Consider the lifted horn map $h$ of $f$, whose definition is recalled in \Cref{sec:horn_maps}.
It commutes with the translation $z\mapsto z+1$ and has an expansion of the form $h(z)=z+c_\pm+o(1)$ as $|\Im(z)|$ tends to $+\infty$.
The \emph{germ} of $h$ at $+ i \infty$ (resp.\ $- i \infty$) is defined as its equivalence class under the relation: $g\sim h$ if and only if there exists $M>0$ such that the restrictions of $g$ and $h$ coincide on the half-plane of equation $\Im(z)>M$ (resp.\ $\Im(z)<-M$).
The pair of germs
\[ (h_g^+, h_g^-) \] (here the index $g$ denotes the initial of the word \emph{germ}) of $h$ at $\pm i \infty$ is well known to be a local conjugacy complete invariant; this is a classical result proved independently by Écalle \cite{b:Ec}, Voronin \cite{b:Vor}, Martinet and Ramis \cite{b:MarRam}, which we state below: 

\begin{theorem} \label{t:conjRevetLoc}
Let $f_1$, $f_2$ be holomorphic functions from open neighborhoods of $0 \in \C$ to $\C$, having the following Taylor expansions at $0$: \[f_i(z) = z + a_iz^2 + o(z^2)\] where $a_i \neq 0$. 

Suppose that $f_1$, $f_2$ are locally conjugate\footnote{I.e.\ that there exists a biholomorphism $\phi$ between open subsets of $\C$ containing $0$, such that $\phi(0)=0$ and such that that $f_2 \circ \phi = \phi\circ f_1$ holds near $0$.} at $0$.
Then the germs at $\pm i \infty$, ${h_1}_g^\pm$, ${h_2}_g^\pm$ of the lifted horn maps $h_1$, $h_2$ are equivalent via translations.
More precisely, there exists $\sigma, \sigma' \in \C$ (independent of the sign $\pm$) such that: 
\[ {h_1}_g^\pm = T_\sigma^{-1} \circ {h_2}_g^\pm \circ T_{\sigma'} \]
where $T_\sigma$ denotes the translation of the plane $\C$ given by the formula $z \rightarrow z+\sigma$.

Conversely, if there exists $\sigma, \sigma' \in \C$ (independent of the sign $\pm$) such that, for the top and bottom germs,
\[ {h_1}_g^\pm = T_\sigma^{-1} \circ {h_2}_g^\pm \circ T_{\sigma'} \]
then $f_1, f_2$ are analytically conjugate in a neighborhood of $0$.
\end{theorem}

Local conjugacy of $f_1$ and $f_2$ near $0$ only depends on their germs ${f_1}_g^0$, ${f_2}_g^0$ at $0$, and the germs ${h_i}_g^\pm$ only depend on ${f_i}_g^0$, for $i\in\{1,2\}$.
For the next result, we cannot anymore identify maps with their restrictions to smaller domains.
A second difference is that it concerns non-lifted horn maps $\h_i$ instead of lifted horn maps $h_i$.
Let $\Dd^\pm \subset \CZ$ be the connected components of the domain of the non-lifted horn map $\h$ containing a punctured neighborhood of $\pm i \infty$.
Let $\h^\pm$ be the restriction $\h: \Dd^\pm \rightarrow \h(\Dd^\pm)$.
Denote by $B^0_{f_i}$ the immediate basin of $f_i$, i.e.\ the connected component of its parabolic basin $B_{f_i}$ that contains a germ of attracting axis.
In Proposition~3.16 of the book \cite{b:LY}, Lanford and Yampolsky
prove the following result.\footnote{Only the top $+$ of the cylinder is stated but the adaptation to the bottom $-$ is straightforward.}

\begin{theorem} \label{t:LY}
  Let $f_1$, $f_2$ be holomorphic functions from open neighborhoods of $0 \in \C$ to $\C$, with the following Taylor expansion at $0$: \[f_i(z) = z + a_iz^2 + o(z^2)\] where $a_i \neq 0$. 
  
  Suppose that $f_1$, $f_2$ are conjugate on their immediate basins via a biholomorphism $\phi: B^0_{f_1}\to B^0_{f_2}$.
  \emph{Assume further that $\phi$ and $\phi^{-1}$ both tend to $0$ at $0$.}
  Then there exists two biholomorphisms $\psi^+: \Dd_1^+\to \Dd_2^+$ and $\psi^-: \Dd_1^-\to \Dd_2^-$ such that
  \[ \h_1^\pm = T_\sigma^{-1} \circ \h_2^\pm \circ \psi^\pm \]
  and such that $\psi^\pm$ have removable singularities at $\pm i \infty$. 
\end{theorem}

In the present paper, we give two generalizations of the above theorem and reciprocals to these generalized statements.
In particular we characterize precisely when the conclusion of \Cref{t:LY} occurs by introducing the notion we call \emph{local pseudo-conjugacy}. In doing so the requirement that $\phi$ and $\phi^{-1}$ must extend continuously at $0$ disappear.
We also give a semi-conjugacy analogue of \Cref{t:LY}, which we will describe first and to which corresponds a notion we call \emph{local semi-conjugacy on immediate basins}.

\subsection{First statement: local semi-conjugacy on immediate basins}

We still assume that $f_1$ and $f_2$ are holomorphic maps defined in open neighborhoods of $0$ and with expansion $f_i(z)=z+a_iz^2+\ldots$ with $a_i\neq 0$.
For any other neighborhood $U_1$ of $0 \in \C$, we denote $U_1^0$ the connected component of $U_1 \cap B^0_{f_1}$ containing a germ of the attracting axis of $f_1$, where $B^0_{f_1}$ is the immediate parabolic basin of $f_1$.

\begin{definition}\label{d:semiConjLocIntro}%
One says that a holomorphic map $\phi$ \emph{locally semi-conjugates $f_1$ to $f_2$ on their immediate basins} if there exists an open neighborhood $U_1$ of $0 \in \C$ such that the domain of $\phi$ is $U_1^0$, $\phi$ takes values in $B^0_{f_2}$ and $\phi$ is a semi-conjugacy, namely: 
\[ \phi \circ f_1(z) = f_2 \circ \phi(z) \mbox{ for all $z \in U_1^0$ such that $f_1(z) \in U_1^0$.} \]
In this case we say that $f_1$, $f_2$ are locally semi-conjugate on their immediate basins, but note that this is not a symmetric relation.
\end{definition}

There are a few points to which we would like to attract the reader's attention: $U_1^0$ is \emph{not} defined as the immediate basin of the restriction of $f$ to $U_1$: the first is usually bigger than the second, i.e.\ in many cases some points in $U_1^0$ escape $U_1$.

It is also important to note that, though the domain of $\phi$ is contained in $U_1$, which may be chosen to be a small neighborhood of $0$, the range of $\phi$ is not required to be contained in a small neighborhood of $0$, i.e.\ $\phi$ is allowed to map points close to $0$ anywhere in the immediate basin of $f_2$.
There are situations where we belive it will have to (see the third example in \Cref{sub:examples}).

There is no requirement of injectivity nor surjectivity for local semi-conjugacies on immediate basins.
A somewhat trivial example is given by restrictions: if $f_1$ is a restriction of $f_2$, and $U_1$ is any open neighborhood of $0$ then $U_1^0\subset B^0_{f_1}\subset B^0_{f_2}$ and the map $\phi: z\mapsto z$ from $U_1^0$ to $B^0_{f_2}$ is a local semi-conjugacy on immediate basins as per \Cref{d:semiConjLocIntro}.

Here is the principal theorem of the paper. It is proven in \Cref{sec:proof} in two separate propositions (four propositions if one counts the complements): 
\begin{theorem} \label{t:conjRevetBasLoc}
Let $f_1$, $f_2$ denote holomorphic maps from open neighborhoods of $0 \in \C$ to $\C$, with Taylor expansions at $0$: $f_i(z) = z + a_iz^2 + o(z^2)$ where $a_i \neq 0$. 
Denote $\mathcal D_i^+$ the connected component containing $+i\infty$ of the domain of the non-lifted horn map $\h_i$ of $f_i$.
Similarly denote $\mathcal D_i^-$ the connected component containing $-i\infty$.
Denote $\h_i^+$ the restriction of $\h_i$ to $\mathcal D_i^+$ and $\h_i^-$ the restriction to $\mathcal D_i^-$.

Suppose that $f_1, f_2$ are locally semi-conjugate on their immediate basins at $0$ as per \Cref{d:semiConjLocIntro}.
Then there exists $\sigma \in \CZ$ and a pair of holomorphic maps $\psi = (\psi^+, \psi^-)$, where $\psi^\pm: \Dd_1^\pm \rightarrow \Dd_2^\pm$, such that: 
\[ \h_1^\pm = T_\sigma^{-1} \circ \h_2^\pm \circ \psi^\pm
\]

Conversely, if there exists $\sigma \in \CZ$ and a pair of holomorphic maps $\psi = (\psi^+, \psi^-)$, where $\psi^\pm: \Dd_1^\pm \rightarrow \Dd_2^\pm$, 
such that:
\[ \h_1^\pm = T_\sigma^{-1} \circ \h_2^\pm \circ \psi^\pm
\]
then $f_1, f_2$ are locally semi-conjugate at $0$ on their immediate basins. 
\end{theorem}

The maps $\psi^\pm$ are not necessarily injective nor surjective, but let us note that:

\begin{lemma}\label{lem:i2}
  Any holomorphic map $\psi^\pm:\wDd^{\pm}_1\to \wDd^{\pm}_2$ such that $\h_1^\pm = T_\sigma^{-1} \circ \h_2^\pm \circ \psi^\pm$ extends holomorphically into a map fixing $\pm i\infty$, and
   \[\psi^\pm = w + \rho^\pm + o(1)\]
  as $w\to \pm i\infty$, for some $\rho^\pm\in\C$.
\end{lemma}

This will be proved in \Cref{sec:horn_maps} as \Cref{lem:i2:copy}.

Moreover, we will see that the local semi-conjugacy $\phi$ and the maps $\psi^\pm$ are related via a commuting diagram:
\[
\xymatrix@R=10pt@C=35pt{
  \wDd_1^\pm \ar[dd]_{\wt\psi^\pm} \ar[r]^{\Psi_R^{1,\ext}} & B_{f_1}^0  \ar[r]^{\Phi_A^{1,\ext}} & \C \ar[dd]^{T_{\tilde \sigma}}
  \\
  & U_1^0 \ar[d]_{\phi} \ar[u]^{\iota} &
  \\
  \wDd_2^\pm \ar[r]_{\Psi_R^{2,\ext}} & B_{f_2}^0 \ar[r]_{\Phi_A^{2,\ext}} & \C
}
\]
where $\iota$ is an inclusion map.




\begin{remark}
  \Cref{t:conjRevetBasLoc} is an equivalence between two statements A and B.
  It is not clear that the hypotheses of A would be stable by composition: the composition $\phi_{3\leftarrow 2}\circ \phi_{2\leftarrow 1}$ of two local semi-conjugacies on immediate basins $\phi_{2\leftarrow 1}$ from $f_1$ to $f_2$ and $\phi_{3\leftarrow 2}$ from $f_2$ to $f_3$ is not necessarily a local semi-conjugacy on immediate basins from $f_1$ to $f_3$.
  (Counterexamples can be given using non-locally connected basins).
  However, the hypotheses of B are easily checked to be stable by composition: if $\h_1^\pm = T_\sigma^{-1} \circ \h_2^\pm \circ \psi_{2\leftarrow 1}^\pm$ and $\h_2^\pm = T_{\sigma'}^{-1} \circ \h_3^\pm \circ \psi_{3\leftarrow 2}^\pm$ then $\psi_{3\leftarrow 1}^\pm := \psi_{3\leftarrow 2}^\pm \circ \psi_{2\leftarrow 1}^\pm : \Dd_1^\pm \to \Dd_3^\pm$ satisfies $\h_1^\pm = T_{\sigma+\sigma'}^{-1} \circ \h_3^\pm \circ \psi_{3\leftarrow 1}^\pm$.
  This implies that being locally semi-conjugate on parabolic basins \emph{is a transitive relation}, which is not obvious from the definition.
  It would be nice to find an alternative notion of local semi-conjugacy $\tilde\phi_{j\leftarrow i}$, whose existence would be equivalent to the existence of a local semi-conjugacy on immediate basins from $f_i$ to $f_j$, but such that $\tilde\phi_{3\leftarrow 2}\circ \tilde\phi_{2\leftarrow 1}$ always satisfies the definition of the alternative notion.
  A track could be in the direction of the notion of \emph{pseudo-invertibility} introduced in \cite{t:LM}, pages~80 to~83.
\end{remark}

\Cref{t:conjRevetBasLoc} naturally adapts for maps $f_1, f_2$ with parabolic points with any number (possibly distinct for $f_1, f_2$) of petals and cycles of petals and for each of their horn maps. This generalization will be explicitly formulated in \Cref{sec:grl}.

\medskip

Notice that a local semi-conjugacy of immediate basins does not only preserve local properties. That is why one may wonder about their extension properties.
For instance, if the two functions $f_i: B^0_{f_i} \rightarrow B^0_{f_i}$ are proper and $B^0_{f_i}$ are simply connected, it is shown in \cite{b:DavidMorris} that a local conjugacy on all of a neighborhood of $0$ extends, even if it means post-composing it by $f_2^n$ (where $n \in \N$) before the extension, into a semi-conjugacy of the immediate basins.
Thus the following question is natural:

\begin{question}
  Let $\phi$ be a local semi-conjugacy of $f_1, f_2$ on their immediate basins. Does there exist $n \in \N$ such that $f_2^n \circ \phi$ extends into a global semi-conjugacy on the immediate basin ? 
\end{question}

The answer is no if we do not require anything more on $f_1$, $f_2$ than \Cref{d:semiConjLocIntro} (a counterexample is for instance given in \Cref{sub:ex:3:sub:1}).
However, for maps with simply connected immediate basin and that are proper on it as in \cite{b:DavidMorris}, it might be yes.

Unlike a local conjugacy on a whole neighborhood of $0$, a local semi-conjugacy on immediate basins has on one hand a rigidity introduced by the preservation of immediate basins, and on the other hand a flexibility introduced by the definition on smaller open sets and by the absence of an assumption of continuity on the boundary of the parabolic basins.
A local semi-conjugacy on immediate basins can modify topologically these boundaries, for instance by turning a non locally connected boundary into a locally connected boundary.
Such a semi-conjugacy, in a generalized definition presented in \Cref{sec:grl}, may also modify the number of petals of $f$.

\subsection{Second statement: local pseudo-conjugacy}

For the complement to the main theorem stated below, we need to introduce a few notions.
The quotient $B_0^{f_1}/f_1$ of $B_0^{f_1}$ by the relation $z\sim f_1(z)$ is known to be conformally isomorphic to the cylinder $\CZ$. The same holds for $f_2$.
A local semi-conjugacy on immediate basins $\phi$ as above (\Cref{d:semiConjLocIntro}) induces a map $\ov\phi : B^0_{f_1}/f_1 \rightarrow B^0_{f_2}/f_2$ in the following way: for $\overline{z} = z \bmod f_1 \in B^0_{f_1}/f_1$, there exists $N \in \N$ such that for all $n \geq N$, $f_1^n(z) \in U_1^0$.
The element $\phi(f_1^n(z)) \bmod f_2 \in B^0_{f_2} / f_2$ is independent of the representative $z$ of $\overline{z}$ and of the chosen integer $n \geq N$.
Then $\ov\phi(\ov z) := \phi(f_1^n(z))\bmod f_2$.
We will see that $\ov\phi $ is necessarily a biholomorphism in \Cref{p:SCQuotientBiholom}.

\begin{definition}\label{d:pseudoConjLocIntro}%
  Let $(\phi, \phi')$ denote a pair of local semi-conjugacies on immediate basins, i.e.\ $\phi: U_1^0 \rightarrow B^0_{f_2}$ is local semi-conjugacy on immediate basins from $f_1$ to $f_2$ and $\phi': U_2^0 \rightarrow B^0_{f_1}$ is a local semi-conjugacy on immediate basins from $f_2$ to $f_1$.
  We say that $(\phi, \phi')$ is a \emph{local pseudo-conjugacy} of $f_1, f_2$ if $\ov\phi : B_{f_1}/f_1 \rightarrow B_{f_2}/f_2$ has for inverse $\ov\phi ': B_{f_2}/f_2 \rightarrow B_{f_1}/f_1$.
\end{definition}

Though it is immediately seen to be reflective and symmetric, it is not obvious whether or not this relation is transitive, i.e.\ that it is an equivalence relation.
It turns out to be, as will follow from \Cref{prop:complement}.

\medskip

A particular case of local pseudo-conjugacy is when $f_1, f_2$ are globally conjugate on their respective parabolic immediate basins.

The complement below says that \Cref{t:conjRevetBasLoc} is still valid when replacing each occurrence of ``local semi-conjugacy on immediate basins'' by ``local pseudo-conjugacy'', and each occurrence of ``$\psi = (\psi^+, \psi^-)$ pair of holomorphic maps $\psi^\pm: \Dd_1^\pm \rightarrow \Dd_2^\pm$'' by ``$\psi = (\psi^+, \psi^-)$ pair of biholomorphisms $\psi^\pm: \Dd_1^\pm \rightarrow \Dd_2^\pm$''.
In this case, $\h_1^\pm$, $\h_2^\pm$ are equivalent via biholomorphism between $\Dd_1^\pm$ and $\Dd_2^\pm$ (more precisely via $\psi^\pm$), and a translation of $\CZ$ (more precisely via $T_\sigma$).
In particular if $\h_1^\pm$ is a ramified covering over the cylinder, then so is $\h_2^\pm$ (and conversely), and these ramified coverings are equivalent up to a translation on the range (the cylinder).

\begin{compl}\label{prop:complement}
  Let $f_1$, $f_2$ denote holomorphic maps from open neighborhoods of $0 \in \C$ to $\C$, with Taylor expansion at $0$: $f_i(z) = z + a_iz^2 + o(z^2)$ where $a_i \neq 0$. 
  
  Suppose that $f_1, f_2$ are locally pseudo-conjugate at $0$.
  Then there exists $\sigma \in \C$ and a pair of biholomorphisms $\psi = (\psi^+, \psi^-)$, $\psi^\pm: \Dd_1^\pm \rightarrow \Dd_2^\pm$, such that: 
  \[ \h_1^\pm = T_\sigma^{-1} \circ \h_2^\pm \circ \psi^\pm \]
  on $\Dd_1^\pm$.
  
  Conversely, if there exists $\sigma \in \C$ and a pair of biholomorphisms $\psi = (\psi^+, \psi^-)$, where $\psi^\pm: \Dd_1^\pm \rightarrow \Dd_2^\pm$ such that
  \[ \h_1^\pm = T_\sigma^{-1} \circ \h_2^\pm \circ \psi^\pm \]
  on $\Dd_1^\pm$
  then $f_1, f_2$ are locally pseudo-conjugate at $0$ on their immediate basins.
\end{compl}

As for the main theorem, the maps $\psi^\pm$ necessarily have at $\pm i \infty$ an expansion of the form $\psi^\pm(w) = w + \rho^\pm + o(1)$. 

\begin{remark}
  The first implication of the above theorem is already interesting
  if $\phi$ is a \textit{global} biholomorphism on the immediate basins: $B^0_{f_1}\to B^0_{f_2}$ and $\phi'$ is its inverse.
  It was proved in Proposition~3.16 of the book \cite{b:LY}, under this assumption and the supplementary hypothesis that
  $\phi$ and $\phi'$ are continuous at $0$.
  The construction of $\psi$ which will be exposed here becomes, without these continuity hypotheses, more difficult.
\end{remark}

\begin{remark}
  Let us make the following remark, whose proof we will omit.
  If we tried to make an analogue of \Cref{prop:complement} where we replace the maps $(\h^+, \h^-)$ by the germs $(\h^+_g, \h^-_g)$, we would get a statement that is a logical equivalence between the two following statements, which are easily seen to be always true:
  \begin{enumerate}
    \item there exists a $\sigma\in\C$ and invertible germs $\psi^{\pm}_g$ fixing the ends of the cylinder, such that ${\h_1}_g^\pm = T_\sigma^{-1} \circ {\h_2}_g^\pm \circ \psi^\pm_g$ (this is true with any $\sigma$, since the germs ${\h_i}_g^{\pm}$ are invertible);
    \item there is a pair of very large attracting petals\footnote{See \Cref{d:pet}} for $f_1$ and $f_2$, included in a small neighborhood of $0$ and on which the restrictions of $f_1$ and $f_2$ are conjugate (this is true since $f_1, f_2$ are conjugate to translations on appropriate very large petals).
  \end{enumerate}
  Whence the importance to work with the actual mappings $(\h^+, \h^-)$ and not only their germs, for the study of this more flexible relation.
\end{remark}

The classical \Cref{t:conjRevetLoc}, up to changing normalization of the Fatou coordinates, gives an equality between ${\h_1}_g^\pm$ and ${\h_2}_g^\pm$.
This enables to study the dynamic of ${\h}_g^\pm$ independently of the choice of the representative of $f$ modulo local conjugacy.
On the other hand, \Cref{prop:complement}, up to changing normalization of the Fatou coordinates, gives an equality modulo pre-composition by $\psi^\pm$ and post-composition with a translation.
This implies a preservation of ramified covering properties of $\h$, but not the dynamical ones.

It has interesting consequences for the parabolic renormalization operator, consisting in associating to $f$ the extension $\mathcal{R}f$ fixing $0$ of $\lambda\cdot E\circ \h\circ E^{-1}$ where $E :\CZ \to \C^*$ is an isomorphism and $\lambda\in\C^*$ is chosen so that $(\mathcal{R}f)'(0)=1$.
Usually we restricting $\mathcal{R}f$ to the component of its domain that contains $0$.
By \Cref{t:conjRevetBasLoc}, two maps that are conjugate on their immediate basins have two equivalent parabolic renormalization as coverings, and under appropriate covering property hypotheses, this will enable us to get a conjugacy between these two renormalizations on their parabolic basins, giving sufficient hypotheses to apply the theorem again.
This allows to build an invariant space by parabolic renormalization, and to think about the dynamic of this operator. 
Its fixed points have been studied.
On these two aspects, see the work of Inou and Shishikura \cite{b:IS}, Lanford and Yampolsky \cite{b:LY} for the quadratic case, and of Chéritat \cite{b:CherRenormParab} for the unicritical case.

\subsection{Examples}\label{sub:examples}

For illustrative purposes, we give an auxiliary lemma of classification which enables to apply \Cref{t:conjRevetBasLoc} to functions that are unisingular on their simply connected immediate basin.

This following lemma is stated as Theorem~6 in \cite{b:CherRenormParab} and references for its proof given there are Theorem 2.9 in \cite{b:LY} and \cite{b:DHOrsay}, Exposé~IX.

\begin{lemma}\label{l:blaschConj}
Let $f$ be holomorphic with a simple parabolic point at $0$. Assume that its immediate basin $B^0_f$ is simply connected and that the restriction $f:B^0_f\to B^0_f$ has a unique singular value $v$.

For integer $d \geq 2$, define the following holomorphic functions from $\D$ to $\D$: 
\[B_d(z) = \left(\frac{z+a}{1+az} \right)^d \mbox{ where $a = a_d = \frac{d-1}{d+1}$ } \]
and:
\[ B_{\infty}(z) = \exp\Big(2 \frac{z-1}{z+1}\Big) \]
Then $f: B^0_f\to B^0_f$ is conjugate to $B_d:\D\to\D$, where $d \geq 2$ (possibly infinite) is equal to the degre of the covering $f: B^0_f-f^{-1}(\{v\}) \rightarrow B^0_f-\{v\}$.
\end{lemma}

This gives the following conjugacies:

\begin{itemize}
\item  The parabolic basin of the map $\tan: \C \rightarrow \hat\C$ is $\C \setminus \R$. The map $\tan$ is conjugate on each of its immediate basins $ \pm \Hh$ to $B_\infty$. More precisely, the map $\tan$ is conjugate to $f: z \rightarrow \exp(2 \frac{z-1}{z+1})$ over all $\wh\C$ via the homography $\mu: z \rightarrow i \frac{z-1}{z+1}$: $\tan = \mu \circ f \circ \mu^{-1}$.
\item The cauliflower map $f_0: z \rightarrow z + z^2$ is conjugate on its immediate basin to the parabolic Blaschke product $B_2$ (passing from an immediate basin with fractal boundary to $\D$).
\item $z \rightarrow \exp(z) - 1$ and $z \rightarrow z \exp(z)$ are respectively analytically conjugate on their immediate basins to $ B_\infty$ and $B_2$ (passing from an immediate basin with non locally connected boundary to $\D$).
\item Let $f$ be a rational map which has only the following non-repelling cycles: a fixed parabolic point with only one petal and a Cremer cycle.
Suppose that the immediate parabolic basin is simply connected, totally invariant (the basin equals the immediate basin) and that $f$ is unicritical of degree $d$ over this immediate parabolic basin.
Then $f$ is conjugate to $B_d$ over its immediate basin, of which the boundary is non locally connected (passing from an immediate basin whith non locally connected boundary to $\D$).
\end{itemize}

\Cref{t:conjRevetBasLoc} allows in each case to reduce the study of the horn map modulo analytic ramified covering equivalence of the aforementioned functions ($\tan$, $f_0$, etc.) to horn maps of the functions $B_d$, $d \in \llbracket 2, + \infty \rrbracket$ ($B_{\infty}$, $B_2$, etc.). These last horn maps have the advantage of having a simple domain of definition: the unit disk. 

The third example above cannot be treated by direct application of \cite{b:LY} due to the absence of 
continuity of $\phi^{-1}:\D\to B^0_f$ at $z=1$, but can be with \Cref{t:conjRevetBasLoc}.
Actually in the third example, the map $\phi^{-1}$ is known to send any neighborhood of $0$, intersected with the the basin $\D$ of $B_d$, to an unbounded subset of $\C$.
Any other local semi-conjugacy on immediate basins will likely have the same behavior.

\subsection{Structure of the paper}

In \Cref{sec:rappel}, we recall a definition of Fatou coordinates and petals which is not fully compliant with \cite{b:Mil}, and their classical consequences.
In \Cref{sec:cyl}, we define the Fatou coordinate extensions and the attracting, repelling cylinders.
In \Cref{sec:horn_maps}, we define the horn map.
The content of \Cref{sec:rappel,sec:cyl,sec:horn_maps} are classical, though its presentation may be unusual in some places.
In \Cref{sec:defs_csq}, we define the pseudo-conjugacies and we prove some properties, which will allow in \Cref{sec:proof} to show the main theorem (\Cref{t:conjRevetBasLoc}), using an auxiliary result covered in \Cref{sub:aux}.
\Cref{sec:ex} provides examples and counterexamples of various statements appearing in the text.
We give in \Cref{sec:grl} the definition of generalized pseudo-conjugacy for parabolic points with several petal cycles, and the generalization of the theorem in this setting. 

Some of the methods and some of the notations are adapted from \cite{b:LY}.

\section{Parabolic point theory}

We recall here the part of the basic theory of parabolic points that we need, adapting the notions when needed. This section also serves to fix notation and terminology.
The results of this section are already known, although the presentation might be original in some parts.

\medskip

In this whole article, when a function $f$ is applied to a set, we do not assume that the set is contained in the domain of $f$: so $f(X)$ refers to $f(X\cap\Dom(f))$.

\subsection{Fatou coordinates, petals}\label{sec:rappel}

Let $f: \Dom (f) \subset \C \to \C$ a holomorphic function defined on a open neighborhood $\Dom (f)$ of $0$. We assume that $f$ has a parabolic point at $0$ with multiplier $1$ and one attracting axis, in other words $f$ has the Taylor expansion:
\[f(z) = z + a z^2 + O(z^3)\] 
with $a \in \C^*$. 

\ %

\begin{definition}%
We will call such a parabolic point a \emph{simple} parabolic point.
\end{definition}

The change of coordinates $w = -1/az$ conjugates $f$ to a map $F$ defined in a neighborhood of $\infty$ and satisfying $F(w) = w+1+o(1)$ as $|w|\to \infty$. From this it is easy: 1.~to build traps for $F$ that are right half-planes of the form $\Re(w)>x_0>0$ and to prove that any orbit of $F$ that tends to $\infty$ without hitting $\infty$ must satisfy $F^n(w) \sim n$,
2.~to prove that every orbit of $F$ tending to $\infty$ either hits $\infty$ or enters the trap.
The image in $z$ coordinate of the trap is a disk $D$ with $0$ in its boundary.
In particular every point $z$ whose orbit under $f$ tends to $0$ but without hitting $0$ must satisfy
\begin{equation}\label{eq:sim}
  f^n(z)\sim \frac{-1}{an}
\end{equation}

The estimate on $F$ implies the following result, which we will use in \Cref{subsub:Dd1Dd2}:
\begin{lemma}\label{lem:cptUnif}
  If $K$ is a compact subset of $\C$ on which all iterates of $f$ are defined, and if $f^n$ converges uniformly to $0$ on $K$ but never hits $0$, then there exists $N\geq 0$ such that $f^N(K)\subset D$.
\end{lemma}
\begin{proof}
  Let $R>0$ such that $|F(w)-(w+1)|<1/2$ for $|w|>R$. In particular $\Re F(w) > \Re(w) + 1/2$ for $|w|>R$.
  Let $N_0$ be such that $f^{N_0}(K)\subset B(0,1/R)-\{0\}$.
  Let $s(z)=1/z$ and $K' = s(f^{N_0}(K))$, which is compact.
  Let $x_1=\min_{w\in K'}{\Re w}$.
  Then $\forall w\in K'$ and $\forall k\geq 0$, $\Re(F^k(w))\geq x_1+k/2$, so $f^{n}(K)\subset D$ as soon as $n=k+n_0$ with $x_1+k/2> x_0$.
\end{proof}

The attracting axis of $f$ is defined by the equation $az^2 /z \in \R_{<0}$, and consists in the half line of points with argument $-\arg(a) + \pi$. Any orbit that tends to $0$ without hitting $0$ must have argument tending to the argument of the attracting axis. Similarly the repelling axis is defined by $az^2 /z \in \R_{>0}$ and has argument $-\arg(a)$.

The \emph{basin} of $0$ is denoted $B_f$ and is defined as the set of $z\in \C$ such that $f^n(z)\tends 0$ without hitting $0$. It is the union of the iterated preimages of $D$ by $f$, so it is an open set.
Let $B^0_f$ be the connected component of $B_f$ containing a germ of the attracting axis. The set $B^0_f$ is also named the \emph{immediate basin}.

\begin{lemma}\label{lem:hyp}
  The basin is a hyperbolic Riemann surface.\footnote{This lemma extends to more general situations that ours (in this article we require the domain and range of $f$ to be contained in $\C$). In general, the only maps with a parabolic point whose basin is not hyperbolic are conjugate to the translation by $1$ on the whole Riemann sphere.}
\end{lemma}
\begin{proof}
  It is enough to prove that it omits at least two point of $\C$.
  By definition, $0\notin B_f$.
  If $0$ were the only omitted point, then $f$ would be defined on $\C$ and $B_f = \C-\{0\}$.
  Consider the unit disk $\D$.
  Its boundary is compact and every of its point eventually enters $D$ and by continuity, a whole neighborhood of this point eventually is contained in $D$ after a finite number of iterates of $f$.
  By a finite covering, one sees the whole boundary of $\D$ is eventually mapped in $D$ by some $f^n$. By the maximum principle $f^n(0)\in D$, contradicting $f^n(0)=0$.
\end{proof}

Here is a definition of attracting petal that comes from \cite{b:LY} with slight modifications. Among others, the condition $\overline{f(P_A)} \setminus \{ 0 \} \subset P_A$ is replaced by $f(P_A) \subset P_A$.
This modified definition will prove to be useful to show that a local semi-conjugacy on immediate basin, even if not assumed continuous at 0, sends an attracting petal on an attracting petal, see \Cref{p:imSCPetAttr0}.
This is not a standard definition, it differs for instance from the more classical definition of Milnor \cite{b:Mil}.
Nevertheless, a petal in the sense of Milnor is a petal in our sense.
The classical existence results are hence still valid concerning our definition.

We denote $T_1(z) = z+1$.
\begin{definition}\label{d:pet}%
An attracting petal $P_A$ of $f$ is an open subset of $\C^*$ such that: 
\begin{enumerate}[label=(\arabic*)]
\item $P_A$ is connected, simply connected (i.e.\ $P_A$ is homeomorphic to $\D$).
\item $f(P_A) \subset P_A$.
\item $\forall z \in P_A, \lim_{n\to \infty} f^n(z) = 0$.
\item Conversely, if $f^n(z)$ converges to $0$, it may be stationnary at $0$, or else $f^n(z)$ must belong to $P_A$ for some $n$.
\end{enumerate}
If the set satisfies conditions (2) to (4), we call it a \emph{quasi-petal}. 
A Fatou coordinate on $P_A$ is a holomorphic function $\Phi_A: P_A \to \C$ such that:
\begin{enumerate}[resume,label=(\arabic*)]
  \item We have the following conjugacy: $\forall z \in P_A, \Phi_A \circ f (z) = T_1 \circ \Phi_A (z)$ for all $z\in\C$.
  \item $\Phi_A$ is injective on $P_A$.
\end{enumerate}
\end{definition}

Petals are special cases of quasi-petal.
Quasi-petals are never empty: indeed parabolic points always have a non-empty basin, so there is at least one non-stationary orbit satisfying (4).

\medskip

The set $\Phi_A(P_A)$ will appear a lot in this article, so we use a special notation for it:
\[Q_A := \Phi_A(P_A),\]
which we will recall every time.


According to \cite{b:Mil}, what we call Fatou coordinates were introduced by Leau and Fatou.
Following \cite{b:Mil} Lemma~10.10 page~114, we have:
\begin{proposition}\label{prop:mil}
  There exists an attracting petal $P_0$ for $f$, with Fatou coordinate $\Phi_0$, such that:
  \begin{enumerate} 
    \item The image of $P_0$ under $z\mapsto u=-1/az$ is the right half plane of equation $\Re(u)>x_0$ for some $x_0>0$.
    \item $\Phi_0(z) \sim - \frac{1}{a z}$ as $z\to 0$, $z\in P_0$.
  \end{enumerate}
\end{proposition}
The second point implies, by the argument principle, that the image $\Phi_0(P_0)$ contains sectors of the form $\arg(z-x(\eps))\in(-\pi/2+\eps,\pi/2-\eps)$ for all $\eps>0$ and some $x(\eps)$.\footnote{$x(\eps)$ may or may not tend to $+\infty$ as $\eps\to 0$.}
By~2.\ above, $z\tends 0$ $\implies$ $\Phi_0(z)\tends\infty$.
By making $P_A$ equal to the smaller set of equation  $\Re(-1/az)>x_0+1$ in 1.\ above, we can moreover ensure that and $\Phi_0(z)\tends\infty$ $\iff$ $z\tends 0$.

\begin{lemma}\label{lem:petint}
 The intersection of two attracting quasi-petals is a quasi-petal.
\end{lemma}
\begin{proof}
  Conditions (2) and (3) are immediate. Let us prove condition (4).
  We saw that quasi-petals are not empty.
  For any $z$ in one quasi-petal, by (3) its orbit tends to $0$ and does not hit 0 by (2).
  By (4) for the other quasi-petal, it eventually enters the other quasi-petal and by (2) for the first quasi-petal again, it is then in the intersection of the two quasi-petals.
\end{proof}

In particular this intersection is non-empty.
In particular two petals have non-empty intersection. However, the intersection is not necessarily a petal.

\medskip

Consider a quasi-petal $P_A$ and Fatou coordinate $\Phi_A : P_A\to\C$ as per \Cref{d:pet}.
Note that properties (2) and (5) imply
\[T_1(Q_A)\subset Q_A\]
where $Q_A=\Phi_A(P_A)$.

\begin{proposition}(Brimming\footnote{We introduce this terminology by comparison of the cylinder to a glass, $\Phi_A(P_A) \bmod \Z$ being the liquid.} property)\label{prop:brim}
  For all $w \in \C$, there exists $n \in \N$ such that $w + n \in Q_A$, where $Q_A=\Phi_A(P_A)$.
  Since $T_1(Q_A)\subset Q_A$, this amounts to: $Q_A \bmod \Z = \CZ$.
\end{proposition}
\begin{proof}
  Consider the petal and Fatou coordinate given in \Cref{prop:mil}.
  Let $Q_0 = \Phi_0(P_0)$.
  We define a map $\Theta: Q_0\to\C$ as follows: for $w\in Q_0$ the point $z=\Phi_0^{-1}(w)\in P_0$ has an orbit that tends to $0$ without hitting $0$, so must eventually enter $P_A$.
  Let $n\in\N$ such that $f^n(z)\in P_A$. We set $\Theta(w) = \Phi_A(f^n(z))-n$. Using that $\Phi_A$ conjugates $f$ to $T_1$, we see that $\Theta(w)$ is independent of the choice of $n$.
  As a consequence, it is holomorphic.
  Also, it is nowhere locally constant, since neither $\Phi_A$ nor the restriction of $f^n$ to $P_A$ are.
  It also satisfies $\Theta(w+1)=\Theta(w)+1$ so it projects to a mapping $\bm{\Theta}: \CZ\to\CZ$, which is also holomorphic.
  We now apply Liouville's theorem to deduce that $\bm{\Theta}$ is surjective: indeed, its domain $\CZ$ is not hyperbolic, but $\CZ$ minus any point is hyperbolic.
  It follows that $\Phi_A(P_A) \bmod \Z = \CZ$.
\end{proof}

\begin{remark} 
  The proof above has a nice interpretation in terms of the Riemann surface isomorphism class of the quotient $B_f/f$ of the basin under $z\sim f(z)$.
  See in \cite{t:LM} the proof of the first point of Proposition~4.1.2.
\end{remark}

\begin{proposition}\label{prop:brim2}
  Consider a quasi-petal $P_A$ with a Fatou coordinate $\Phi_A$.
  Consider any subset $X\subset\C$, not necessarily contained in $P_A$, and that $X$ satisfies the following condition (which is condition~(4) of \Cref{d:pet}): if $f^n(z)$ converges to $0$, it may is either stationnary at $0$, or else $f^n(z)\in X$ for some $n$.
  Then,\footnote{Recall that by our conventions, if $X\not\subset P_A$, $\Phi_A(X)$ refers to $\Phi_A(X\cap P_A)$.} $\Phi_A(X) \bmod \Z = \CZ$.
\end{proposition}
\begin{proof}
  Let $w\in \CZ$. By \Cref{prop:brim}, there exists $z\in P_A$ such that $\Phi_A(z) \bmod \Z = w$.
  Iterating forward $z$, by (3) for $P_A$ and (4) for $X$, we eventually enter $X$, while staying in $P_A$.
  By (5) for $P_A$, $\Phi_A(f^n(z))\bmod \Z = w$.
\end{proof}

\begin{proposition}[Uniqueness of attracting Fatou coordinates] \label{p:uniFatou} 
  Let $P^1$, $P^2$ be two attracting quasi-petals of $f$, and $\Phi_{1}$, $\Phi_{2}$ their Fatou coordinates. Then $\Phi_{1}$, $\Phi_{2}$ differ by a constant on $P^1 \cap P^2$.
\end{proposition}

\begin{proof}
  Let $P = P^1 \cap P^2$.
  Then $P$ is a quasi-petal.
  Note that $\Phi_1|_P$, $\Phi_2|_P$ are Fatou coordinates on $P$.
  
  We define a map $\alpha: \C\to\C$ as follows: given $w\in \C$, there exists by \Cref{prop:brim2} applied to $X=P$ some $n\in\N$ such that $w+n\in \Phi^1(P)$.
  Let $\alpha(w) = \Phi_2(\Phi_1^{-1}(w+n))-n$.
  It is independent of the choice of $n$, so the map $\alpha$ is holomorphic.
  It is also injective: if $\alpha(w)=\alpha(w')$ then we can choose a common $n$ such that $w+n$ and $w'+n \in \Phi^1(P)$,
  so $\Phi_2(\Phi_1^{-1}(w+n))-n = \Phi_2(\Phi_1^{-1}((w'+n))-n$ which by injectivity of $\Phi^2$ means $w=w'$.
  
  The only injective entire maps commuting with $T_1$ are the translations. This means $\Phi_2(\Phi_1^{-1}(w+n))-n = w+\sigma$ for some $\sigma$ independent of $w$ and $n$.
  For $z\in P$, apply this to $w=\Phi_1(z)$ and $n=0$ and get
  $\Phi_2(z) = \Phi_1(z)+\sigma$.
\end{proof}

Again there is a nice formulation of this proof in terms of quotients, see \cite{t:LM}, proof of Proposition~4.1.9.

\smallskip

In particular, Fatou-coordinates are unique on quasi-petals up to addition of a constant. Concerning existence:

\begin{proposition}\label{prop:exExistsFatou}
  On any quasi-petal $P_A$ such that $f$ is injective on $P_A$, there exists a Fatou coordinate.
\end{proposition}
\begin{proof}
  Uniqueness follows from \Cref{p:uniFatou}.

  For the existence, consider the petal $P_0$ and Fatou coordinate $\Phi_0$ given in \Cref{prop:mil}.
  Extension of $\Phi_0$ to a map $\Phi$ defined on the whole parabolic basin $B_f$ is classical: for any $z\in B_f$, there exists by (4) for $P_0$ some $n\in\N$ such that $f^n(z)\in P_0$ and we let $\Phi(f^n(z)) = z-n$, which is independent of the choice of $n$, hence holomorphic.
  It satisfies $\Phi(f(z))=\Phi(z)+1$ for all $z\in B_f$.
  
  We take $\Phi_A$ to be the restriction of $\Phi$ to $P_A$.
  There remains to check that $\Phi_A$ is injective. For this we use injectivity of $f$ on $P_A$: if $\Phi_A(z)=\Phi_A(z')$ take any $n\in\N$ such that both $f^n(z)$ and $f^n(z')$ belong to $P_0$. Then $\Phi_0(f^n(z))=\Phi_A(z)+n = \Phi_A(z')+n = \Phi_0(f^n(z'))$ so $f^n(z)=f^n(z')$. Since $f(P_A)\subset P_A$ and $f$ injective on $P_A$, $f^n(z)$ is injective on $P_A$, so $z=z'$.
\end{proof}

In particular for any petal, the Fatou coordinates on this petal are unique up to addition of a constant.

\begin{proposition} \label{p:covPetCompact}
  For any compact subset $K$ of $\C$, there exists $n \in \N$ such that $K+n \subset Q_A := \Phi_A(P_A)$.
\end{proposition}

\begin{proof}
The set $Q_A$ is open, so the map which from $z$ associates the least integer $\tau(z)$ such that $z + \tau(z) \in Q_A$ is upper semi-continuous. This implies that $\tau(z)$ admits a maximum $n$ in the compact set $K$. As $Q_A$ is forward invariant by the translation $T_1$, we have $K + n \subset Q_A$.
\end{proof}

\begin{proposition} \label{p:qpcp}
  Every quasi-petal contains a petal.
\end{proposition}
\begin{proof}
  Let $P_A$ be an attracting petal and denote $Q_A=\Phi_A(P_A)$.
  The set $Q_A$ is open, $T_{1}(Q_A)\subset Q_A$ and $\piZ(Q_A)=\C/\Z$.
  For every $k\in\Z$, consider the compact set $C_k = [0,1]\times[k,k+1]\subset \C$ (where we identified $\C$ to $\R\times \R$).
  By \Cref{p:covPetCompact}, there exists $n(k)\in\N$ such that $T_n(C_k)\subset Q_A$.
  Then actually $[n(k),+\infty)\times[k,k+1]\subset Q_A$.
  Choose now $h:\R\to\R$ such that $\forall k\in\Z$, $\forall y\in[k,k+1]$, $h(y)>n(k)$.
  Then the set $Q_A' = \setof{x+iy}{y>h(x)}$ is open, simply connected, contained in $Q_A$, $T_1(Q'_A)\subset Q'_A$ and: ($\ast$)every point of $\C$ has a forward $T_1$-orbit that eventually enters $Q'_A$.
  It follows that the set
  $P_A' =\Phi_A^{-1}(V)$ is an attracting petal: the first three points of \Cref{d:pet} are straightforward to check, and for the 4th one, once an orbit has entered $P_A$, it will enter $P'_A$ because $\Phi_A$ is a conjugacy and because of $\ast$.
\end{proof}

The proof of \Cref{prop:mil} classically extends to bigger domains (see for instance \cite{b:LY}, Proposition~2.5):
\begin{proposition}\label{prop:ePAa}
  For all $\alpha\in(0,\pi)$,
  there exists an attracting petal $P_A$ for $f$, with Fatou coordinate $\Phi_A$, such that:
  \begin{enumerate} 
    \item The image of $P_A$ under $z\mapsto -1/az$ contains the sector $S=\setof{w\in\C}{|\arg (w-w_0)|<\alpha}$ for some $w_0\in\C$,
    \item $\Phi_A(z) \sim - \frac{1}{a z}$ as $z\to 0$, $z\in P_A$ and $\Phi_A(z)\tends\infty$ $\iff$ $z\tends 0$.
  \end{enumerate}
\end{proposition}

As for \Cref{prop:mil}, the second point implies that the image $Q_A:=\Phi_A(P_A)$, at least \emph{for these specific petals}, contains sectors of the form  $\arg(z-x(\eps))\in(-\alpha+\eps,\alpha-\eps)$ for all $\eps>0$.

\begin{definition}\label{d:petale}%
Let $\alpha \in{} (0, \pi)$, and $P_A$ an attracting petal. We say that $P_A$ is an $\alpha$-petal, resp.\ large, very large, if: 
\begin{itemize}
\item ($\alpha$-petal) The image of $P_A$ by $\Phi_A$ equals a sector, more precisely: 
\[ \exists w_0 \in \C, \; \setof{w \in \C}{|\arg(w-w_0)| < \alpha } = \Phi_A(P_A)\]
\item (large) For all $\alpha \in ( 0, \pi )$, 
$\Phi(P_A)$ contains a sector as above.
\item (very large) The image of $P_A$ by $\Phi_A$ contains an upper half-plane, a lower half-plane, and a right half-plane. 
\end{itemize}
\end{definition}

\begin{definition}\label{def:rp}%
\emph{Repelling (quasi-)petals} of $f$ are defined as attracting (quasi-)petals for a local inverse $f^{-1}$ of $f$ which sends $0$ on $0$.
But concerning the Fatou coordinate, we post-compose it by $z \mapsto -z$, in such a way that it conjugates $f$ (instead of $f^{-1}$) to $T_1$, i.e.
\[ \Phi_R(z) = - \Phi_A^{f^{-1}}(z) \]
\end{definition}

All the properties of attracting petals above extend to repelling petals with straightforward modifications (in the dynamical plane of $f$, sectors in $P_R$ are bisected by the repelling axis instead of the attracting axis, images of $\alpha$-petals in Fatou are sectors whose bisector has argument $\pi$, etc.).

\begin{lemma}[Small petals]\label{lem:asp}
  For every neighbourhood $U$ of $0$ there are attracting and repelling petals of $f$ included in $U$.
\end{lemma}
\begin{proof}
  Just apply \Cref{prop:ePAa} to the restriction of $f$ to $U$ for attracting petals and to a local inverse for repelling petals. 
\end{proof}

Similarly, for all $\alpha\in(0,\pi)$, existence of arbitrarily small attracting $\alpha$-petal is ensured by \Cref{prop:ePAa} applied to $\alpha'\in(\alpha,\pi)$ together with the argument principle.
The same holds for repelling petals, using $f^{-1}$.

\begin{remark}
  By the uniqueness statement (\Cref{p:uniFatou}), a petal is large if and only if it contains an $\alpha$-petal for all $\alpha\in(0,\pi)$.
\end{remark}

\begin{lemma}[Asymptotics]\label{lem:asymp1}
  Let $P_A$ be any attracting petal and $\alpha\in(0,\pi)$.
Consider a sector $S = \setof{z\in \C}{ |\arg(-az)| < \alpha}$ (whose bisector is the attracting axis). Then
\[\Phi_A(z) \sim - \frac{1}{a z}\]
when $z\to 0$, $z\in S\cap P_A$.

  Let $P_R$ be any repelling petal and $\alpha\in(0,\pi)$.
Consider a sector $S = \setof{z\in \C}{|\arg(az)| < \alpha }$ (whose bisector is the repelling axis). Then
\[\Phi_R(z) \sim - \frac{1}{a z}\]
when $z\to 0$, $z\in S\cap P_R$.
\end{lemma}

\begin{proof}
  Choose any $\alpha'\in(\alpha,\pi)$.
  Consider a petal $P'_A$ as provided by \Cref{prop:ePAa} for this $\alpha'$, and its Fatou coordinate $\Phi_{A'}$ satisfying $\Phi_{A'}(z)\sim -1/az$ as $z\to 0$ in $P'_A$.
  The image of $P'_A$ by $z\mapsto -1/az$ contains a sector $\{w\in\C,|\arg (w-w_0)|<\alpha'\}$ which itself contains a sector of the form $S=\setof{w\in\C}{|w|>R \text{ and }|\arg w|<\alpha}$.
  So the set $P'_A$ contains the sector $S' = \setof{z\in \C}{|z| < 1/R \text{ and } |\arg(-az)| < \alpha}$.
  By the uniqueness statement (\Cref{p:uniFatou}), $\Phi_{A'}$ and $\Phi_A$ differ by a constant on $P'_A\cap P_A$, thus in particular on $S'\cap P_A$.
  The claim follows.
  The argument is identical for repelling petals.  
\end{proof}

We would like to deduce that $\Phi_A^{-1}(w) \sim -1/aw$ as $w\tends \infty$ within sectors. The deduction requires to pay attention to domains of validity, so we state this as a lemma too.

\begin{lemma}[Asymptotics, 2]\label{lem:asymp2}
Let $w_0 \in \C$, $R>0$ and $\beta \in (0, \pi)$. 

Let $P_A$ be an attracting petal and suppose that the sector $S = \setof{w \in \C}{|w-w_0| > R \text{ and } |\arg(w-w_0)| < \beta}$ is included in $Q_A := \Phi_A(P_A)$.
Then
\[\Phi_A^{-1}(w) \sim - \frac{1}{a w}\]
when $w \to \infty$, $w \in S$.

Let $P_R$ be a repelling petal and suppose that the sector $S = \setof{w \in \C}{|w-w_0| > R \text{ and } |\arg(-(w-w_0))| < \beta }$ is included in $Q_R = \Phi_R(P_R)$.
Then
\[\Phi_R^{-1}(w) \sim - \frac{1}{a w}\]
when $w \to \infty$, $w \in S$.  
\end{lemma}
\begin{proof}
  Choose a real $\alpha$, such that $\beta<\alpha<\pi$.
  Consider a petal $P'_A$ as provided by \Cref{prop:ePAa} for this $\alpha$, such that the image of $P'_A$ by $z\mapsto -1/az$ contains the sector $\setof{w\in\C}{|w-w'_0|>R\text{ and }|\arg (w-w'_0)|<\alpha}$ for some $w'_0\in \C$ and $R>0$, and its Fatou coordinate $\Phi_{A'}$ satisfying $\Phi_{A'}(z)\sim -1/az$ as $z\tends 0$ in $P'_A$.
  By \Cref{p:uniFatou}, $\Phi_{A'}$ and $\Phi_A$ differ by a constant on $P'_A\cap P_A$.
  By subtracting this constant to $\Phi_{A'}$ (and to $w'_0$), we can assume that the constant is $0$.
  By the argument principle, the set $Q'_A = \Phi_{A'}(P'_A)$ contains $S':=S \cap \setof{w\in\C}{|w-w_0|>R'}$ for some $R'\geq R$.
  The set $Q_A = \Phi_A(P_A)$ also contains $S'$.

  The crucial point is to prove the deceivingly obvious claim that $\Phi_{A'}^{-1}= \Phi_{A}^{-1}$ on $S'$.
  There are connected open susbets $U$ of $S'$ such that $T_1(U)\subset U$, for instance appropriate sectors.
  Choose any such $U$.
  Let $V=\Phi_A^{-1}(U)$.
  Then $f(V)\subset V$.
  Choose a point $z_0\in V$.
  Since $P'_A$ is a petal, the orbit $z_0$ eventually enters $P'_A$.
  So $V$ and $P'_A$ have a non-empty intersection.
  On it, $\Phi_A = \Phi_{A'}$ since $V\subset P_A$.
  So $\Phi_A^{-1} = \Phi_{A'}^{-1}$ on a non-empty open susbet of $S'$.
  By holomorphic continuation, $\Phi_A^{-1} = \Phi_{A'}^{-1}$ on $S'$.

  Finally, since $\Phi_{A'}(z)\sim -1/az$ on $P'_A$ and $\Phi_{A'}(z)\tends 0$ $\iff$ $z\tends 0$, we have $w\sim -1/a\Phi_{A'}^{-1}(w)$ as $w\tends \infty$ within $Q'_A$, i.e.\ $\Phi_{A'}^{-1}(w) \sim -1/aw$.
  It holds in particular on $S'$, on which $\Phi_{A'}^{-1} = \Phi_{A}^{-1}$.
  The proof for repelling Fatou coordinates is the same.
\end{proof}

Using the argument principle and calling \emph{small $\alpha'$-sector bisected by the attracting axis} any set of the form $\setof{z\in \C}{ 0 < |z| < r\text{ and } | \arg(- a z) | < \alpha'} \subset P_A$ for some $r>0$, we get:

\begin{corollary}\label{cor:sec1}
  An attracting/repelling $\alpha$-petal contains a small $\alpha'$-sector bisected by the attracting/repelling axis for all $0 < \alpha' < \alpha$.
  In particular, two large petals, one of which is attracting and the other repelling, always share a non-empty intersection.
\end{corollary}

\begin{corollary} \label{r:sectLarge}
An attracting petal $P_A$ is large if and only if it contains a small $\alpha'$-sector bisected by the attracting axis, for all $0 < \alpha' < \pi$.
A similar statement holds for repelling petals, replacing $\arg(- a z)$ by $\arg(a z)$.
\end{corollary}
\begin{proof}
  For the direct implication, proceed as in the previous statement. For the reciprocal, use \Cref{lem:asymp1} and again the argument principle.
\end{proof}

\begin{figure}[!t]
  \begin{center}
  \begin{tikzpicture}
    \node at (0,0) {\includegraphics[scale=0.6]{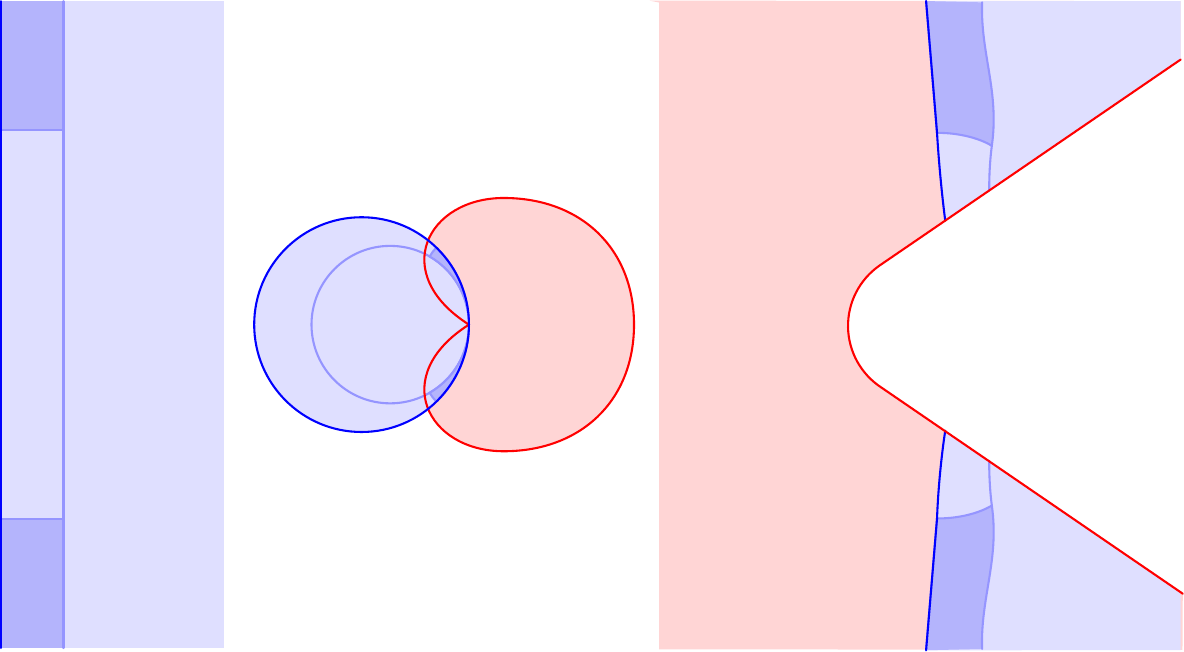}};
    \draw[Stealth-] (-4.3,0) -- node[above] {$\Phi_A$} (-3.2,0);
    \draw[-Stealth] (0,0) -- node[above] {$\Phi_R$} (1.2,0);
    \node at (-2.45,-1.7) {$P_A$};
    \node at (-0.6,-1.7) {$P_R$};
    \draw[|-|] (-6.4,0) -- node[left]{$M'$} (-6.4,2);
  \end{tikzpicture}
  \end{center}
  \caption{Illustration of \Cref{p:petalesAR}}
  \label{fig:A}
\end{figure}

\begin{lemma} \label{p:petalesAR}
Let $P_A$ be an attracting petal containing a $\frac{\pi}{2}$-petal, $C$ the inverse image of a strip $\setof{z\in \C}{ x_0 \leq \Re(z) < x_0 + 1 }$ by $\Phi_A$, so that $C$ is a fundamental domain of $P_A$ (i.e.\ a set containing exactly one point of each grand orbit of $f$ restricted to $P_A$), and $P_R$ a repelling petal containing a $\beta$-petal for some $\beta>\pi/2$.
Let $M > 0$.
Then there exists $M' > 0$ satisfying the following property: if $z \in C$ is such that $\Im( \Phi_A(z) ) > M'$, then $z \in P_R$ and $\Im(\Phi_R(z)) > M$. There also exists $M'$ satisfying the property: if $z \in C$ is such that $\Im( \Phi_A(z) ) < - M'$, then $z \in P_R$ and $\Im(\Phi_R(z)) < -M$.
See \Cref{fig:A}.

The property is also true if one permutes simultaneously each occurrence of repelling by attracting and vice versa. 
\end{lemma}

\begin{remark}
The constant $M'$ depends on the choice of $P_A$, $P_R$, and the choice of fundamental domain $C$.
\end{remark}

\begin{proof}
Let $w \in \C$ such that $x_0 \leq w < x_0+1$. The asymptotic-equivalent of $\Phi_A^{-1}$ (\Cref{lem:asymp2}) shows that $\Phi_A^{-1}(w) \sim - \frac{1}{a w}$ when $w \to +\infty$ within the strip.

When $\Im(w)$ is large enough (resp.\ large enough with negative sign), then we have $|\Phi_A^{-1}(w)|$ small and $\arg(\Phi_A^{-1}(w))$ close to $- \arg(a) - \pi/2$ (resp.\ close to $-\arg(a) + \pi/2$).

Recall that the attracting axis admits as argument $- \arg(a) + \pi$, the repelling axis as argument $- \arg(a)$. The quantity $\Phi_A^{-1}(w)$ is close to $0$, in a little cone containing a germ of the line $X$, where $X$ is the line passing through $0$ perpendiculary to both attracting and repelling axis. Since $P_R$ contains a $\beta>\pi/2$ petal, we may apply \Cref{r:sectLarge} to deduce that $\Phi_A^{-1}(w)$ must belong to $P_R$.

Now, if $z \in D$ is such that the quantity $\Im(\Phi_A(z))$ is big, resp.\ big with negative sign, we may apply the preceding argument to deduce that $z \in P_R$. Since $z$ is close to $0$ and avoids a sector centered on the attracting axis and a sector centered on the repelling axis, we have the asymptotic-equivalents $\Phi_A(z) \sim - \frac{1}{az} \sim \Phi_R(z)$, hence $\Im(\Phi_R(z))$ must be big (resp.\ big with negative sign) when $\Im(\Phi_A(z))$ is big (resp.\ big with negative sign). This finishes the proof.

The reasoning obtained inverting the role of attracting and repelling petals is identical. 
\end{proof}

\begin{proposition}\label{prop:petEx}
  The map $f$ admits attracting $\alpha$-petals for all $\alpha\in(0,\pi)$, large petals and very large petals.
\end{proposition}

\begin{proof}
Existence of the $\alpha$-petals follows from \Cref{prop:ePAa} as we already explained.
For the other types, it is sufficent to prove the existence of a very large petal, since it is a particular case of a large petal.
We only detail here the construction of very large attracting petals, the construction of very large repelling ones being similar.

\textit{Existence of very large attracting petals: }
Let $P_R$ be a repelling $\beta$-petal with $\beta>\pi/2$.
Let $P_A$ be an attracting $\pi/2$-petal of $f$.
By definition, $D = \Phi_A(P_A)$ is a right half-plane: $D = \setof{z\in \C}{ a < \Re(z) }$.
Let $B = \setof{z\in \C}{a < \Re(z) \leq a+1 }$.
For $M > 0$ we note $\Delta =  \setof{z\in \C}{ | \Im(z) | \leq M }$, $P_A^M = \Phi_A^{-1}(D \setminus \Delta)$ and $X^M = \Phi_A^{-1}(B \setminus \Delta)$.
If $M$ is big enough, $X^M\subset P_R$ by \Cref{p:petalesAR}.
We denote $g=f^{-1}: P_R \to f^{-1}(P_R)$.

Let $X^M_{-n} = g^{n}(X^M)$ for all $n \in \N$, and 
\[\begin{cases} 
P_A^L = P_A \cup \bigcup_{n \in \N^*} X^M_{-n}\\ 
D^L = D \cup (\C \setminus \Delta) = D \cup \bigcup_{n \in \N^*} \Big( \big(B - n \big) \setminus \Delta \Big)
\end{cases}
\]
Let $n \in \N^*$. Note that $X_{-n}^M$ is disjoint from $P_A$. Indeed, $\Phi_A( f^n(X_{-n}^M)) \subset B$ and $\Phi_A(f^n(P_A)) \subset D + n$, and $B$, $D+n$ are disjoint. By a similar reasoning, $X_{-n}^M$ and $X_{-m}^M$ are also disjoint.

Let for $z \in P_A^L$, $\Phi_A(z) = \Phi_A(f^n(z)) - n$ where $n$ is such that $f^n(z) \in P_A$. The value of the expression does not depend on the integer $n$ chosen so $\Phi_A$ is holomorphic.
Let us check that $\Phi_A$ is a biholomorphism from $P_A^L$ to $D^L$.
It is surjective since $\Phi_A(X^M_{-n}) = (B\setminus\Delta)-n$. The map $\Phi_A$ is injective on $P_A$, and on $X_{-(n+1)}^M$ for $n \in \N^*$. These sets are a partition of $P_A^L$ and have disjoint images by $\Phi_A$, this finally implies that $\Phi_A$ is injective. 

It remains to check that $P_A^L$ and $\Phi_A$ satisfy the definition of a petal and associated Fatou coordinate.

\begin{itemize}
\item $P_A^L$ is homeomorphic via $\Phi_A$ to $D^L$, which is homeomorphic to the unit disk.
\item $f(P_A^L) \subset P_A^L$, because $f(P_A) \subset P_A$ and $f(X_{-n}^M) = X_{-n+1}^M \subset P_A^L$ for $n \in \N^*$.
\item $\forall z \in P_A^L, \; \lim_{n\to \infty} f^n(z) = 0$. This property is indeed valid if $z \in P_A$, and if $z \in X_{-n}^M$, then $f^n(z) \in P_A$.
\item Conversely, if $f^n(z)$ converges to $0$, either it is stationary at $0$, or $f^n(z)$ belongs to $P_A$ for $n$ big enough, hence to $P_A^L \supset P_A$.
\item $\Phi_A$ is injective on $P_A^L$
\item (very large) The image of $P_A^L$ by $\Phi_A$ is a union of, so contains, an upper half-plane, a lower half-plane and a right half-plane.
\item We have the following conjugacy: $\forall z \in P_A, \Phi_A \circ f (z) = T_1 \circ \Phi_A (z)$, with $T_1(z) = z+1$ for all $z \in \C$.
\end{itemize}
\end{proof}

\begin{remark}\label{rk:smallVLPetals}
  Let $U$ be an open neighbourhood of $0$.
  By applying the previous proposition to a restriction of $f$ to $U$, we get that $f$ has very large attracting and repelling petals contained in $U$, i.e.\ $f$ has arbitrarily small very large attracting and repelling petals.
\end{remark}

\subsection{Extension of Fatou coordinates}\label{sec:cyl}

We say that two functions are \emph{equal in the strong sense} if they have the same domain and take the same value at every point of this domain.\footnote{Formally this is what equality of functions is supposed to mean. Here we use this terminology to insist on the importance of the equality of the domains.}

\subsubsection{Extension of attracting Fatou coordinates}

\begin{definition}\label{d:FatouEt}%
  Given an attracting Fatou coordinate $\Phi_A$ on an attracting petal $P_A$, one defines the corresponding \emph{extended attracting Fatou coordinate} of $f$, denoted $\Phi^\ext_A : B_f \to \C$, in the following way. For $z \in B_f$, there exists $n \in \N$ such that $f^n(z) \in P_A$. The quantity $w:= \Phi_A(f^n(z)) - n$ does not depend on the chosen integer $n$ and we let $\Phi_A^{\ext}(z) = w$.
\end{definition}

This also is the only function that extends $\Phi_A$ to $B_f$ and satisfies $\Phi_A^\ext \circ f = T_1\circ \Phi_A^\ext$. It is holomorphic.

\begin{lemma}\label{lem:paeor}
  For any two $z,z'\in B_f$,
  $\Phi_A^\ext(z) = \Phi_A^\ext(z')$ $\iff$ $\exists n\in\N$ $f^n(z)=f^n(z')$.
\end{lemma}
\begin{proof}
  If $f^n(z)=f^n(z')$ then $\Phi_A^\ext(z)+n = \Phi_A^\ext(f^n(z)) = \Phi_A^\ext(f^n(z')) = \Phi_A^\ext(z')+n$ so $\Phi_A^\ext(z) = \Phi_A^\ext(z')$.
  Conversely if $\Phi_A^\ext(z) = \Phi_A^\ext(z')$ then let $n$ be big enough so that $f^n(z)$ and $f^n(z')$ both belong to $P_A$.
  We have $ \Phi_A(f^n(z')) = \Phi_A^\ext(z')+n = \Phi_A^\ext(z)+n = \Phi_A(f^n(z))$.
  By injectivity of $\Phi_A$ on $P_A$, $f^n(z)=f^n(z')$.
\end{proof}

\begin{remark}
  Let $P_A^1, P_A^2$ be two attracting petals of $f$, with Fatou coordinates that coincide on their intersection.
  Then they give the same extended Fatou coordinate.
  This comes from the fact that $\Phi_A^{\ext}(z)$ equals $\Phi_A(f^n(z)) - n$ for all $n$ such that $f^n(z) \in P_A^1 \cap P_A^2$.
  In this sense, the extended Fatou coordinate does not depend on the choice of the petals.
  Two extended Fatou coordinates differ by a constant (post-composition by a translation).
  It also follows that an extended attracting Fatou coordinate is injective on any attracting petal.
\end{remark}

\begin{lemma}\label{lem:W}
  Denote $\piZ:\C\to\CZ$ the canonical projection and assume that:
  \begin{itemize}
    \item $W\subset B_f$ is open and simply connected,
    \item $f(W)\subset W$, 
    \item $\Phi_A^\ext$ is injective on $W$,
    \item and $\piZ \circ \Phi_A^\ext(W) = \CZ$.
  \end{itemize}
  Then $W$ is an attracting petal as per \Cref{d:petale}.
\end{lemma}
\begin{proof}
  Every point of the definition of a petal is straightforward to verify, except maybe point (4), i.e.\ that all orbit in $B_f$ eventually falls in $W$.
  Choose $P_A$ to be any attracting petal.
  Let $z_0\in B_f$, $z_n = f^n(z_0)$ and $w_0 = \Phi_A^\ext(z_0)$, so $(\Phi_A^\ext)^n(z_n) = w_0+n$.
  Since $P_A$ is a petal, by property (4) for $P_A$, for all $n$ big enough, $z_n\in P_A$.
  By hypothesis $\pi_\Z(w_0)\in\pi_\Z\circ \Phi^\ext_A(W)$, i.e.\ $w_0 \in T_{k}\circ \Phi_A^\ext(W)$ for some $k\in\Z$.
  Since $f(W)\subset W$, $\Phi_A^\ext(W)$ is $T_1$-stable.
  So $\exists N\in\N$ such that $\forall n\geq N$, $w_0 + n \in \Phi_A^\ext(W)$.
  For $n\geq N$ let $z'_n$ be the unique element of $W\cap (\Phi_A^\ext)^{-1}(w_0 + n)$.
  Then $f(z'_n)\in W$ and $\Phi_A^\ext(f(z'_n))=1+\Phi_A^\ext(z'_n) = 1+w_0+n = \Phi_A^\ext(z'_{n+1})$ so by injectivity of $\Phi_A^\ext$ on $W$: $f(z'_n) = z'_{n+1}$.
  So $(z'_n)_{n\geq N}$ is an orbit. Since it is contained in $B_f$, it tends to $0$. So it is in $P_A$ for $n$ big enough.
  So for $n$ big enough, both $z_n$ and $z'_n$ belong to $P_A$, and both have an image by $\Phi_A^\ext$ equal to $w_0+n$.
  By injectivity of $\Phi_A^\ext$ on $P_A$, $z_n=z'_n$.
  So $z_n\in W$.
\end{proof}

On $B_f$ we consider the equivalence relation generated by $f(z)\sim z$, i.e.\ the smallest equivalence relation $\sim$ such that $\forall z\in B_f$, $z\sim f(z)$.
Note that $f(z)\in B_f$.
It is a classical fact that, in this situation, $z\sim z' \iff (\exists n,m\in N)$ $f^n(z)=f^m(z')$ (see \cite{t:LM}, Remark~4.2.6 and more generally from Definition~4.2.5 to Remark~4.2.14.).

\begin{lemma}\label{lem:quot}
  $\piZ\circ \Phi_A^\ext(z)=\piZ\circ \Phi_A^\ext(z')$ $\iff$ $z\sim z'$.
\end{lemma}
\begin{proof}
  We can either imitate the proof of \Cref{lem:paeor} or use it.
  Let us choose the second option:
  $\piZ\circ \Phi_A^\ext(z)=\piZ\circ \Phi_A^\ext(z')$
  $\iff$ $(\exists k\in \Z)$ $\Phi_A^\ext(z)= \Phi_A^\ext(z')+k$ $\iff$ $(\exists a,b\in\N)$ $\Phi_A^\ext(z)+a= \Phi_A^\ext(z')+b$
  $\iff$ $(\exists a,b\in\N)$ $\Phi_A^\ext(f^a(z))= \Phi_A^\ext(f^b(z'))$ 
  $\iff$ (by \Cref{lem:paeor}) $(\exists a,b,c\in\N)$ $f^{c+a}(z)= f^{c+b}(z')$.
  $\iff$ $(\exists n,m\in\N)$ $f^n(z)= f^m(z')$.
\end{proof}

\begin{corollary}\label{cor:quot}
  The map $\piZ\circ \Phi_A^\ext$ quotients to a bijection from $B_f/f$ to $\CZ$.
\end{corollary}

\begin{remark}
  One can justify that the quotient topology of $B_f/f$ is Hausdorff and that one gets a Riemann surface, and that the bijection is an isomorphism to the Riemann surface $\CZ$.
  See \cite{t:LM} from Definition~4.2.5 to Remark~4.2.14.
  We will not use these facts here.
\end{remark}

\subsubsection{Extension of inverse repelling Fatou coordinates}

\begin{definition}\label{d:paramFEt}%
  Given a repelling Fatou coordinate $\Phi_R$ on a repelling petal $P_R$, one defines the \emph{extended repelling Fatou parametrization} of $f$ in the following way. Let $w \in \C$. By \Cref{prop:brim}, there exists $n \in \N$ such that $w-n \in Q_R:=\Phi_R(P_R)$.
  The fact that $z = f^n( \Phi_R^{-1}(w-n))$ is defined, and its value, are independent of the chosen non-negative integer $n$.
  We set $\Psi_R^{ \ext }(w) = z$.
  Otherwise $\Psi_R^{ \ext }$ is not defined at $w$.
\end{definition}

\begin{remark}\label{rk:1}
In the same way as the preceding paragraph, the extended repelling Fatou parametrization does not depend on the choice of repelling petal.
If $f$ is defined on all $\C$, the holomorphic map $\Psi_R^{ \ext }$ is defined on $\C$.
Else, $\Psi_R^{ \ext }$ is holomorphic on its domain of definition, which is open.
In all cases, two extended Fatou parametrizations differ only, in the strong sense, by pre-composition by a translation. 
\end{remark}

\begin{lemma}\label{lem:V}
  Every repelling quasi-petal $P'_R$ for $f$ takes the form $\Psi_R^\ext(V)$ where $V$ is an open subset of $\Dom (\Psi_R^\ext)$ and such that:
  \begin{itemize}
    \item $T_{-1}(V)\subset V$,
    \item $\Psi_R^\ext$ is injective on $V$,
    \item $(\Psi_R^\ext|_V)^{-1}$ is a Fatou coordinate for $P'_R$,
    \item $\piZ(V)=\CZ$,
    \item (in the case of a petal) $V$ is connected and simply
    connected.
  \end{itemize}
\end{lemma}
\begin{proof}
  Given a quasi-petal $P'_R$ and its Fatou coordinate $\Phi'_R$, the corresponding extended inverse repelling Fatou coordinate ${\Psi'_R}^\ext$ satisfies ${\Psi'_R}^\ext = \Psi_R^\ext\circ T_t$ in the strong sense for some $t\in \C$.
  The set $V=T_{t}\circ\Phi'_R(P_R)$ satisfies the required conditions.
%
\end{proof}

\begin{remark}
  With our definitions, the converse does not hold, i.e.\ if a set $V$ satisfies the above condition, $\Psi_R^\ext$ is not necessarily a repelling petal: in \Cref{ss:ex:4:1} we present a counter example with the Cauliflower map, contradicting our requirement that for a repelling petal, the conjugate of $T_{-1}:V\to V$ by $\Psi_R^\ext|_V$ must coincide in a neighborhood of $0$ with the local branch of $f^{-1}$.
\end{remark}

Because repelling petals are defined as attracting petals of inverse branches of $f$, and since $f$ may have several inverse branches, the intersection of two repelling petals is not necessarily a quasi-petal (see \Cref{d:pet}): \Cref{ss:ex:4:2} gives a counterexample.
What is true is that the intersection of any two repelling petals or quasi-petals always contains a repelling petal, but we will not use this fact.

\begin{corollary}\label{cor:invPsionPrep}
  Consider an extended repelling Fatou coordinate $\Psi_R^{\ext}$ for $f$, obtained by extending the repelling Fatou associated of some repelling petal $P_R$.
  Let $P'_R$ be another repelling petal.
  Then there exists a unique repelling Fatou coordinate $\Phi'_R$ on $P'_R$ that is also an inverse branch of $\Psi_R^\ext$, i.e.\ such that $\Psi_R^\ext \circ \Phi'_R = \id|_{P'_R}$.
\end{corollary}
\begin{proof}
  Existence is an immediate consequence of point of \Cref{lem:V}, more precisely the fact that $(\Psi_R^\ext)|_V)^{-1}$ is a Fatou coordinate for $P'_R$.
  
  Uniqueness: recall any two Fatou coordinates on $P'_R$ differ by a constant, call them $\Phi_R$ and $T_t\circ \Phi_R$ for some $t\in\C$.
  If both satisfy the hypothesis, i.e.\ $\Psi_R^\ext \circ \Phi'_R = \id|_{P'_R} = \Psi_R^\ext \circ T_t\circ \Phi'_R$ then we get $\Psi_R^\ext = \Psi_R^\ext \circ T_t$ on $Q_R:=\Phi_R(P_R)$.
  However, $\Psi_R^\ext$ is injective on $Q_R$ and $Q_R\cap T_t Q_R \neq \emptyset$
  (\footnote{Take any $w_0\in Q_R$; by $\piZ(Q_R)=\CZ$ we know that $\piZ(w_0-t)\in\piZ(Q_R)$, i.e.\ $\exists k_0\in \Z$ such that $w_0-t+k_0\in Q_R$; by $T_{-1}$-stability of $Q_R$ we have $w_0-n\in Q_R$ for all $n\geq0$ and $w_0-t-n\in Q_R$ for all $n\geq n_0=\max(0,-k_0)$; then $w_0-n_0\in Q_R\cap T_t Q_R$}),
  so $t=0$
  (\footnote{Take $w\in Q_R\cap T_t Q_R$ and denote $w'=w-t$. Then $w'\in Q_R$ and $\Psi_R^\ext(w') = \Psi_R^\ext(w)$, so by injectivity of $\Psi_R^\ext$ on $Q_R$, $w'=w$.}).
\end{proof}

We have $\Im(\Psi_R^{ \ext }) = \cup_{n \in \N} f^n(P_R)$, independent of the chosen petal $P_R$.
If $f$ is defined on all $\C$, the fact that $\Psi_R^{ \ext }$ is a nonconstant holomorphic function on $\C$ (since it is injective on the sets of the form $\Phi_R(P_R)$) implies that $\Im(\Psi_R^{ \ext  })$ avoids at most one point of $\C$. 

The equality
\[\Psi_R^\ext \circ T_1 = f\circ\Psi_R^\ext\]
holds in the strong sense.
In other words, $\forall w\in\C$, we have $w+1\in \Dom  \Psi_R^\ext$ $\iff$ ($w\in \Dom  \Psi_R^\ext$ and $ \Psi_R^\ext(w)\in \Dom  f$) and then $\Psi_R^\ext (w+1)= f(\Psi_R^\ext(w))$. 

\begin{proposition}\label{prop:sepals}
  Let $P_A$ be a very large attracting petal with a Fatou coordinate $\Phi_A:P_A\to \C$ and let $\Psi_R^\ext$ be any extended repelling Fatou parametrization.
  Then there exists $M'>0$ such that $\forall w\in\C$, $|\Im(w)|>M'$ implies $w\in\Dom \Psi_R^\ext(w)$ and $\Psi_R^\ext(w)\in P_A$.
\end{proposition} 
\begin{proof}
  The case of $\Im(w)<-M'$ is analogue to that of $\Im(w)>M'$ so we only treat the second.
  The set $Q_A=\Phi_A(P_A)$ contains a set of the form
  $\setof{w\in \C}{\Im(w)>M}$ for some $M>0$.
  Let us choose a repelling petal $P_R$ which is a $\pi/2$-petal and choose for repelling Fatou coordinates $\Phi_R$ on $P_R$, the one such that $\Psi_R^\ext\circ\Phi_R$ is the identity on $P_R$ (see \Cref{rk:1}).
  Let $S=\setof{w\in\C}{x_0\leq \Re(w)< x_0+1}$ be contained in $Q_R:=\Phi_R(P_R)$.
  By \Cref{p:petalesAR}, there exists $M'>0$ such that if $w\in S$ and $\Im(w)>M'$ then $z:=\Phi_R^{-1}(w)\in P_A$ and $\Im(\Phi_A(z))>M$.
  In particular $w\in\Dom \Psi_R^\ext$.
  By the $T_1$-invariance of $\Dom \Psi_R^\ext$ and the fact that $S$ is a fundamental domain for $T_1$, we get that $w\in \C$ and $\Im(w)>M'$ imply $w\in \Dom \Psi_R^\ext$.
  Let us prove the deceivingly obvious fact that $\Psi_R^\ext(w)\in P_A$.

  Consider any $w_0\in S$ with $\Im(w_0)>M'$ and let $w_n=w_0+n$ for $n\in \Z$.
  Let $z_0=\Psi_R^{\ext}(w_0)$.
  By the above paragraph, $z_0\in P_A$ and $\Im(\Phi_A(z_0)) > M$.
  \\
  \emph{For $n\geq 0$:} then we get by induction that $f^n(z)$ is defined and belongs to $P_A$ and
  $\Psi^\ext(w+n) = f^n(z)$ so $\Phi_A(\Psi_R^\ext(w+n))$ is defined and equals $\Phi_A(f^n(z))=\Phi_A(z)+n$.
  In particular its imaginary part is $>M$.
  \\
  \emph{For $n<0$:} we use an analytic continuation argument together with an induction.
  More precisely, let $H=\setof{w\in\C}{\Im w>M}$ and $W=\Phi_A^{-1}(H) \subset P_A$.
  The map $w\mapsto w-1$ on $H$ is conjugate by $\Phi_A^{-1}$ to a branch $g:W\to W$ of $f^{-1}$ that is a bijection of $W$, in particular $\forall w\in W$, $g\circ f(w)=w$.
  Denote $C = \setof{w\in\C}{\Im w>M'\text{ and }\Re w>x_0}$.
  We have $\forall w\in T_1 C$, $g(\Psi_R^\ext(w))=g(f(\Psi_R^\ext(w-1)))=\Psi_R^\ext(w-1)$.
  Note that $T_{-k}C \subset \Dom \Psi_R^\ext$.
  Our induction hypothesis $H_k$ for $k\geq 0$ is: $\Psi_R^\ext(T_{-k} C) \subset W$, which is satisfied for $k=0$.
  If $H_k$ holds then the domain of $g\circ \Psi_R^\ext$ contains 
  $T_{-k} C$.
  The identity $\Psi_R^\ext(w-1) = g(\Psi_R^\ext(w))$, valid on $T_1 C$ as we saw above, then holds on $T_{-k}C$ by analytic continuation.
  As a consequence, $\Psi_R^\ext(T_{-k-1}C) = g(\Psi_R^\ext(T_{-k}C)) \subset g(W)=W$.
\end{proof}

\subsection{Horn maps}\label{sec:horn_maps}

In this section, as in the previous ones, we remind the elements of a standard theory, with a presentation which might be original in some places.

We denote
\[\piZ: \C \to \CZ\]
the canonical projection $z\mapsto z+\Z$.

\begin{PropositionDefinition} \label{d:cornes}
  The lifted horn map is defined as \[ h = \Phi^\ext_A \circ \Psi^\ext_R \]
  Its domain
  \[\wDd = \Dom  h = (\Psi_R^\ext)^{-1}(B_f)\]
  is invariant by $T_1$:
  \[ T_1(\widetilde\Dd) = \widetilde\Dd,\]
  and 
  \[ h\circ T_1 = T_1\circ h .\]
  There is thus a quotient map $\h$ called the (non-lifted) horn map, satisfying
  \[ \xymatrix{
  \widetilde\Dd\ar[d]_{\piZ} \ar[r]^{h} & \ar[d]^{\piZ} \C \\
  \mathcal D\ar[r]_{\h} & \CZ
  } \]
  where
  \[\mathcal D=\Dom  \h = \piZ(\wDd) =  \wDd/\Z.\]
\end{PropositionDefinition}
\begin{proof}
  We saw that $\Psi_R^\ext \circ T_1 = f\circ\Psi_R^\ext$ holds with equality of domains.
  This is also the case for the relation $\Phi_A^\ext \circ f = T_1\circ \Phi_A^\ext$ since $f^{-1}(B_f) = B_f$.
  Each relation in the following sequence is then valid with equality of domains:
  $(\Phi_A^\ext \circ \Psi_R^\ext) \circ T_1
  = \Phi_A^\ext \circ (\Psi_R^\ext \circ T_1)
  = \Phi_A^\ext \circ (T_1 \circ \Psi_R^\ext)
  = (\Phi_A^\ext \circ T_1) \circ \Psi_R^\ext
  = (T_1 \circ \Phi_A^\ext ) \circ \Psi_R^\ext
  = T_1 \circ (\Phi_A^\ext \circ \Psi_R^\ext)
  $.
\end{proof}

Though we will not do it here, this is is best interpreted in terms of quotients and grand orbits, see \cite{t:LM}.

Note that since $\Phi_A^\ext$ and $\Psi_R^\ext$ do not depend on the choice of petals $P_A$ and $P_R$, but only respectively post and pre composition by translations, the map $h$ and its domain $\wDd$ only depend on these translations.
The horn map is said to be \emph{normalized} if both attracting and repelling Fatou coordinates are normalized.

\begin{remark}
The set $\wDd$ contains $\Phi_R (P_A)=\Phi_R (P_A \cap P_R)$. Let $D_R$ be a fundamental domain of $f$ restricted to a repelling petal $P_R$.
Then $\Phi_R^\ext(D_R)$ is a fundamental domain for the action of $T_1$ on $\C$ and the set $\Dd$ is exactly equal to 
$\piZ\circ \Phi_R( B_f \cap D_R)$.
\end{remark}

We denote $\mathbb{H} = \setof{z\in \C}{\Im(z) > 0}$ the Poincaré half-plane. For $\lambda \in \C^*$, the notation $\lambda \mathbb{H}$ refers to the image of $\Hh$ by the map $z \mapsto \lambda z$.

By \Cref{p:petalesAR}:
\begin{corollary}\label{cor:hDomCHP}
   The domain of $\h$ contains a punctured neighborhood 
   of $+ i \infty$, and a punctured neighborhood 
   of $- i \infty$.
\end{corollary}

\begin{definition}%
  We denote $\Dd^+$ the connected component of $\Dd$ containing a neighborhood of $+ i \infty$ and $\Dd^-$ the one containing a neighborhood of $- i \infty$.
  We denote $\wDd^\pm = \piZ^{-1}(\Dd^\pm)$.
\end{definition}

\begin{remark}
  There are maps $f$ for which $\Dd^+ = \Dd^-$, for instance any Blaschke product, seen as a mapping $\widehat{\C}\to\widehat{\C}$, that has a parabolic fixed point with only one attracting direction (which must then be tangent to the unit circle $\partial \D$, and $J(f)$ will be a disconnected subset of $\partial \D$).
  On the other hand, if $\Dd^+ \neq \Dd^-$, then $\wDd^+$ and $\wDd^-$ are disjoint.
\end{remark}

Recall that we denote $B^0_f$ the immediate basin of $f$ at $0$.

\begin{proposition} \label{p:tildeDpm}
$\Psi_R^{ \ext }(\wDd^\pm) \subset B^0_f$
\end{proposition}

\begin{proof}
First, $\Psi_R^{ \ext  }(\wDd^\pm) \subset \Psi_R^{ \ext  }(\wDd) \subset\Dom  \Phi^\ext_{A} = B_f$.
The set $\Psi_R^{ \ext  }(\wDd^\pm)$ is a continuous image of a connected set so is connected, so it is contained in a connected component of $B_f$, so we just need to prove that it contains at least one point of $B_f^0$.
Let $P_R$ be a repelling $\pi/2$-petals. Its image $Q_R=\Phi_R(P_R)$ is a left half plane and we have $P_R = \Psi_R^{ \ext  }(Q_R)$.
Let $P_A$ be a large attracting petal.
By \Cref{p:petalesAR}, $\Phi_R(P_R \cap P_A)$ contains points of arbitrary high/low imaginary part, hence it contains points $w$ in $\wDd^\pm$, since the latter contains an upper/lower half plane.
Since $\Psi_R^\ext(w)\in P_A\subset B^0_f$, the lemma is proved.
\end{proof}

\begin{remark}
  It follows that $\wDd^\pm$ can also be defined as a connected component of the set $\wDd^0 = \piZ(\Dd^0)$ where $\Dd^0 = (\Psi_R^\ext)^{-1}(B_f^0)$.
  However, we will not use that notation in the sequel.
  Note that $T_1(\Dd^0)\subset \Dd^0$ but that the converse inclusion is not necessarily true, see the example in \Cref{sub:ex:2}.
\end{remark}

\begin{proposition}\label{p:hi}
The map $\h$ satisfies
\[ \lim_{w \to \pm i \infty} \h(w) = \pm i \infty \]
It thus has erasable singularities at $\pm i \infty$, fixing each of these two points.
\end{proposition}
\begin{proof}
  In order to apply \Cref{p:petalesAR}, more precisely the version with $A$ and $R$ permuted, we choose as in that statement $P_A$ and $P_R$ to be respectively an attracting $\beta>\pi/2$ petal and a repelling $\pi/2$ petal, and $D' = \setof{z\in \C}{ x_0 \leq \Re(z) < x_0 + 1 }$ contained in $Q_R:=\Phi_R(P_R)$.
  For a given $M>0$, it gives us an $M'>0$ such that $\forall w\in D'$, if $\Im w>M'$ then $z:=\Phi_R^{-1}(w)\in P_A$ and $\Im \Phi_A(z)> M$, i.e.\ $\Im(h(w)) >M$ since in this case, $h(w) = \Phi_A(\Psi_R^{-1}(w))$.
  Consider now any $\mathbf{w}\in\CZ$ with $\Im(\mathbf{w})>M'$.
  Let $w\in D$ be its unique representative modulo $\Z$.
  Then $\h(\mathbf{w}) = \piZ (h(z))$ and by the above, $\Im h(w)>M$.
  
  The analysis is similar near the bottom end of the cylinder
\end{proof}

\begin{remark}\label{r:hlim2}
  The fact that $\h$ has a continuous extension fixing both ends of the cylinder $\CZ$ can be rephrased as follows:
  \begin{itemize}
  \item For all $M>0$ there exitst $M'>0$ such that $\forall w\in\C$, if $\Im w>M$ then $w\in\Dom \Psi_R^{\ext}$ and $z:=\Psi_R^{\ext}(w) \in B_f$ and $\Im(\Phi_A^\ext(z))>M'$.
  \item For all $M>0$ there exitst $M'>0$ such that $\forall w\in\C$, if $\Im w<-M$ then $w\in\Dom \Psi_R^{\ext}$ and $z:=\Psi_R^{\ext}(w) \in B_f$ and $\Im(\Phi_A^\ext(z))<-M'$.
  \end{itemize}
\end{remark}

\begin{corollary} \label{cor:dlCornes}
1. Let $P_A, P_R$ be attracting, repelling petals of $f$ included in a neighborhood $W$ of $0$ where $f$ admits a local inverse branch $f^{-1}$.
We denote $\Phi_A, \Phi_R$ Fatou coordinates on these petals.
Then for $z \in P_A \cap P_R$, $\Im (\Phi_A(z) )$ is high (resp.\ low) \underline{if and only if} $\Im (\Phi_R(z) )$ is high (resp.\ low).
 
2. Nevertheless for $z \in B_f \cap P_R$ we have in general only one implication: \underline{if} $\Im( \Phi_R(z))$ is high (resp.\ low) \underline{then} $\Im (\Phi_A^{ \ext  }(z) )$ is high (resp.\ low), where $\Phi_A^{ \ext  }$ is the extended Fatou coordinate.
\end{corollary}

\begin{proof}
To show the equivalence in 1., we use the horn maps $h, h'$ of $f$ and $f^{-1}$. The formulas, when $z \in P_A \cap P_R$:
\[ \Phi_A(z)= h( \Phi_R(z) ) \]
\[ - \Phi_R(z) = h'( - \Phi_A(z)) \]
give the claim by applying \Cref{p:hi} to $h$ and $h'$.

To show the implication in 2., we use the fact that the first formula above still holds when $z \in B_f \cap P_R$, by replacing $\Phi_A$ by $\Phi_A^{ \ext  }$:
\[ \Phi_A^\ext(z) = h( \Phi_R(z)) \]
and again we apply \Cref{p:hi} to $h$.

Note nonetheless that the second formula analogue, $- \Phi^\ext_R(z) = h'( - \Phi_A(z))$, would only be avalaible for $z \in B_{f^{-1}} \cap P_A$ where $P_A$ is an attracting petal of $f$ (and consequently a repelling petal of $f^{-1}$), and this is not enough to show the converse statement.

To give an effective counter-example to the converse statement in 2., we choose the cauliflower map $f: z \mapsto z + z^2$.
A small enough attracting petal $P_A$ admits iterated inverse images $P'$ by $f$ which accumulate densely over the boundary of $B_f$, with Euclidean diameters converging to $0$. Since a repelling petal $P_R$ contains in its interior boundary elements of $B_f$, it must contain the closure of one of the $P'$ totally. That is why $\Im(\Phi_R)$ is bounded on $P'$, although it is not the case for $\Im(\Phi_A)$.
Incidentally it implies that the converse statement of \Cref{p:hi} is also false : for $z\in \Dd$, if $\h(w)$ has a high resp.\ low imaginary part, it is not necessary that $w$ also has a high resp.\ low imaginary part. 
\end{proof}

%
%
%

We end this section with the proof of \Cref{lem:i2} in the introduction, whose statement we copy here:

\begin{lemma}\label{lem:i2:copy}
  Any holomorphic map $\psi^\pm:\wDd^{\pm}_1\to \wDd^{\pm}_2$ such that $\h_1^\pm = T_\sigma^{-1} \circ \h_2^\pm \circ \psi^\pm$ extends holomorphically into a map fixing $\pm i\infty$, and
  \[\psi^\pm = w + \rho^\pm + o(1)\]
  as $w\to \pm i\infty$, for some $\rho^\pm\in\C$.
\end{lemma}
\begin{proof}
  Consider the isomorphism $E:\CZ\to\C^*$ defined by $E(z)=2\pi i z$.
  Denote $D_i = E(\Dd_i^+)$, which are two punctured open neighborhoods of $0$ contained in $\C^*$.
  Let
  \[\zeta(z)=E\circ\psi^+\circ E^{-1}:D_1\to D_2.\]
  The relation $\h_1^+ = T_\sigma^{-1} \circ \h_2^+ \circ \psi^+$ translates into $g_2 \circ \zeta = g_1$
  holds on $D_1$, where $g_1 = E\circ \h_1\circ E^{-1}$ and $g_w = E\circ T_{-\tilde\sigma}\circ \h_2\circ E^{-1}$ are holomorphic functions on $E(\Dd_1)$ and $E(\Dd_2)$, which both extend at $z=0$ into holomorphic maps fixing $0$ and $0$ is not a critical point of $g_1$ nor $g_2$, and the same at $z=\infty$.
  The isolated singularity at $0$ of $\zeta$ cannot be essential, for otherwise $g_2\circ \zeta$ would not extend continuously at $0$: indeed by Picard's theorem, $\zeta$ woud map any pointed neighbourhood of $0$ to a subset of $\C$ that avoids at most $2$ points, so $g_2\circ \zeta$ would fail to have a limit $0$, but $g_1$ has.
  We distinguish three cases.
  
  If $\zeta$ has a pole at $0$ then\footnote{This means that $D_2$ is a neighborhood of $\infty$, i.e.\ that $\Dd_2^+$ is a neighborhood of $-i\infty$, so $\Dd_2^+=\Dd_2^-$. This actually happens. Take a rational map having a parabolic fixed point at $0$ an whose Julia set is a Cantor set.
  Conjugate it by a homography fixing $0$ and so that $\infty$ in the Cantor set. Restrict it to the preimage of $\C$, so that we are in the setting of our statements. Then $\Dd$ is the complement of a Cantor subset of $\CZ$, in particular $\Dd=\Dd_2^+=\Dd_2^-$.}
  $\zeta(z)\sim a z^{-d}$ for some $d\in\N$ with $d>0$ and some $a\in\C^*$.
  Since $g_2(z)\sim bz$ as $z\to\infty$ for some $b\in\C^*$ and $g_1(z)\sim cz$ as $z\to 0$ for some $c\in\C^*$, the relation $g_2 \circ \zeta = g_1$ implies $baz^{-d}\sim cz$ as $z\to 0$, which is obviously wrong.
  
  If $\zeta(z)$ has a non-zero limit $a$ at $0$, then the relation 
  $g_2 \circ \zeta = g_1$ implies, since $\zeta$ cannot be constant near $0$ (because $g_1$ is not and $g_1=g_2\circ \zeta$ near $0$), that the domain of $g_2$ contains a punctured disk around $a$ and $g_2$ tends to $0$ at $a$.
  Let us prove this leads to a contradiction.
  Reinterpreting in terms of $\h_i$, the domain of $\h_2$ contains a punctured disk around some $w_0 \in \Dd_2$ and $\h_2$ tends to $+i\infty$ at $w_0$.
  Passing to the non-lifted horn maps, the domain of $h_2$ contains a punctured disk around some $\tilde w_0$ such that $\Im h_2(\tilde w)$ tends to $+\infty$ as $\tilde w$ tends to $\tilde w_0$.
  Such a behavior is impossible for a holomorphic function: consider a first circle $C$ centered on $\tilde w_0$ and consider a second one $C'$ small enough so that $y':=\min_{z\in C'} \Im h_2(z) > y:=\max_{z\in C} \Im h_2(z)$.
  Let $A$ be the annulus between $C$ and $C'$.
  Since $h_2(\ov{A})$ is connected, there must be some $\tilde w\in A$ such that $\Im h_2(z)\in (y,y')$.
  This contradicts the winding number theorem: the winding numbers of $h_2(C)$ and $h_2(C')$ around $w$ are equal to $0$, while their difference must be equal to the number of preimages of $w$ in $A$ by $h_2$, so not $0$.
   
  If $\zeta$ has a zero at $0$ then let $d$ be its order: $\zeta(z)\sim az^d$ for some $a\in\C^*$.
  Since $g_2(z)\sim bz$ as $z\to 0$ for some $b\in\C^*$ and $g_1(z)\sim cz$ as $z\to 0$ for some $c\in\C^*$, the relation $g_2 \circ \zeta = g_1$ implies $baz^d\sim cz$ as $z\to 0$, which is only possible if $d=0$.
  It follows that $\zeta(z)\sim az$ as $z\to 0$ and this implies the estimate on $\psi^+$.
  
  The proof is the same for $\psi^-$.
\end{proof}

\section{Definition of local semi/pseudo-conjugacy on immediate basins and immediate consequences}\label{sec:defs_csq}

We start this section by proving some properties about domains of horn maps and the sets we denoted $U^0$ in the introduction.

\subsection{About domains of horn maps and $U^0$}\label{sub:dhu0}

As in the previous section, we assume that $f: \Dom (f)\subset\C\to \C$ is holomorphic with $f(z)=z+az^2+\ldots$ and $a\in\C^*$.
We recall that $B_f$ denotes the basin of its parabolic point at the origin and $B_f^0$ the immediate basin,
$h$ and $\h$ denote its lifted and non-lifted horn maps (see \Cref{d:cornes}), $\Dd^\pm$ is the connected component of the domain $\Dd$ of $\h$ containing a neighborhood of $\pm i \infty$ and that $\wDd^\pm = (\Psi_R^{\ext})^{-1} (B_f) = \piZ^{-1}(\Dd^\pm)$, where $\piZ: \C \to \CZ$ is the canonical projection.

\begin{definition}%
Given any open subset $U$ of $\C$ containing $0$, we let $U^0$ denote the connected component of $U\cap B_f$ (equivalently, of $U\cap B_f^0$) that contains a germ of the attracting axis of $f$.
\end{definition}

The two versions in the definition are equivalent because any connected subset of $B_f$ containing a germ of the attracting axis is necessarily contained in $B_f^0$.

\begin{remark} \label{r:petitPetAttr}
  The attracting petals are connected and contain a germ of the attracting axis. Thus the set $U^0$ might be understood as the connected component of $B_f \cap U$ containing a given attracting petal of $f$ included in $U$, or equivalently all the attracting petals that are included in $U$.
\end{remark}

\begin{remark}\label{rk:incU0}
  The following is straightforward:
  \[U'\subset U \implies {U'}^0\subset U^0.\]
\end{remark}

To motivate the introduction of the set $U^0$, let us give some of its properties in connection with domains of horn maps.

\subsubsection{A general inclusion}\label{ss:gi}

The title of this section refers to \Cref{p:petitPetRep}.
For this we have to introduce two new sets:

\begin{definition}\label{def:wtus}%
Given a repelling petal $P_R$, let $\Wt{U}$ denote the set of points of $w\in\C$ such that for all $n$ big enough, $w-n\in \Phi_R(U^0)$.
Let $\W{U} = \piZ(\Wt{U})$.
\end{definition}

We have $T_1(\Wt{U})=\Wt{U}$ so $\Wt{U}=\piZ^{-1}(\W{U})$.
The following inclusion is immediate:
\begin{equation}\label{eq:wsu0}
  \W{U}\subset\piZ(\Phi_R(U^0))
\end{equation}
but the converse inclusion $\piZ(\Phi_R(U^0)) \subset \W{U}$ does not necessarily hold, as is shown by the example of \Cref{sub:ex:2}.
On the other hand, in the particular (not uncommon) cases where $\Phi_R(U^0)$ is stable under $T_{-1}$, then $\W{U}=\piZ(\Phi_R(U^0))$.

\begin{lemma}\label{lem:WindepU}
  The sets $\Wt{U}$ and $\W{U}$ are independent of the choice of repelling petal $P_R$.
\end{lemma}
\begin{proof}
  Given two repelling petals, according to \Cref{lem:V}, they can be written as $P_R=\Psi_R^\ext(V)$ and $P'_R=\Psi_R^\ext(V')$ where $V$ and $V'$ are two  $T_{-1}$-stable open sets on which $\Psi_R$ is injective and such that $\piZ(V)=\CZ$ and $\piZ(V')=\CZ$.
  The map $\Phi_R$ is the inverse of $\Psi_R^\ext|_V$ and $\Phi_R'$ is the inverse of $\Psi_R^\ext|_{V'}$.
  Let $W$ and $W'$ denote the respective sets $\Wt{U}$ obtained from $P_R$ and $P'_R$.
  Consider $w\in W$: ($\exists N\in \N$) ($\forall n\geq N$) $w-n\in\Phi_R(U^0)$. Note that $\Phi_R(U^0) = V\cap (\Psi_R^\ext)^{-1}(U^0)$.
  Since $\piZ(V')=\CZ$ and by $T_{-1}$-stability of $V'$, there is some $N'\geq N$ such that $\forall n\geq N'$, $w-n\in V'$.
  Since $w-n$ is also in $(\Psi_R^\ext)^{-1}(U^0)$ we get $w-n\in V'\cap (\Psi_R^\ext)^{-1}(U^0) = \Phi'_R(U^0)$.
  Hence $W\subset W'$. By a symmetric argument $W'\subset W$.
\end{proof}

The set $\W{U}$ can be written as $\W{U}=\piZ(\bigcap_{n\in\N}T_n \circ\Phi_R(U^0))$, so is a intersection of a decreasing sequence of open sets.
Such an intersection is not necessarily open but in the case of $\W{U}$, equivalently of $\Wt{U}$, this is the case:
\begin{proposition}
  $\Wt{U}$ is open
\end{proposition}
\begin{proof}
  Let $w\in\Wt{U}$, and $N\in\N$ such that $\forall n\geq N$, $w-n\in\Phi_R(U^0)$.
  The open set $\Phi_R(U^0)$ contains $w-N$ hence contains $B-N$ where $B:=B(w,\eps)$ for some $\eps>0$.
  The set $\Phi_R^{-1}(B-N)$ is contained in $U^0$ hence in $B_f$.
  From $B-N\subset Q_R:=\Phi_R(P_R)$ we get $\forall n\geq N$, $B-n\subset Q_R$.
  The set $\Phi_R^{-1}(B-n)$ tends uniformly to $0$ (it follows for instance from \Cref{lem:asymp2}) so is contained in $U$ for $n$ big enough.
  It contains a point in $U^0$: the point $\Phi_R^{-1}(w-n)$.
  As a consequence: $\Phi_R^{-1}(B-n)\subset U^0$ (still for $n$ big enough).
  It follows that $B\subset \Wt{U}$.
\end{proof}

\begin{proposition} \label{p:petitPetRep}
Let $P_R$ be a repelling petal of $f$ and $\Phi_R$ the associated Fatou coordinate.
Let $U$ be an open neighbourhood of $0$.
Then
\[ \Dd^+ \cup \Dd^- \subset \W{U} \] 
\end{proposition}

\begin{proof}
Let $P_A$ be a very large attracting petal of $f$ included in $U$, hence in $U^0$.
A consequence of \Cref{prop:sepals} is that there is some $M>0$ such that $P_A$ contains $\Phi_R^{-1}(H\cup B)$, where $H, B$ are the upper, lower half planes such that $H \cup B = \setof{w \in \C}{| \Im (w)| > M }$.

Consider $w \in \Dd^\pm$ and a path $\gamma : [0,1]\to\Dd$ joining $w$ with an element $w_0 \in \Dd^\pm$ such that $|\Im (w_0)| > M$.
We lift this path by $\piZ$ to obtain a path $\widetilde{\gamma}$ of $\wDd^\pm$ joining $\widetilde{w}$ to $\widetilde{w}_0 \in \C$ such that $\pi(\widetilde{w}) = w$ and $\pi(\widetilde{w}_0) = w_0$.
The set $\widetilde{\gamma}$ is compact, so by \Cref{p:covPetCompact}, which is also valid for repelling petals, for all $n \in \N$ big enough, $\widetilde{\gamma} - n$ is included in $Q_R:=\Phi_R(P_R)$.
Then $\delta_n := \Phi_R^{-1}\circ T_{-n}\circ \wt \gamma = \Psi_R^\ext\circ T_{-n}\circ \wt \gamma$ is defined on $[0,1]$.
Moreover as $n$ tends to infinity, $\delta_n$ tends uniformly to $0$.
So for $n$ big enough (the support of) $\delta_n$ is contained in $U$.
Since $\wt\gamma \subset\wDd$ and $T_{-1}(\wDd) = \wDd$, we have $\wt\gamma-n \subset \wDd$, and since $\wDd = (\Psi_R^{\ext})^{-1}(B_f)$, we get $\delta_n \subset B_f$.
Also, $\delta_n(1)= \Phi_R^{-1}(\wt w_0-n)$, which by the first paragraph belongs to $P_A$, hence to $U^0$.
Since $U^0$ is the connected component of $U\cap B_f$ that contains $P_A$, we have $\delta_n\subset U^0$.

So $\delta_n(0) = \Phi_R^{-1}(\widetilde{w} - n) \in U^0$, and by construction $\delta_n(0)\in P_R$.
Hence $w \in \piZ \circ \Phi_R(U^0)$.
But better: the above tells us that for all $n$ big enough, say $n\geq N$, $\wt w \in T_n(\Phi_R(U^0))$.
So $\wt w-N\in \Wt{U}$.
\end{proof}

The converse inclusion $\W{U}\subset \Dd^+\cup D^-$ is not true in full generality: a counterexample is given in \Cref{sub:ex:1}.

\subsubsection{A converse inclusion in favorable cases}\label{ssub:ci}

If $U$ is small enough, we have the opposite inclusion to \Cref{p:petitPetRep}.
More precisely, let $\alpha,\beta \in (0 , \pi )$ such that
\[ \alpha+\beta>\pi \]
Let $P_A$ be an attracting $\alpha$-petal and
$P_R$ be a repelling $\beta$-petal
of $f$ such that $P_A \cap P_R$ has two connected components $C_+, C_-$. This is possible by the Leau-Fatou flower theorem.
In this case, we must have (permuting $C_+$ and $C_-$ if necessary):
\begin{lemma}\label{lem:phiCinD}
\[ \piZ\circ\Phi_R(C_\pm) \subset \Dd^\pm \]
\end{lemma}
\begin{proof}
  Indeed, $\piZ\circ\Phi_R(C_\pm)$ is connected, included in $\Dd$ and intersects $\Dd^\pm$.
  The first claim is immediate.
  For the second claim, use $\Dd = \piZ(\wDd)$ and $\wDd = \Dom \Phi_A^{\ext}\circ\Psi_R^\ext \supset \Dom \Phi_A\circ(\Phi_R)^{-1} \supset \Phi_R(P_A\cap P_R)$.
  The last claim is proved as follows: by \Cref{cor:sec1}, $P_A\cap P_R$ contains a small sector $V$ bisected by a half-line transversal to the attracting and repelling axis, so that $\Im -1/az$ tends to $\pm\infty$ when $z$ tends to $0$ within $V$.
  Using \Cref{lem:asymp1} we get that for such $z$, $\Im \Phi_R(z)$ tends to $\pm\infty$.
  Since $\Dd^\pm$ contains an upper/lower half-plane, the claim follows.
\end{proof}

In other words, under these conditions: 
\[ \piZ\circ\Phi_R(P_A) = \piZ\circ\Phi_R(P_A\cap P_R) \subset \Dd^+ \cup \Dd^- \]
We now assume that $U$ is an open subset of $\C$ containing $0$ and such that
\[U \subset P_A \cup P_R\cup \{0\}\]
(still under the same conditions on the petals $P_A$ and $P_R$).

\begin{proposition}\label{p:petitPetRep2} 
Under these asumptions
\[ \Phi_R(U^0) = \Phi_R(U^0\cap P_R)\subset \wDd^+ \cup \wDd^-\]
Together with \Cref{p:petitPetRep} and \cref{eq:wsu0} this implies
\[ \W{U} = \piZ\circ\Phi_R(U^0) = \Dd^+ \cup \Dd^-\]
\end{proposition}

\begin{proof} 
  It is sufficent to show the inclusion when $U = P_A \cup P_R \cup \{0\}$, since $U' \subset U$ implies ${U'}^0 \subset U^0$.
  Since $U = P_A \cup P_R \cup \{0\}$, we have $U \cap B_{f} = P_A \cup (P_R \cap B_{f})$.
  The set $U^0$ is the connected component of $U \cap B_{f}$ containing $P_A$.
  It is then the union of $P_A$ with the connected components $C$ of $P_R \cap B_{f}$ intersecting $P_A$.
  We already know that $\Phi_R(P_A)=\Phi_R(P_A\cap P_R)\subset \wDd^+ \cup \wDd^-$.
  A component $C$ as above intersects $P_A \cap P_R = C_+ \cup C_-$.
  The nonempty set $(C_+ \cup C_-) \cap C$ is included in the domain of $\Phi_R$, this implies that the set $\Phi_R(C)$, which is connected and included in $\Dd$, intersects 
  $\Phi_R( C_+ \cup C_-) \subset \wDd^+ \cup \wDd^-$.
  Thus $\Phi_R(C) \subset \wDd^+ \cup \wDd^-$.
\end{proof}

\begin{figure}[!t]
  \begin{center}
    \begin{tikzpicture}
      \node at (0,0) {\includegraphics[scale=0.7]{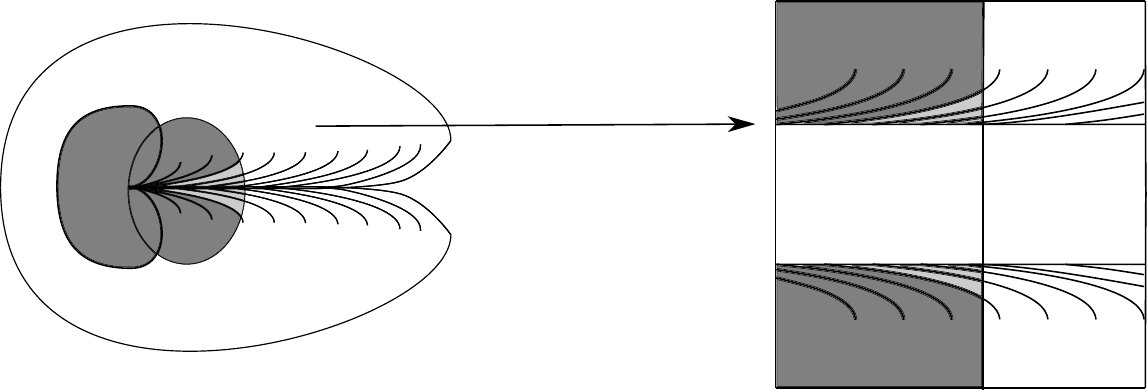}};
      \node at (-4.8,1.4) {$U^0$};
      \node at (-4,-1) {$P_R$};
      \node at (-6,-1) {$P_A$};
      \node at (0,1.2) {$\Phi_R$};
      \node at (4,0) {$Q_R$};
    \end{tikzpicture}
  \end{center} 

  \caption{Illustration of \Cref{p:petitPetRep2}}
  \label{f:U10}
  \par\bigskip
  \footnotesize{\textit{Case $U=P_A\cup P_R\cup \{0\}$.
The attracting petal is represented by a cardioid, and the repelling petal by a circle.
The set $Q_R:=\Phi_R(P_R)$ is represented by a left half plane. The dark gray region at the left represents $U^0$, while the union of the dark gray and the light gray is equal to $(U\cap B^0_f) - U^0$.
The gray regions on the right represents the images by $\Phi_R$ of the gray zones on the left.}}
\end{figure}

See \Cref{f:U10}.

\subsection{Semi-conjugacy from petals}\label{sub:scfp}

Let $f_i : \Dom  f_i \to \C$, $i=1,2$ be two holomorphic functions, where $\Dom  f_i \subset \C$ is an open neighborhood of $0$, with simple parabolic points at $0$, i.e.\  
\[ f_i(z) = z + a_i z^2 + o(z^2) \]
with $a_i\neq 0$.

\begin{hypothesis}\label{h:1} We assume in this section that $P_A^1$ is a petal of $f$ and that \[\phi: P_A^1\to B_{f_2}\]
is a semi-conjugacy from $f_1$ to $f_2$, i.e.\ it satisfies $\forall z\in P_A^1$, $\phi(f_1(z)) = f_2(\phi(z))$.
\end{hypothesis}

The first remark is that under this hypothesis, $\phi$ actually maps in the immediate basin $B^0_{f_2}$:
\[\phi(P_A^1)\subset B^0_{f_2}\]
Indeed, $\phi(P_A)$ is a connected subset of $B_{f_2}$ and has non-empty intersection with an attracting petal $P_A^2$ of $f_2$: pick any $z\in \phi(P_A^1)$; then $f_2^n(z)$ eventually enters $P_A^2$ while staying in $\phi(P_A^1)$.

Consider as in \Cref{sec:cyl} the quotient $B_f/f$ of the basin by the equivalence relation generated by $f(z)\sim z$ and denote the canonical projection:
\[ \pi_f : B_f \to B_f/f\]

\begin{PropositionDefinition}[Induced map]
  The map $\phi: P_A^1\to B^0_{f_2}$ descends to a map $\ov\phi:B_{f_1}/{f_1} \to B_{f_2}/{f_2}$, i.e.\ the following diagram commutes.
  \[
  \xymatrix{
    P_A^1 \ar[r]^{\phi} \ar[d]_{\pi_{f_1}} & B^0_{f_2} \ar[d]^{\pi_{f_2}}
    \\
    B_{f_1}/{f_1} \ar[r]_{\ov\phi} & B_{f_2}/{f_2}
  }
  \]
\end{PropositionDefinition}
\begin{proof}  
  The relation $\phi\circ f_1 = f_2\circ \phi$ is valid on $P_A^1$ and $f(P_A^1)\subset P_A^1$.
  This implies $\forall z\in P_A^1$ and $\forall n\in\N$, $\phi\circ f_1^n(z) = f_2^n\circ \phi(z)$.
  If $z,z'\in P_A^1$ are $f_1$-equivalent, i.e.\ $\exists n,m\in\N$ such that $f_1^n(z)=f_1^m(z')$ then by the previous sentence, $f_2^n\circ \phi(z) = f_2^m\circ \phi(z')$, so $\phi(z)$ and $\phi(z')$ are $f_2$-equivalent.
  This defines a quotient map $\ov \phi$ on the set $\pi_{f_1}(P_A^1)$.
  This set is equal to $B_{f_1}/f_1$
  since every orbit in $B_{f_1}$ eventually enters $P_A^1$.
\end{proof}

We recall (see \Cref{cor:quot}) that extended attracting Fatou coordinate $\Phi_A^\ext : B_f\to \C$ induces a bijection from $B_f/f$ to $\CZ$.
We decide to denote this bijection $\ov\Phi_A^\ext$:
\[
\xymatrix{
  B_f \ar[d]_{\pi_f} \ar[r]^{\Phi_A^\ext} & \C \ar[d]^{\piZ}
  \\
  B_{f}/{f} \ar[r]_{\ov\Phi_A^\ext} & \CZ
}
\]
commutes.
We will use the bijection $\ov\Phi_A^\ext$ in order to work with $\CZ$ instead of $B_f/f$.
In particular, we consider the mapping $\ov\phi$ from the point of view of $\CZ$ by letting $\wh\phi = \ov\Phi_A^{2,\ext} \circ \ov\phi \circ (\ov\Phi_A^{1,\ext})^{-1}$, so 
\[
\xymatrix{
  B_{f_1}/{f_1} \ar[r]^{\ov\phi} \ar[d]_{\ov\Phi_A^{1,\ext}} & B_{f_2}/{f_2} \ar[d]^{\ov\Phi_A^{2,\ext}}
  \\
  \CZ \ar[r]_{\wh\phi} & \CZ
}
\]
commutes.
Consider the last three diagrams.
We get from the first and the third (and rotating the diagrams clockwise), that
\[
\xymatrix{
  P_A^1 \ar[d]_{\phi} \ar[r]^{\pi_{f_1}} & B_{f_1}/f_1 \ar[d]_{\ov\phi} \ar[r]^{\Phi_A^{1,\ext}} & \CZ \ar[d]^{\wh\phi}
  \\
  B^0_{f_2} \ar[r]^{\pi_{f_1}} & B_{f_2}/f_2 \ar[r]^{\Phi_A^{2,\ext}} & \CZ
}
\]
commutes; then using the second, that
\[
\xymatrix{
  P_A^1 \ar[d]_{\phi} \ar[r]^{\Phi_A^{1,\ext}} & \C \ar@{.>}[d]_{\text{?}} \ar[r]^{\pi_{\Z}} & \CZ \ar[d]^{\wh\phi}
  \\
  B^0_{f_2} \ar[r]_{\Phi_A^{2,\ext}} & \C \ar[r]_{\pi_{\Z}} & \CZ
}
\]
commutes and it is natural to ask whether there is a middle arrow closing the diagram (indicated by a question mark).
The answer is yes:

\begin{PropositionDefinition}\label{prop:inter}
  There exists a unique map $\wt\phi:\C\to\C$ such that $\wt\phi\circ T_1 = T_1 \circ \wt \phi$ and such that the following diagram commutes:
  \[
  \xymatrix{
    P_A^1 \ar[d]_{\phi} \ar[r]^{\Phi_A^{1,\ext}} & \C \ar[d]_{\wt\phi} \ar[r]^{\pi_{\Z}} & \CZ \ar[d]^{\wh\phi}
    \\
    B^0_{f_2} \ar[r]_{\Phi_A^{2,\ext}} & \C \ar[r]_{\pi_{\Z}} & \CZ
  }
  \]
  This map $\wt\phi$ holomorphic.
\end{PropositionDefinition}
\begin{proof}
  We first define $\wt\phi$ on the image of $P_A^1$ by $\Phi_A^{1,\ext}$ by noticing that $\Phi_A^{1,\ext}$ is injective on $P_A^1$ because it coincides on it with a Fatou coordinate $\Phi_A^1$ for the petal, and we can define $\wt\phi = \Phi_A^{2,\ext} \circ \phi \circ (\Phi_A^{1})^{-1}$ and the left square commutes.
  
  But note that we only defined $\wt\phi$ on $Q_A^1:=\Phi_A^{1}(P_A^1)$, which is not necessarily all $\C$.
  However $T_1(Q_A^1)\subset Q_A^1$ and we now check that
  \begin{equation}\label{eq:2}
    \wt\phi(w+1) = \wt\phi(w)+1\text{ for all }w\in Q_A^1,
  \end{equation}
  Indeed, since $w=\Phi_A^{1}(z)$ for some $z\in P_A^1$, we have $w+1=\Phi_A^{1}(f_1(z))$ so 
  $\wt\phi(w+1) = \wt\phi(\Phi_A^{1}(f_1(z))) = \Phi_A^{2,\ext}(\phi(f_1(z))) = \Phi_A^{2,\ext}(f_2(\phi(z))) = \Phi_A^{2,\ext}(\phi(z))+1 = \wt\phi(\Phi_A^{1}(z))+1 = \wt\phi(w)+1$.
  
  The relation \cref{eq:2} together with $T_1(Q_A^1)\subset Q_A^1$ imply that $\wt \phi$ uniquely extends on $\Z+Q_A^1$ to a map satisfying $\wt\phi\circ T_1 = T_1 \circ \wt \phi$.
  The right square of the diagram is then satisfied.
  By \Cref{prop:brim}, $\Z+Q_A^1 =\C$, so now $\wt\phi$ is defined on $\C$.
  Of course the left part of the diagram still commutes since it only concerns points $z\in P_A^1$, which are already dealt with.
\end{proof}

Let us call \emph{small sector} sets of the form $S[r] := S\cap r\D$ where $S$ is an open sector based on $0$ and of opening angle $<2\pi$.

\begin{proposition}\label{p:SCQuotientBiholom}
  Assume that $P_A^1$ contains an  $\alpha$-petal for some $\alpha\in(0,\pi)$ (resp.\ is a large petal).
  Then
  \begin{enumerate}
    \item For all\footnote{For any half-opening angle $<\alpha$ (resp. $<\pi$) there is such a small sector, see \Cref{cor:sec1,r:sectLarge}.} small sector $S[r]$, $r>0$, of half-opening angle $<\alpha$ (resp. $<\pi$) contained in $P_1$ and bisected by the attracting axis of $f_1$, we have \[\phi(z)\sim \frac{a_1}{a_2} z\]
    as $z$ tends to $0$ within $S[r]$.
    \item The map $\wh\phi $ is an isomorphism of the cylinder $\CZ$ fixing both ends, i.e.\  it is a translation $w\mapsto w+\sigma$ for some $\sigma\in \CZ$.
    The map $\wt\phi$ is a translation too, of vector $\wt \sigma$ satisfying $\piZ(\wt\sigma)=\sigma$.
  \end{enumerate}
\end{proposition}

Before giving its proof, let us name the quantities appearing in the proposition.

\begin{definition}\label{def:phase}%
  Under the assumption of the previous proposition, we call $\sigma$ the \emph{phase shift} and $\wt \sigma$ the \emph{lifted phase shift}, associated to $\phi$.
\end{definition}

We sum this up in the following commuting diagram (that naturally fits on a cube but for readability we draw it flat):
\begin{equation}\label{eq:cd}
\xymatrix{ 
  B_{f_1}/f_1 \ar[d]^{\ov\phi} \ar@/^8mm/[rrr]^{\ov\Phi_A^{1,\ext}} & P_A^1 \ar[l]_{\pi_{f_1}} \ar[d]_{\phi} \ar[r]^{\Phi_A^{1,\ext}} & \C \ar[d]^{T_{\tilde \sigma}} \ar[r]^{\pi_{\Z}} & \CZ \ar[d]^{T_{\sigma}}
  \\
  B_{f_2}/f_2 \ar@/_8mm/[rrr]_{\ov\Phi_A^{2,\ext}} & B^0_{f_2} \ar[l]^{\pi_{f_2}} \ar[r]_{\Phi_A^{2,\ext}} & \C \ar[r]_{\pi_{\Z}} & \CZ
}
\end{equation}

\begin{proof}[Proof of \Cref{p:SCQuotientBiholom}]

  Let $P_A^2$ be a large petal.
  Let $z_0 \in P_A^1$.
  For $n$ big enough $f_1^n(z_0)\in U$
  We denote 
  \[ X = \setof{ f_1^n(z_0) }{ n \in \N } \subset P_A^1\] 
  and: 
  \[ Y = \phi(X) = \setof{ f_2^n(\phi(z_0)) }{ n \in \N}  \subset B_{f_2} \] 
  We are goint to apply \Cref{p:equivSect}, stated and proved later in this article, to these sets $X$ and $Y$, and to $U =
  P_A^1$, $V = B_{f_2}$, and some $\kappa_0>0$, to deduce that $\phi(z)$ is asymptotic to a linear map on appropriate sectors.
  Actually, to apply the proposition, we will have to rotate $V$ so that the attracting axis coincide with the negative real axis.
  
  For this we notice that (see \Cref{eq:sim} in \Cref{sec:rappel} or use the asymptotic-equivalent of the inverses of Fatou coordinates),
  $f_1^n(z_0) \sim -\frac{1}{a_1 n}$ and
  $f_2^n(\phi(z_0)) \sim - \frac{1}{a_2 n}$.
  We have $f_2^n(\phi(z_0)) = \phi(f_1^n(z_0))$ so
  $\phi(f_1^n(z_0)) \sim - \frac{1}{a_2 n} \sim \frac{a_1}{a_2} \times -\frac{1}{a_1 n} \sim \frac{a_1}{a_2} f_1^n(z_0)$, i.e.\ 
  \[ \phi(x) \sim \frac{a_1}{a_2} x\text{ when }x \in X,\ x\to 0 \]
  Also note that $\lim_{n\to \infty} \frac{f_1^{n+1}(z)}{f_1^n(z)} = 1$.
  Then using \Cref{lem:bz} 
  we can apply \Cref{p:equivSect} to $X$,  $a_2 Y$, $a_2\phi$, $U = P_A^1$, $V = a_2 B_{f_2}$, $S_X=S$ and $\kappa_0 = a_1$ and deduce that:
  for all infinite punctured closed sector $S$ at $0$ contained in $S_X$: 
  \[ a_2\phi(z) \sim a_1 z\ \text{ when } z\tends 0,\ z\in S[r] \]
  This proves the first point.
  
  Consider a sector $S[r]$ as in the first point (which exists).
  By decreasing the angle and $r$ we may assume that $\ov{S[r]}-\{0\}\subset P_A^1$.
  Then under the condition $z\in S[r]$, we have $\Phi_A^1(z)\tends \infty$ $\iff$ $z\tends 0$.
  By the asymptotic-equivalent of \Cref{lem:asymp1} and, $(\Phi_A^{1})^{-1}(w)\sim -w/a_1$ as $w\to \infty$ within $\Phi_A^1(S[r])$.
  By what we just proved, $\phi((\Phi_A^{1})^{-1}(w)) \sim -w/a_2$.
  By \Cref{lem:asymp1} again, $\Phi_A^{2,\ext}(\phi((\Phi_A^{1})^{-1}(w)))\sim w$.
  Recall that $\wt\phi \circ \Phi_A^{1,\ext}(z) = \Phi_A^{2,\ext} \circ \phi(z)$ for $z\in P_A^1$, so $\forall w\in \Phi_A^1(S[r])$,   $\wt\phi(w) = \Phi_A^{2,\ext} \circ \phi((\Phi_A^1)^{-1}(w))$. We have thus proved:
  \[\wt\phi(w)\sim w\text{ as }w\in\Phi_A^1(S[r])\text{ tends to }\infty.\]
  Still using \Cref{lem:asymp1}, we have $\piZ\circ\Phi_A^1(S[r]) = \CZ$.
  So the estimate above implies in particular that 
  \[ \Im \wh\phi (w) \to \pm \infty \text{ when } w\in\CZ\text{ tends to }\pm i \infty \]
 %
%
  The singularities of $\wh\phi$ at $\pm i \infty$ are thus removable. Via the biholomorphism $\CZ \cup \{ \pm i \infty\} \to \wh\C$ given by $w \mapsto \exp(2i\pi w)$, the map $\wh\phi$ induces a map from $\wh\C$ to itself admitting $\{0\}$ and $\{\infty\}$ as completely invariant sets (i.e.\ invariant by direct and inverse images).
  It is thus of the form $z \mapsto C z^n$ with $C \in \C^*$ and $n \in \N^*$.
  Hence the map $\wt\phi$ must be a restriction of $z \mapsto n z + \sigma$ for some $\sigma\in\C$.
  Since $\wt\phi$ commutes with $z \mapsto z + 1$, we must have $n = 1$, thus $\wt\phi$ is a translation.
  This allows us to finish the proof: $\wh\phi $ is a biholomorphism of the cylinder preserving the ends.
\end{proof}

The proposition is not true for all petals: see \Cref{sub:ex:5}.

\begin{corollary}\label{c:SCInjPetAttr} 
  Under \Cref{h:1} and the assumption that $P_A^1$ contains an $\alpha$-petal for some $\alpha>0$, 
  the map $\phi$ is injective.
\end{corollary}

\begin{proof}
  By point 2 of \Cref{p:SCQuotientBiholom}, $\wt\phi(w)=w+\wt\sigma$.
  The corollary then follows from: $\forall z\in P_A^1$, $ \Phi_A^{2,  \ext } \circ \phi (z) = \wt\phi \circ \Phi_A^{1,\ext}(z) = T_{\wt\sigma}\circ \Phi_A^1(z)$ and injectivity of $\Phi_A^1$ (whose domain is $P_A^1$).
\end{proof}

\begin{proposition}(Semi-conjugacies send big enough petals to petals of the same nature)\label{p:imSCPetAttr0}
  Under \Cref{h:1} and the assumption that $P_A^1$ contains an $\alpha$-petal for some $\alpha>0$, then $P_A^2:=\phi(P_A^1)$ is an attracting petal for $f_2$, and if $P_A^1$ is respectively an $\alpha$-petal, large, very large, then $P_A^2$ is of the same type.
\end{proposition}
\begin{proof}
  The result follows from \Cref{lem:W} applied to $f=f_2$ and $W=\phi(P_A^1)$. Let us check each point of the lemma:
  \begin{itemize}
    \item $\phi(P_A^1)\subset B_{f_2}$ by hypothesis.
    \item $\phi(P_A^1)$ is homeomorphic to $P_A^1$, since $\phi$ is a biholomorphism by \Cref{c:SCInjPetAttr}. So $\phi(P_A^1)$ is homeomorphic to the unit disk.
    \item The semi-conjugacy property of $\phi$ and the forward invariance of $P_A^1$ by $f_1$ shows that $f_2(\phi(P_A^1)) \subset \phi(P_A^1)$.
    \item $\Phi_A^{2,\ext} \circ \phi = T_{\tilde \sigma} \circ \Phi_A^1$ on $P_A^1$, $\phi$ is injective on $P_A^1$ by \Cref{c:SCInjPetAttr} and $\Phi_A^{1,\ext}$ is injective on $P_A^1$, so $\Phi_A^{2,\ext}$ is injective on $\phi(P_A^1)$.
    \item $\piZ \circ \Phi_A^{2,\ext}(\phi(P_A^1)) = \piZ \circ T_{\tilde\sigma} \circ \Phi_A^{1,\ext}(P_A^1)) = \CZ$ by the brimming property \Cref{prop:brim}.
  \end{itemize}  
Moreover, the petal $P_A^2=\phi(P_A^1)$ is of the same nature as $P_A^1$ since $\Phi_A^{2,\ext}(\phi(P_A^1)) = T_{\tilde\sigma}(\Phi_A^1(P_A^1))$.
\end{proof}

\subsection{Local semi-conjugacy on immediate basins}\label{sub:lscib}

Let $f_i : \Dom  f_i \to \C$, $i=1,2$ be two holomorphic functions, where $\Dom  f_i \subset \C$ is an open neighborhood of $0$, with simple parabolic points at $0$. We write the Taylor expansion at $0$: 
\[ f_i(z) = z + a_i z^2 + o(z^2) \]

As earlier, for all open set $U_i$ of $\C$ containing $0$, we denote by $U_i^0$ the connected component of $U_i \cap B_{f_i}$ containing a germ of the attracting axis of $f_i$.
Recall that $U_i^0\subset B_{f_i}^0$.

\begin{definition}\label{d:semiConjLoc}%
  A holomorphic map $\phi$ is said to \emph{locally semi-conjugate $f_1$ to $f_2$ on their immediate basins} if there exists an open neighborhood $U_1$ of $0 \in \C$ such that the domain of $\phi$ is $U_1^0$, $\phi$ takes values in $B_{f_2}$ and $\phi$ is a semi-conjugacy, namely: 
  \[ \phi \circ f_1(z) = f_2 \circ \phi(z) \mbox{ for all $z \in U_1^0$ such that $f_1(z) \in U_1^0$.} \]
\end{definition}

\begin{lemma}\label{lem:sclMapToImm}
  Under these conditions, $\phi$ maps in $B_{f_2}^0$.
\end{lemma}
\begin{proof}
  The set $U_1^0$ contains an orbit $z_n$ of $f_1$.
  Let $P_A^2$ be any petal for $f_2$.
  Then $P_A^2\cap \phi(U_1^0)\neq \emptyset$: indeed $\phi(z_n)\in \phi(U_1^0)$ and since $\phi(z_0)\in B_{f_2}$, the sequence $\phi(z_n)=f_2^n(\phi(z_0))$ eventually enters $P_A^2$.
  Since $\phi(U_1^0)$ is connected, contained in $B_{f_2}$, and contains points in $P_A^2$,
  the lemma follows.
\end{proof}

\begin{remark}
  We require neither that $\phi$ admits a continuous extension at $0$, nor that $\phi(U_1^0)$ is included in a small neighborhood of $0$.
\end{remark}

The following is immediate:

\begin{lemma}\label{lem:restrictLC}
  Consider a neighborhood $U_1'$ of $0$ such that $U_1'\subset U_1$. Then the restriction of $\phi$ to ${U_1'}^0$ is a local semi-conjugacy on immediate basins too.
\end{lemma}

By \Cref{rk:smallVLPetals} there exists a large attracting petal $P_A^1$ for $f_1$ that is contained in $U_1$.
Then $P_A^1\subset U_1^0$ and restriction of $\phi$ to $P_A^1$ satisfies \Cref{h:1}.
By \Cref{c:SCInjPetAttr}, $\phi$ is injective on $P_A^1$.

\begin{remark}
Note that if $P_A^1$ is not itself large but is included\footnote{If one insists in petals being simply connected, we do not know if every petal is included in a large petal.} in a large petal, then the above applies to the large one and hence $\phi$ is injective on $P_A^1$.
\end{remark}

Since $P_A^1$ is assumed large, \Cref{p:SCQuotientBiholom} applies, so by the second point of this proposition, there is a phase shift $\tilde\sigma=\tilde\sigma(\phi)$ (which depends on $\phi$ but also on the normalization of the attracting Fatou coordinates of $f_1$ and $f_2$) such that the following diagram commutes (see \Cref{def:phase} and the commutative diagram~\eqref{eq:cd}):
\[
\xymatrix{
  P_A^1 \ar[d]_{\phi} \ar[r]^{\Phi_A^1} & \C \ar[d]^{T_{\tilde\sigma}}
  \\
  B^0_{f_2} \ar[r]_{\Phi_A^{2,\ext}} & \C
}
\]

On the other hand, the phase shift is independent of the choice of the large petal $P_A^1$: indeed, any two petal intersect (see \Cref{lem:petint}), and at a common point $z$, the diagram implies that $\wt\sigma + \Phi_A^{1,\ext}(z) = \Phi_A^{2,\ext}\circ \phi(z)$.

By analytic continuation, commutation of the diagram above extends to the connected set $\Dom \phi=U_1^0$:
\begin{equation}\label{eq:cd2}
\xymatrix{
  U_1^0 \ar[d]_{\phi} \ar[r]^{\Phi_A^{1,\ext}} & \C \ar[d]^{T_{\tilde\sigma}}
  \\
  B^0_{f_2} \ar[r]_{\Phi_A^{2,\ext}} & \C
}
\end{equation}

\begin{corollary}[Petals map to petals]\label{cor:phiinjpatt}
  The map $\phi$ is injective on any attracting petal $P_A^1$ included in its domain $U_1^0$ and $\phi(P_A^1)$ is a petal for $f_2$.
\end{corollary}
\begin{proof}
  By the diagram, $T_{\tilde\sigma} \circ \Phi_A^{1,\ext} = \Phi_A^{2,\ext} \circ \phi$ holds on $P_A^1$, and the left hand side is injective on $P_A^1$, so $\phi$ is injective on $P_A^1$ and $\Phi_A^{2,\ext}$ is injective on $\phi(P_A^1)$.
  The set $\phi(P_A^1)$ satisfies all the conditions of \Cref{lem:W} (the last one, brimming, comes from $\Phi_A^{2,\ext} \circ \phi(P_A^1) = T_{\tilde\sigma} \circ \Phi_A^{1,\ext}(P_A^1)$ and the brimming property for $P_A^1$, \Cref{prop:brim}) so is a petal.
\end{proof}

\begin{corollary} \label{p:exprPsPetA}
  Assume $P_A^1$, $P_A^2$ are petals with $P_A^1\subset U_1^0$ and that $\phi(P_A^1) \subset P_A^2$.
  Then we have for $z \in P_A^1$:
  \[ \phi(z) = (\Phi_A^2)^{-1} \circ T_{\tilde{\sigma}} \circ \Phi_A^1(z) \]
  where we denote $\Phi_A^2 = \Phi_A^{2,\ext}|_{P_A^2}$.
\end{corollary}
\begin{proof}
  This follows from the diagram and $\Phi_A^{i,\ext}$ being injective on petals. 
\end{proof}

\Cref{p:exprPsPetA} implies by uniqueness of analytic extension that $\phi$ is uniquely determined by $\tilde{\sigma}$. More precisely: 

\begin{proposition}[Uniqueness of the local semi-conjugacy] \label{p:uniPsConj}
  Let $\Phi_A^1, \Phi_A^2$ be Fatou coordinate of $f_1, f_2$. Let $\phi$ and $\psi$ be local semi-conjugacies from $f_1$ to $f_2$ on their immediate basins admitting the same lifted phase shift. Denote $U_1$ the open set corresponding to $\phi$ and $V_1$ the one corresponding to $\psi$. Set $W_1 = U_1 \cap V_1$, and let $W_1^0$ be the connected component of $W_1$ containing a germ of the attracting axis of $f_1$. Then: 
  \[ 
  \phi_{|W_1^0} = \psi_{|W_1^0} 
  \]
  (and this map is a local semi-conjugacy on immediate basins).
\end{proposition}

For the next proposition we first make two remarks.
The post-composition of $\phi$ by positive powers of $f_2$ leaves the domain of definition of $\phi$ invariant.
It is not the case for pre-composition by integer powers of $f_1$.

\begin{proposition} \label{p:translSemiConj}
Let $\phi$ be a local semi-conjugacy from $f_1$ to $f_2$ on their immediate basins.
For $m \in \N$, the map $f_2^m \circ \phi$ is also a local semi-conjugacy on immediate basins. Furthermore,\footnote{Leaving the normalization of the attracting Fatou coordinates of $f_1$ and $f_2$ unchanged.}
\[ 
\tilde{\sigma}(f_2^m \circ \phi) = \tilde{\sigma}(\phi) + m
\]
Let $U_1'=f_1^{-m}(U_1^0)$.
Then $\phi \circ {f_1^m} = f_2^m \circ \phi$ holds on ${U_1'}^0$ and
\[
\tilde{\sigma}(\phi \circ {f_1^m}|_{{U_1'}^0}) = \tilde{\sigma}(\phi) + m
\]

Let $m \in \N$ and $W_1, W'_1$ be open neighborhoods of $0$ such that $f_1^m: W_1 \to W'_1$ is a biholomorphism and such that $f_1^m$ is injective on $W_1\cup f(W_1)$. Denote by $f_{1,\mathrm{loc}}^{-m}: W'_1 \to W_1$ the inverse of this biholomorphism.
Set $V'_1 = f_1^m(U_1 \cap W_1)$. 
It is an open neighborhood of $0$.
Let ${V'_1}^0$ be the connected component of $V'_1 \cap B_{f_1}$ containing an attracting petal.
The set ${V'_1}^0$ is contained in the domain of $\phi \circ f_{1,\mathrm{loc}}^{-m}$ and the restriction $\phi \circ {f_{1,\mathrm{loc}}^{-m}}|_{{V'_1}^0}$ is a local semi-conjugay from $f_1$ to $f_2$ on their immediate basins. Furthermore:
\[
\tilde{\sigma}(\phi \circ {f_{1,\mathrm{loc}}^{-m}}|_{{V'_1}^0}) = \tilde{\sigma}(\phi) - m
\]
\end{proposition}

\begin{proof}
  We have ${V_1'}^0 \subset V_1'\subset f_1^m(W_1) = W_1' = \Dom(f_{1,\mathrm{loc}}^{-m})$.
  Consider any attracting petal $P_A^1$ of $f_1$ contained in $U_1\cap W_1$, in particular in $U_1^0$.
  The set $f_1^m(P_A^1)$ is connected and contained in $V_1'$ and in $P_A^1$, hence it is contained in ${V_1'}^0$.
  So $P_A^1 \subset f_{1,\mathrm{loc}}^{-m}({V_1'}^0)$.
  The set $f_{1,\mathrm{loc}}^{-m}({V_1'}^0)$ is connected, contained in $U_1$ and contains $P_A^1$, so it is actually contained in $U_1^0$. It follows that ${V_1'}^0$ is contained in the domain of $\phi \circ f_{1,\mathrm{loc}}^{-m}$.
  
  It is clear that $f_2^m\circ\phi$ and $\phi\circ f_1^m|_{{U_1'}^0}$ are semi-conjugacies, since $f_i^m$ commutes with $f_i$.
  For $\phi\circ f_{1,\mathrm{loc}}^{-m}$, this follows similarly from the claim ``$\forall z\in W_1'$ such that $f_1(z)\in W_1'$, $f_{1,\mathrm{loc}}^{-m} \circ f_1(z) = f_1 \circ f_{1,\mathrm{loc}}^{-m} (z)$'', which we now justify: $f_1^m \circ f_{1,\mathrm{loc}}^{-m}$ is defined and is the identity on $W_1'$ and this applies to $z$ and $f_1(z)$: $f_1^m ( f_{1,\mathrm{loc}}^{-m}(z)) = z$ and $f_1^m ( f_{1,\mathrm{loc}}^{-m}(f_1(z))) = f_1(z)$.
  So $f_1^m (f_{1,\mathrm{loc}}^{-m}(f_1(z))) = f_1(f_1^m ( f_{1,\mathrm{loc}}^{-m}(z)) )= f_1^m ( f_1 ( f_{1,\mathrm{loc}}^{-m}(z)))$.
  Now $f_1( f_{1,\mathrm{loc}}^{-m}(z))\in f_1(W_1)$ and $f_{1,\mathrm{loc}}^{-m}(f_1(z)) \in W_1$ so we get the claim by the hypothesis of injectivity of $f_1^m$ on $W_1\cup f_1(W_1)$.
   
  Pick any $z\in U_1^0$:
  $\Phi_A^{2,\ext}\circ f_2^m\circ \phi(z) = T_m\circ \Phi_A^{2,\ext}\circ \phi(z) = m+\wt\sigma + \Phi_A^{1}(z)$.
  So $\tilde{\sigma}(f_2^m \circ \phi) = \tilde{\sigma}(\phi) + m$.

  Pick any $z\in {U_1'}^0$:
  $\Phi_A^{2,\ext}\circ \phi \circ f_1^m(z) = \Phi_A^{2,\ext}\circ f_2^m\circ \phi(z) = m + \Phi_A^{2,\ext}\circ \phi(z) = m+\wt\sigma + \Phi_A^{1}(z)$.
  So $\tilde{\sigma}(\phi\circ f_1^m|_{{U_1'}^0}) = \tilde{\sigma}(\phi) + m$.

  Pick any $z\in {V_1'}^0$ such that $f_{1,\mathrm{loc}}^{-m}(z)\in {V_1'}^0$ ($z$ in a small enough petal will do):
  $\Phi_A^{2,\ext}\circ \phi \circ f_{1,\mathrm{loc}}^{-m}(z) = \wt\sigma + \Phi_A^{1,\ext} \circ f_{1,\mathrm{loc}}^{-m}(z) = \wt\sigma + \Phi_A^{1,\ext}(z) -m$.
  So $\tilde{\sigma}(\phi \circ {f_{1,\mathrm{loc}}^{-m}}|_{{V'_1}^0}) = \tilde{\sigma}(\phi) - m$.
\end{proof}

As a reciprocal:

\begin{proposition}
  Let $\phi$ and $\psi$ be two local semi-conjugacies from $f_1$ to $f_2$ on their immediate basins such that $\ov\phi  = \overline{\psi}$ (recall that $\ov\phi , \overline{\psi} : B_{f_1} / f_1 \to B_{f_2} / f_2$ are the quotient maps induced by $\phi, \psi$), i.e.\ having the same (non-lifted) phase shift in $\CZ$. Let $U_1$ be the open set corresponding to $\phi$ and $V_1$ the one corresponding to $\psi$. Set $W_1 = U_1 \cap V_1$ and let $W_1^0$ be the connected component of $W_1$ containing a germ of the attracting axis of $f_1$. Then there exists $m \in \N$ such that: 
  
  \[ 
  f_2^m \circ \phi_{|W_1^0} = \psi_{|W_1^0} \mbox{ \; or \; } \phi_{|W_1^0} = f_2^m \circ \psi_{|W_1^0} 
  \]
\end{proposition}

\begin{proof}
  We use \Cref{p:uniPsConj} and the lifted phase shift calculation given by \Cref{p:translSemiConj}.
\end{proof}

\begin{proposition}\label{p:inclSCPetA}
Let $P_A^1, P_A^2$ be attracting petals of $f_1, f_2$ with $P_A^1\subset \Dom \phi = U_1^0$. Then:
\begin{itemize}
\item  $\phi(P_A^1) \subset P_A^2$ is equivalent to $T_{\tilde{\sigma}}( \Phi_A^1(P_A^1) ) \subset \Phi_A^2(P_A^2)$
\item $\phi(P_A^1) \supset P_A^2$ is equivalent to $T_{\tilde{\sigma}}( \Phi_A^1(P_A^1) ) \supset \Phi_A^2(P_A^2)$
\end{itemize} 
\end{proposition}

\begin{figure}[!t]
\vspace{-3cm}
\def\svgwidth{\textwidth}
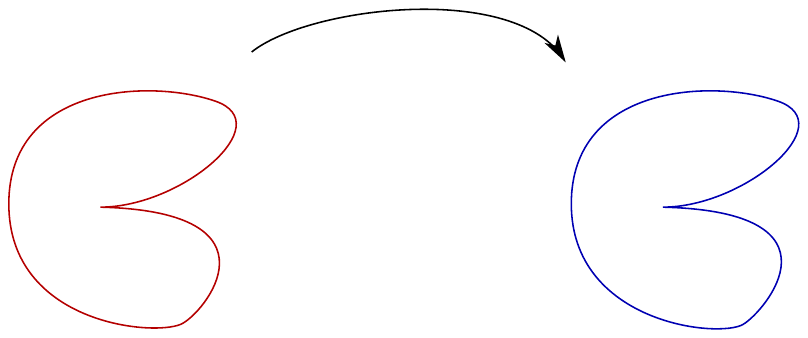
\vspace{-10cm}
\caption{Petal inclusion. In this example the petals $P_A^i$ are very large and their images by $\Phi_A^i$ are each the complement of a horizontal half strip. In this example, there is no inclusion between the sets $\phi(P_A^1)$ and $P_A^2$.}
\label{f:inclusionPetales}
\end{figure}

See \Cref{f:inclusionPetales}.

\begin{proof}
We know that $T_{\tilde{\sigma}} \circ \Phi_A^{1,\ext} = \Phi_A^{2,  \ext } \circ \phi$ holds on $U_1^0$.
So: 
\[ T_\sigma( \Phi_A^1(P_A^1) ) = \Phi_A^{2,  \ext }( \phi( P_A^1 ) ) \]
and $\Phi_A^{2,  \ext }$ induces a biholomorphism from $\phi(P_A^1)$ to $\Phi_A^{2,  \ext }( \phi( P_A^1 ) )$.
The question is thus to compare inclusion between $\phi(P_A^1)$ and $P_A^2$ and inclusion between their images by $\Phi_A^{2,\ext}$.
For this we use:

\begin{lemma} \label{l:phiUsousEnsemblePhiV}
Let $X, Y$ be two topological spaces, $U, V$ two nondisjoint open connected subsets of $X$. Let $\Phi: X \to Y$ be a map and assume that $\Phi: V \to \Phi(V)$ is proper and $\Phi(V)$ is open. Then $U \subset V$ is equivalent to $\Phi(U) \subset \Phi(V)$.
\end{lemma}
\begin{proof}
The implication is clear. For the converse, we use contraposition. Since $U$ is connected and intersects $V$, the fact that $U$ is not included in $V$ implies that $U$ intersects the boundary of $V$. Let $z \in U \cap \partial V$. Then $\Phi(z) \in \Phi(U)$ and $\Phi(z) \in \partial \Phi(V)$ since $\Phi : V \to \Phi(V)$ is proper, so $\Phi(z) \notin \Phi(V)$ since $\Phi(V)$ is open.
\end{proof}
The proposition then follows from the above lemma with $U = \phi(P_A^1)$, $V = P_A^2$ and $\Phi = \Phi_A^{2,  \ext  }$, and the same lemma permuting the role of $U$ and $V$. Let us check that the conditions are satisfied.
The sets $U$ and $V$ are open and connected.
Since $U$ contains a forward orbit of $f_2$ in its basin, it is not disjoint from $V$.
We also note that $\Phi_A^{2,  \ext  }: U \to \Phi_A^{2,  \ext  }(U)$ and $\Phi_A^{2,  \ext  }: V \to \Phi_A^2(V)$ are biholomorphisms (since, for the first one, $\Phi_A^{2,  \ext  } \circ \phi \circ (\Phi_A^1)^{-1}_{|\Phi_A^1(P_A^1)}$ is a translation, and for the second one the extended Fatou coordinates coincides with Fatou coordinates on petals), hence open and proper.
\end{proof}

\begin{corollary}\label{p:inclSCPetA4}
For all open subset $U_2$ of $\C$ containing $0$, there exists a very large attracting petal $P_A^1 \subset U_1$ such that $\phi(P_A^1) \subset U_2$.  
\end{corollary}
\begin{proof}
Let $P_A^2 \subset U_2$ be a very large attracting petal for $f_2$, which exists by \Cref{rk:smallVLPetals}.
There exists a very large attracting petal $P_A^1$ for $f_1$, and shrinking it by pushing the boundary of its three defining half planes, we can ensure that $T_{\tilde{\sigma}}(\Phi_A^1(P_A^1)) \subset \Phi_A^2(P_A^2))$, from which \Cref{p:inclSCPetA} implies that $\phi(P_A^1) \subset U_2$.
\end{proof}

\subsection{Local pseudo-conjugacies}

\begin{definition}\label{d:paireSemiConjLoc}%
A pair of maps $(\phi, \phi')$ is said to be a \emph{pair of local semi-conjugacies} of $f_1, f_2$ on their immediate basins if: 
\begin{itemize}
\item $\phi$ is a local semi-conjugacy from $f_1$ to $f_2$ on their immediate basins, 
\item $\phi'$ is a local semi-conjugacy from $f_2$ to $f_1$ on their immediate basins.
\end{itemize}
\end{definition}

\begin{definition}\label{d:pseudoConjLoc}%
Let $(\phi, \phi')$ be a pair of local semi-conjugacies on immediate basins. The pair $(\phi, \phi')$ is said to be a \emph{local pseudo-conjugacy} of $f_1, f_2$ if $\ov\phi : B_{f_1}/f_1 \to B_{f_2}/f_2$ (we already know by \Cref{p:SCQuotientBiholom} that $\ov\phi $ is a biholomorphism) has for inverse $\ov\phi ': B_{f_2}/f_2 \to B_{f_1}/f_1$.
\end{definition}

Note that there is \emph{no requirement} that $\phi$ would map to $U_2$ or $\phi'$ to $U_1$. As in \Cref{lem:restrictLC}, the following is immediate:

\begin{lemma}\label{lem:restrictPC}
  Consider neighborhoods $U_1'$ and $U_2'$ of $0$ such that $U_1'\subset U_1$ and $U_2'\subset U_2$.
  Then the pair $(\phi|_{{U_1'}^0},\phi|_{{U_2'}^0})$ is a local local pseudo-conjugacy.
\end{lemma}

The term pseudo-conjugacy is motivated by the fact that the maps $\phi, \phi'$ are mutually inverse only after passing to the quotient.

Recall that $B_{f_i}/f_i$ is identified to $\CZ$ via an isomorphism $\ov\Phi_A^{i,\ext}$ induced by $\Phi_A^{i,\ext}$, see \Cref{cor:quot}.
We saw in \Cref{p:SCQuotientBiholom} that $\ov\phi$ expresses under these isomorphism as a translation we denoted $T_{\sigma}$ for some $\sigma=\sigma(\phi)\in\CZ$, which we called the phase shift.
The following statement is thus immediate:

\begin{proposition}\label{prop:pciffssz}
  Let $(\phi, \phi')$ be a pair of local semi-conjugacies on immediate basins. The pair $(\phi, \phi')$ is a pseudo-conjugacy if and only if $\sigma(\phi) + \sigma(\phi') = \ov 0$.
\end{proposition}

Note that $\sigma(\phi)$ and $\sigma(\phi')$ depends on the normalization of the Fatou coordinates.
We will see below, as a consequence of the coming \Cref{p:sigmaPlusSigmaPrime}, that $\sigma(\phi) + \sigma(\phi')$ is independent of the normalization.

Recall that $\sigma(\phi)$, $\sigma(\phi')$, the phase shifts of $\phi, \phi'$, are defined in $\CZ$, and $\widetilde{\sigma}(\phi)$, $\widetilde{\sigma}(\phi')$, the lifted phase shifts, are defined in $\C$.
See \Cref{p:SCQuotientBiholom}.

\begin{definition}\label{def:sync}%
  A local pseudo-conjugacy is said to be \emph{synchronous} if $\tilde{\sigma}(\phi) + \tilde{\sigma}(\phi') = 0$. 
\end{definition}

\begin{example} \label{e:conjLoc}
  A semi-conjugacy $\phi$ inducing a biholomorphism $\phi: U_1^0 \to U_2^0$ is a particular case of member of a synchronous local pseudo-conjugacy, by taking $\phi' = \phi^{-1}$. In this case, the map $\phi$ could be called a \emph{local conjugacy of $f_1, f_2$ on their immediate basins}. 
\end{example}

\begin{proposition} \label{p:sigmaPlusSigmaPrime}
  The value $\tilde{\sigma}(\phi) + \tilde\sigma(\phi')$ does not depend on the choice of the normalization of Fatou coordinates. 
\end{proposition}

\begin{proof}
  We have, denoting $\Phi_A^i=\Phi_A^{i,\ext}|_{P_A^i}$:
  \begin{align*}
    {T_{\tilde\sigma(\phi)}}_{|\Phi_A^1(P_A^1)} &= \Phi_A^{2,\ext} \circ \phi \circ (\Phi_A^1)^{-1} \\ 
    {T_{\tilde\sigma(\phi')}}_{|\Phi_A^2(P_A^2)} &= \Phi_A^{1,\ext} \circ \phi' \circ (\Phi_A^2)^{-1}
  \end{align*}
  Choose $P_A^1$, $P_A^2$ such that $\phi(P_A^1) \subset P_A^2$.
  This is possible by, for instance, \Cref{p:inclSCPetA}.
  By composition, we get that ${T_{\tilde{\sigma}(\phi) + {\tilde\sigma}(\phi')}}_{|\Phi_A^1(P_A^1)} = \Phi_A^1 \circ \phi' \circ \phi \circ (\Phi_A^1)^{-1}$.
  This translation does not depend on the normalization of Fatou coordinates.
  More precisely, if we had chosen an other Fatou coordinate $T_{\alpha} \circ \Phi_A^1$, the resulting function would be the same (up to restriction) since ${T_{\tilde{\sigma}(\phi) + {\tilde\sigma}(\phi')}}$ commutes with $T_\alpha$.
  This is equivalent to saying that $\tilde{\sigma}(\phi) + {\tilde\sigma}(\phi')$ is independent of the choice of normalization.
\end{proof}

\begin{remark}
In the definition of pseudo-conjugacy, the condition that $\overline{\phi'} \circ \ov\phi  = \id$ on $B_{f_1}/f_1$ may also be reformulated in the following equivalent way. For all $z \in U_1^0$, there exists $n_0 \in \Z$ and $N \geq \max(0, - n_0)$ such that for all $n \geq N$ we have the three following conditions (the first two are automatically true for big enough $n$):  
\begin{enumerate}
\item $f_1^n(z) \in U_1^0$
\item $\phi(f_1^n(z)) \in U_2^0$
\item $\phi' \circ \phi(f_1^n(z)) = f_1^{n+n_0}(z)$
\end{enumerate}
\end{remark}

\begin{proof}
We first show the direct implication.
Let $N_0$ such that for all $n \geq N_0$ the points 1. and 2. are true.
Denote $\pi_{f_i}:B_{f_i}\to B_{f_i}/f_i$ the quotient map.
Notice that for $n \geq N_0$ we have by definition of $\ov\phi $ and $\overline{\phi'}$: 
\[
 \overline{\phi'} \circ \ov\phi \circ \pi_{f_1} (z) 
 = \overline{\phi'} \circ \ov\phi \circ \pi_{f_1} (f_1^n(z))
 = \overline{\phi'} \circ \pi_{f_2} \circ \phi (f_1^n(z))
 = \pi_{f_1} \circ \phi' \circ \phi(f_1^n(z))
\]
If $\overline{\phi'} \circ \ov\phi  = \id $, we must have $\pi_{f_1}(z) = \pi_{f_1} \circ \phi' \circ \phi(f_1^{N_0}(z))$.
Hence there exists $a, b \in \N$ such that $f_1^a(z) = f_1^b\Big( \phi' \circ \phi(f_1^{N_0}(z)) \Big) = \phi' \circ \phi(f_1^{N_0+b}(z))$.
Let $N = N_0 + b$. We have $\phi' \circ \phi(f_1^N(z)) = f_1^a(z)$, hence for all $n \geq N$ we get $\phi' \circ \phi(f_1^n(z)) = f_1^{n + a - N}(z)$.
Thus we take $n_0 = a - N$ and this finishes the proof since $N \geq - n_0 = N-a$.

Let us show the converse implication. Let $z \in U_1^0$. The equality $\phi' \circ \phi(f_1^n(z)) = f_1^{n+n_0}(z)$ for all integer $n$ bigger than $N$ implies that $\overline{\phi'} \circ \ov\phi (z \bmod f_1) = z \bmod f_1$. This is valid for all $z \in U_1^0$, hence $\overline{\phi'}\circ \ov\phi  = \id $.
\end{proof}

\begin{remark} 
  As we can see from the proof, if one $z\in U_1^0$ satisfies the conditions of the previous remark, then $\ov\phi'\circ\ov\phi=\on{id}$.
  Indeed, $\ov\phi'\circ\ov\phi$ is a translation of the cylinder $\C/\Z$ to which $B_{f_1}/f_1$ is isomorphic, as we have seen.
\end{remark} 

\begin{remark} 
  We will not use this property, but one can show that $n_0$ may be chosen independently of $z$ and equal to $\wt\sigma(\phi)+\wt\sigma(\phi')$ (which is in $\Z$ since $\sigma(\phi)+\sigma(\phi')=0$).
\end{remark} 

\begin{proposition} \label{p:exPSSync}
The existence of a local pseudo-conjugacy between $f_1, f_2$ is equivalent to the existence of a synchronous local pseudo-conjugacy.
\end{proposition}

\begin{proof}
Let $(\phi, \phi')$ be a local pseudo-conjugacy. Since $\sigma(\phi') = - \sigma(\phi)$, we obtain $\tilde{\sigma}(\phi') = - \tilde{\sigma}(\phi)$ up to post-composing (or pre-composing) $\phi, \phi'$ by appropriate integer powers of $f_1, f_2$ by \Cref{p:translSemiConj}.
\end{proof}

\begin{proposition} \label{p:PsSyncPetA}
  Let $(\phi, \phi')$ be a synchronous pseudo-conjugacy.
  Let $P_A^1, P_A^2$ be attracting petals of $f_1, f_2$ included in $U_1, U_2$. Then $\phi(P_A^1) = P_A^2$ is equivalent to $P_A^1 = \phi'(P_A^2)$.
\end{proposition}

\begin{proof}
  This is an immediate consequence of \Cref{p:inclSCPetA}.
\end{proof}

Recall that an attracting petal of $f_i$ is contained in $U_i$ if and only if it is contained in $U_i^0$.

\begin{proposition} \label{p:caracPsSync}
  A pseudo-conjugacy $(\phi,\phi')$ is synchronous if and only if there exists an attracting petal $P_A^1\subset U_1$ such that $\phi: P_A^1 \to \phi(P_A^1)$ (which is a biholomorphism) has inverse $\phi': \phi(P_A^1) \to P_A^1$ (this means in particular that $\phi(P_A^1) \subset U_2)$.
  In this case, for all attracting petal $P_A^1$ such that $P_A^1 \subset U_1$ and $\phi(P_A^1) \subset U_2$, the map $\phi: P_A^1 \to \phi(P_A^1)$ is a biholomorphism with inverse $\phi': \phi(P_A^1) \to P_A^1$.
\end{proposition}

\begin{proof}
Let $(\phi, \phi')$ be a pseudo-conjugacy.
Given any two attracting petals $P_A^1\subset U_1$ and $P_A^2\subset U_2$ of $f_1$ and $f_2$, under the hypothesis that $\phi(P_A^1) \subset P_A^2$:
a first applications of \Cref{p:exprPsPetA} 
gives $\forall z\in P_A^1$, $\phi(z) = (\Phi_A^2)^{-1} \circ T_{\tilde{\sigma}} \circ \Phi_A^1(z)$
where $\Phi_A^2 := \Phi_A^{2,\ext}|_{P_A^2}$.
A second, gives $\forall z\in P_A^2$, $\phi'(z) = (\Phi_{A'}^1)^{-1} \circ T_{\tilde{\sigma}'} \circ \Phi_A^2(z)$
where $\Phi_{A'}^1 := \Phi_A^{1,\ext}|_{P_{A'}^1}$ where $P_{A'}^1:=\phi'(P_A^2)$.
It follows that
\begin{equation} \label{e:phi'phi}
  \forall z \in P_A^1, \; \phi' \circ \phi(z) = (\Phi_{A'}^1)^{-1} \circ T_{\tilde{\sigma} + {\tilde\sigma}'} \circ \Phi_A^1(z)
\end{equation}

Assume now that the pseudo-conjugacy $(\phi,\phi')$ is synchronous.
By \Cref{p:inclSCPetA4}, 
there exists a petal $P_A^1 \subset U_1^0$ such that $\phi(P_A^1) \subset U_2^0$ (we do not use the synchronous hypothesis here but only that $\phi$ is a local semi-conjugacy on immediate basins).
Denote $P_A^2:= \phi(P_A^1)$, which is an attracting petal for $f_2$ by \Cref{cor:phiinjpatt}.
Then by \cref{e:phi'phi}, 
$\forall z \in P_A^1$, $\phi' \circ \phi(z) = (\Phi_{A'}^1)^{-1} \circ \Phi_A^1(z)$.
Since $P_{A'}^1$ is a petal, it has common points with $P_A^1$.
Near such a point, $(\Phi_{A'}^1)^{-1} \circ \Phi_A^1(z)=z$.
By analytic continuation, $(\Phi_{A'}^1)^{-1} \circ \Phi_A^1(z)$ is the identity on $P_A^1$.
So $\phi'\circ\phi$ is the identity on $P_A^1$, so $P_A^1 = \phi'(P_A^2)$ and the restrictions $\phi:P_A^1\to P_A^2$ and $\phi':P_A^2\to P_A^1$ are mutual inverses.

Conversely, assume that there is a petal $P_A^1 \subset U_1$ such that $\phi(P_A^1)\subset U_2$ and such that $\phi'\circ\phi$ is the identity on $P_A^1$.
Then, letting $P_A^2:=\phi(P_A^1)$ we have $P_A^1=\phi'(P_A^2)$ and by \cref{e:phi'phi}: $\forall z\in P_A^1$, $\Phi_A^1 (z)=\Phi_A^1 \circ \phi' \circ \phi(z) = T_{\tilde{\sigma} + {\tilde\sigma}'} \circ \Phi_A^1(z)$ so taking any $z\in P_A^1$ proves that $\tilde{\sigma} + {\tilde\sigma}'=0$.
\end{proof}

\section{Hausdorff limits and an asymptotic estimate}\label{sec:hlae}

Here we grouped a few results used in the previous section.

Recall that in a metric space $A$ with distance function $d$, for any non-empty subset $B\subset A$ and $x\in A$ one defines $d(x,B) = \inf_{y\in B} d(x,y)$.
The \emph{$\eps$-neighbourhood} of $B$ denotes the set $\setof{x\in A}{d(x,B)<\eps}$.

Recall that there is a notion of Hausdorff distance on the set $\mathcal C$ of non-empty compact subsets of $\wh\C$: let $d$ be a metric on $\wh\C$ inducing its topology;
for $C, C'\in \mathcal C$ let 
\[ \delta(C,C') = \sup_{z\in C} d(z,C')\]
\[ d(C,C') = \max (\delta(C,C'),\delta(C',C)) \]
Note that $\delta(C,C')$ is also the infimum of $\eps>0$ such that $C'$ is contained in the $\eps$-neighborhood of $C'$.
We use this distance on complements to define a distance on the set of open subsets of $\wh\C$ that are not equal to $\wh\C$:
\begin{definition}\label{def:hmo}%
\[ d(U,V) := d(\wh\C-U,\wh\C-V) \]  
\end{definition}
If $U_n\tends U$ for this distance, then in particular every compact subset $K$ of $U$ is contained in $U_n$ for $n$ big enough.


%
%

\subsection{Limits of zoomed-in basins}\label{sub:bz}

\begin{lemma}\label{lem:bz}
  Let $f:\Dom (f)\to \C$ be holomorphic with $\Dom(f)$ an open subset of $\C$ containing $0$ and assume that $f(z) = z + a z^2 + \mathcal O(z^2)$ at $0$, with $a\in\C^*$.
  Let $B_f$ be the parabolic basin of $0$ and $B^0_f$ the immediate basin.
  For $t\in\R$ denote $s_t(z) = taz$.
  Then as $t$ tends to $+\infty$, the rescaled set $s_t(B^0_f)$ tends for the distance of \Cref{def:hmo} to $\C\setminus\R_{\geq 0}$.
\end{lemma}

Note that $f$ ``lives'' in $\C$. If we allowed for domain and range to be $\wh\C$ then the conjugate of $z\mapsto z+1$ by $z\mapsto 1/z$ would be a counterexample, for then $B_f = \wh\C-\{0\}$.
Note also that our hypotheses do not imply that $B^0_f$ is simply connected. A counterexample can be obtained by conjugating by $z\mapsto 1/z$ the restriction of the map $T: z\mapsto z+1$ to $\wh\C\setminus(D\cup T_{-1}(D))$ where $D=\ov B(0,1/10)$: then $B_f$  is the image by $\iota$ of $\C\setminus\bigcup_{n\in\N} T_{-n} (D)$.

\begin{proof}
  By compactness of the set of non-empty compact subsets of $\wh\C$ with the Hausdorff distance, it is enough to consider a subsequence $t_n$ such that the domains $s_{t_n}(B^0_f)$ converge.
  
  Let $P_R$ be a repelling $\pi/2$-petal and $P_A$ be a large attracting petal. Then $P_R\cup P_A\cup\{0\}$ is a neighborhood of $0$.
  
  We first prove that $P_R$ contains at least one point $p$ not in $B_f$.
  By contradiction, if not, then there is a circle $C$ of center $0$ and small radius and contained in $P_A\cup P_R$ hence in the basin. Since $C$ is compact, $f^n$ tends uniformly to $0$ on $C$.
  By the maximum principle, $f^n$ tends uniformly to $0$ on the disk bounded by $C$.
  So $(f^n)'(0)$ tends to $0$ as $n$ tends to $+\infty$, which contradicts $f'(0)=1$.
  
  Knowing this, the full inverse orbit $p_n$ of $p$ in $P_R$ is in the complement of the basin too.
  The asymptotics of this orbit is as follows: $p_n\sim 1/an$.
  It follows that Hausdorff limit of $s_t(\wh\C - B_f)$ contains $[0,+\infty]$.
  
  For the converse inclusion, \Cref{r:sectLarge} provides for all $\alpha'=\pi-\eps$ with $\eps>0$ a small sector $S$ contained in $P_A$, hence in $B_f^0$, bisected by the attracting axis, and of half-opening $\alpha'$.
  The image $s_t(S)$ is a sector bisected by the negative real axis and when $t$ tends to $+\infty$, its complement tends to the full closed sector $S_\eps = [0,+\infty]e^{i[-\eps,\eps]}\subset\wh\C$.
  This shows that the limit of $s_t(\wh\C - B^0_f)$ is contained in $S_\eps$ for all $\eps>0$, hence in $[0,+\infty]$.
\end{proof}

The lemma above is valid in more generality (for instance for all rational maps of degree $\geq 2$) but we will not need such a generalization here.

\subsection{Extending some asymptotic equivalent}\label{sub:aux}

\begin{notation} \label{n:AlambdaMr}
  Let $\Ss$ be the set of (infinite) punctured closed sectors at $0$ of $\C$, i.e.\ of the form $\setof{z\in \C^*}{\arg(z) \in [\alpha,\beta]}$ where $\alpha<\beta$ and $\beta-\alpha<2\pi$.
  We denote for $A \subset \C$ and $r>0$: 
  \[A[r] = A \cap r \D\]
  For two sectors $S$ and $S'$ as above we denote
  \[S \subset_c S'\]
  whenever $\partial\D \cap S'$ is compactly contained in $\partial\D \cap S$, or equivalently $\alpha<\alpha'<\beta'<\beta$.
\end{notation}

The following proposition was used in the proof of \Cref{p:SCQuotientBiholom}.

\begin{proposition} \label{p:equivSect}
  Let $U$ be an open subset of $\C$ containing $S_U[r_0]$ for some $S_U \in \Ss$ and some $r_0>0$.
  
  Let $V$ be an open subset of $\C$ such that:
  when $\lambda\in\R$ tends to $+\infty$, the set 
  $\lambda V$ converges, for the metric of \Cref{def:hmo}, to an open infinite sector $V_\infty = \C-S_\infty$ where either $S_\infty \in \mathcal S$ of $S_\infty$ is a closed half-line from $0$.
  This implies that $0$ belongs to the boundary of $V$.
  
  Let $X = \setof{ x_n }{ n \in \N }$ be a subset of $U-\{0\}$ ($U$ may or may not contain $0$), such that $(x_n)$ converges in the Euclidean sense to $0$ and such that
  \[\frac{|x_{n+1}|}{|x_n|}\underset{n\to \infty}\tends 1\]
  We also assume that $X$ is included in a sector $S_X \in \Ss$ such that $S_X \subset_c S_U$.
  
  Let $\phi: U \to V$ be a holomorphic function. Suppose that there exists $\kappa_0 \in \C^*$ such that $\phi(z) \sim \kappa_0 z$ when $z \in X$ tends to $0$ (this amounts to saying that $\phi(x_n)\sim \kappa_0 x_n$ when $n \to \infty$).
  Then:
  \begin{itemize}
    \item for all sector $S \in \Ss$ with $S\subset_c S_U$, we have $\phi(z) \sim \kappa_0 z$ when $z \to 0$, $z \in S$.
  \end{itemize}
  As a consequence:
  \begin{itemize}
    \item for all sectors $S, S' \in \Ss$, with $S' \subset_c S \subset_c S_U$, there exists $r,r'> 0$ such that $S[r] \subset U$ 
    and $\phi(S[r]) \supset \kappa_0 S'[r']$.
  \end{itemize}
\end{proposition}

Note that under these hypotheses, $V$ is a hyperbolic subset of $\wh\C$: indeed it is enough to omit three points to be hyperbolic.

\begin{figure}[!t]
  \def\svgwidth{\textwidth}
  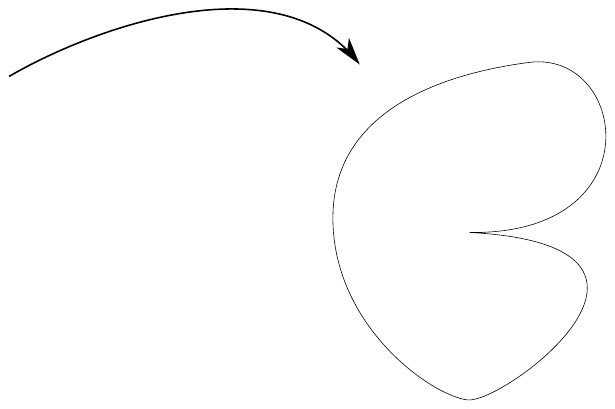
  \vspace{-14cm}
  \caption{Intermediate proposition}
  \label{f:propositionIntermediaire}
\end{figure}

\begin{remark}
  One may replace $X$ by any set which contains such a sequence $(x_n)$, for instance any curve converging towards the origin and included in a sector $S_X \in \Ss$, $S_X\subset_c S_U$, or $X = S_X$ itself.
\end{remark}

We now proceed to prove \Cref{p:equivSect}.


Denote for $\lambda \geq 1$, $\phi_\lambda: \lambda U \to \lambda V$ defined by $\phi_\lambda(z) = \lambda \phi(\frac{z}{\lambda})$. 
Let us show that the family $(\phi_\lambda)_{\lambda \geq 1}$ of holomorphic maps is normal on the open subset $G = S[r_0]$ of $\C$. 
Since $\lambda V \tends \C \setminus S_\infty$ when $\lambda \to \infty$, we have in particular that for all $z_0 \in S_\infty$, $\epsilon > 0$, there exists $\lambda_1 > 0$ such that for $\lambda \geq \lambda_1$, $B(z_0, \epsilon) \cap (\lambda V)^c \neq \emptyset$.
Taking some $z_0$ with $|z_0|=1$, $\epsilon = \frac{1}{2}$, we get for $\lambda \geq \lambda_1$ that there exists $z(\lambda) \in B(z_0, \frac{1}{2}) \cap (\lambda V)^c$.
Note that the function $z$ is not necessarily continuous.
The family $(\phi_\lambda)_{\lambda \geq a}$ is normal on $G$ since $\phi_{\lambda}(\lambda U)$ avoids three points, $0, z(\lambda), \infty$ where $z(\lambda)$ belongs to a fixed compact set avoiding $0, \infty$.\footnote{This is a classical consequence of Montel's theorem. We recall that for any pair of injective functions $a:\{1,2,3\} \to \wh\C$ and $b:\{1,2,3\} \to \wh\C$ there exists a unique (explicit) homography sending each $a(i)$ to $b(i)$ and that it depends continuously on the $a(i)$ and $b(i)$.}

Let $\psi : G\to\C$, $\psi = \lim \phi_{\lambda_n}|_G$ with $\lambda_n\to +\infty$ be any accumulation point of $(\phi_\lambda|_G)$ as $\lambda$ tends to $+\infty$, for the uniform local convergence on $G$. 
Consider the sequence $\lambda_n X\cup\{0\}$ of compact subsets of $\C$.
Consider them as compact subsets of the Riemann sphere $\widehat\C$ endowed with the spherical metric.
Let $C$ be an accumulation point of this sequence for the Hausdorff metric and $F=C\cap \C$.
Up to a further subsequence extraction, we may assume that $C$ is the limit of the sequence $\lambda_n X\cup\{0\}$.
Note that $F$ is a closed subset of $\C$ included in $S_X\cup \{0\}$.

\begin{lemma}
  For all $r>0$, the set $F$ contains an element $z$ of modulus $r$. In particular, $F$ admits accumulation points.
\end{lemma}

\begin{proof}
  Let $\epsilon > 0$. Let $A$ be the closed annulus delimited by circles of center $0$ and respective radii $r$, $r + \epsilon$. 
  Since $\lim_{k \to +\infty} \frac{|x_{k+1}|}{|x_k|} = 1$, there exists $k_0 \geq 0$ such that for $k \geq k_0$ we have $\frac{r}{r+\epsilon} \leq  \frac{|x_{k+1}|}{|x_k|}  \leq \frac{r+\epsilon}{r}$.
  Let $N \in \N$ be such that for all $n \geq N$ we have $|\lambda_n x_{k_0}| \geq r$.
  Recall that $\lambda_n x_k\tends 0$ as $k\to +\infty$ while $n$ is fixed.
  As $|\lambda_n x_k|$, $k\geq k_0$ passes from a value above $r$ to a value below, it must enter $A$.
  It follows that $\lambda_n X \cap A$ is nonempty for all $n\geq N$.
  When $n \to \infty$, this set converges in the Hausdorff sense to a set included in $F \cap A$.
  This set must be nonempty as a limit of nonempty sets included in a compact subset of $\C$.
  The property $F \cap A$ nonempty must be true for all $\epsilon > 0$.
  The set $F$ is closed, it thus intersects $r\, \partial \D$.
\end{proof}


Let $z \in F$. There exists an infinite subset $N$ of $\N$ and for $n\in N$ an index $k_n\in\N$ such that the point $z'_n=\lambda_n x_{k_n} \in \lambda_n X$ tends to $z$ when $n\in N$ tends to infinity.
Note that the $z'_n$ belong for $n$ big enough to a compact neighborhood $K$ of $z$ in $\C\setminus S_\infty$.
We have $ \phi_{\lambda_n}(z'_n) = \lambda_n \phi(x_{k_n}) \sim \lambda_n \kappa_0 x_{k_n} = \kappa_0 z'_n$, so
$ \phi_{\lambda_n}(z'_n) \underset{n \to \infty}{\longrightarrow} \kappa_0 z$
where $\kappa_0$ is given in the hypothesis of \Cref{p:equivSect}.
So $\psi(z) = \kappa_0 z$.

We just proved that $\psi(z) = \kappa_0 z$ for all $z \in F$.
By the preceding lemma and the identity theorem, we must have $\psi(z) = \kappa_0 z$ for all $z \in G$. Since the family $(\phi_\lambda)$ is normal on $G$ and has only one accumulation point, it follows that $(\phi_\lambda)$ locally converges to $\kappa_0 \,\mathrm{id}$ on $G$.

This implies the first part of \Cref{p:equivSect}.

\medskip

The first part implies the second part by the argument principle.

\section{Proof of \Cref{t:conjRevetBasLoc} and \Cref{prop:complement}}\label{sec:proof}

\subsection{Direct implication}

We first show the direct implication of \Cref{t:conjRevetBasLoc} and \Cref{prop:complement}. Recall that the subsets $\Dd_1, \Dd_2$ of $\CZ$ are the domains of $\h_1, \h_2$, and that $\Dd_1^{\pm}, \Dd_2^{\pm}$ are the connected components of these domains containing a punctured neighborhood of $\pm i \infty$.
The maps $\h_i^\pm$ are the horn maps of $\h_i$ restricted to $\Dd_i^\pm$.
We state again below the direct implication of \Cref{t:conjRevetBasLoc} and \Cref{prop:complement}:

\begin{proposition} \label{p:conjRevet}
Let $f_1$, $f_2$ denote holomorphic maps from open neighborhoods of $0 \in \C$ to $\C$, with Taylor expansions at $0$: $f_i(z) = z + a_iz^2 + o(z^2)$ where $a_i \neq 0$. 

Suppose that $f_1, f_2$ are locally semi-conjugate on their immediate basins at $0$ as per \Cref{d:semiConjLocIntro}.
Then there exists $\sigma \in \CZ$ and a pair of holomorphic maps $\psi = (\psi^+, \psi^-)$, where $\psi^\pm: \Dd_1^\pm \rightarrow \Dd_2^\pm$, such that: 
\[ \h_1^\pm = T_\sigma^{-1} \circ \h_2^\pm \circ \psi^\pm \]
where $T_\sigma$ is the translation of the cylinder $\CZ$ given by the formula $z \mapsto z+\sigma$.
\end{proposition}

\begin{compl} \label{cpl:conjRevet}
Suppose that $f_1, f_2$ are locally pseudo-conjugate at $0$ (see \Cref{d:pseudoConjLoc}). 
Then there exists $\sigma \in \C$ and a pair of biholomorphisms $\psi = (\psi^+, \psi^-)$, $\psi^\pm: \Dd_1^\pm \rightarrow \Dd_2^\pm$, such that: 
\[ \h_1^\pm = T_\sigma^{-1} \circ \h_2^\pm \circ \psi^\pm \]
\end{compl}

Below we prove the proposition in \Cref{subsub:sum,subsub:Dd1Dd2}, then the complement in \Cref{subsub:psipsiinv}.

\begin{remark}
  As already mentionned, the quotient of a petal by the equivalence relation $z\sim f(z)$ is naturally a Riemann surface $\Con$ isomorphic to $\CZ$ via the Fatou coordinate (see for instance \cite{t:LM}, from Definition~4.2.5 to Remark~4.2.14).
  In our setting, the heuristics is that $\phi$ induces a translation map $T_\sigma$ from $\Con_A^{f_1}$ to $\Con_A^{f_2}$ (this is the point of \Cref{sub:scfp,sub:lscib}), and a partial map $\psi$ from $\Con_R^{f_1}$ to $\Con_R^{f_2}$ that can be restricted to the maps $\psi^\pm : \Dd_1^\pm \to \Dd_2^\pm$.
  These maps $T_\sigma$ and $\psi^\pm$ will enable us to show the desired equality $T_\sigma \circ \h_1^\pm = \h_2^\pm \circ \psi^\pm$.
\end{remark}

\subsubsection{Construction}\label{subsub:sum}

Note that $T_\sigma$ is associated the the attracting Fatou coordinates, while $\psi$ to the repelling ones.
 
\medskip

\noindent\textbf{Construction of $T_\sigma$:} The map $T_{\sigma}:\CZ\to\CZ$ that we will use turns out to be the map $\wh\phi$, quotient of the map $\wt\phi=T_{\wt\sigma}:\C\to\C$ induced by $\phi$ and constructed in \Cref{sub:scfp}, for which we showed in \Cref{sub:lscib} that we have the commutative diagram~\eqref{eq:cd2}, which we copy below:
\begin{equation}\label{eq:cd2:copy}
\xymatrix{
  U_1^0 \ar[d]_{\phi} \ar[r]^{\Phi_A^{1,\ext}} & \C \ar[d]^{\wt\phi = T_{\tilde \sigma}}
  \\
  B^0_{f_2} \ar[r]_{\Phi_A^{2,\ext}} & \C
}
\end{equation}
where $\wt\sigma$ is the lifted phase shift of $\phi$ and $\sigma = \piZ(\wt\sigma)$ is its phase shift, see \Cref{def:phase}.

\medskip

\noindent\textbf{Construction of $\psi$:} The construction on the repelling side is harder.
Any point $w$ in $\CZ$ has an antecedent $\wt w\in Q_R^1 := \Phi_R^1(P_R^1)$ by the canonical projection $\piZ : \C \to \Z$ (\Cref{prop:brim}).
The inverse sequence $z_{-n} := \Phi_R^{-1}(\wt w-n)$, for $n\in\N$, is a backward orbit of $f_1$.
For some $w\in\C/\Z$, but not all, we have the condition
\begin{itemize}
  \item[($\alpha$)]\label[$(\alpha)$]{item:alpha} there exists $N$ such that $\forall n\geq N$, $z_{-n}\in U_1^0$,
\end{itemize}
i.e.\ the backward orbit is in $U_1^0$ for all $n$ big enough.
In this case, we get an $f_2$-backward orbit $\phi(z_{-n})$, where $n\geq N$, but it is not obvious whether or not the following condition holds:
\begin{itemize}
  \item[($\beta$)]\label{item:beta} there exists $N'\geq N$ such that $\forall n\geq N'$, $\phi(z_{-n})\in P_R^2$.
\end{itemize}
If both \hyperref[item:alpha]{$(\alpha)$} and \hyperref[item:beta]{$(\beta)$} hold, we set
\[ \psi(w) := \piZ(\Phi_R^2(\phi(z_{-n}))) \]
for some $n\geq N'$.
Recall that $\Dd_2 = (\Psi_R^{2,\ext})^{-1}(B_{f_2})$.
So $\psi$ takes values in $\Dd_2$ since $\phi$ takes values in $B^0_{f_2}\subset B_{f_2}$.
Note also that the domain of $\psi$ is contained in $\Dd_1 = \Dom  \h_1 = \piZ((\Psi_R^{1,\ext})^{-1}(B_f))$, since $z_{-n}\in U_1^0\subset B^0_{f_1}\subset B_{f_1}$ for $n\geq N$.

\begin{remark}
  The backward orbit $\phi(z_{-n})$ for $f_2$ might enter and then go out of the repelling petal $P_R^2$ several times, in such a way that the intersection of the orbit and the repelling petal may in general define distinct elements in $P_R^2 / f_2$ at each step, i.e.\ different values for $\piZ(\Phi_R^2(z_{-n}))$.
  In view of independence from the petal, it is natural to restrict, as we do above, to the case when $\phi(z_{-n})$ remains in $P_R^2$ for $n$ big enough and to select this element for the definition of $\psi(w) = \piZ(\Phi_R^2(z_{-n}))$.
\end{remark}

\begin{lemma}\label{lem:indepPsi}
  For a given $w\in\CZ$, the fact that $\psi$ is defined at $w$ and its value $\psi(w)$, are independent (i) of the representative $\wt w$, (ii) of $n\geq N'$, and (iii) of the choice of petals $P_R^1$ (iv) and $P_R^2$. 
\end{lemma}
\begin{proof}
  (i) Any two representatives $\wt w$ differ by an integer, so up to shifting indices in one of the associated inverse orbits, these inverse orbits coincide.
  (ii) For any $n\geq N'$, $\Phi_R^2(\phi(z_{-n})) = \Phi_R^2(\phi(f_1(z_{-n-1}))) = \Phi_R^2(f_2(\phi(z_{-n-1}))) = 1+\Phi_R^2(\phi(z_{-n-1}))$ so they have the same image by $\piZ$.
  (iii) For any two choices $P_R^1$ and $P_R^{1'}$, the respective representatives $\wt w$ and $\wt w'$ of $w\in\CZ$ differ by an integer.
  By $T_{-1}$ stability of the image of repelling petals in repelling Fatou coordinate, one of $\wt w$ or $\wt w'$ belongs to both sets $\Phi_R^1(P_R^1)$ and $\Phi_R^1(P_R^{1'})$.
  (iv) Consider any two choices $P_R^2$ and $P_R^{2'}$.
  Take $n$ bigger than the maximum $N''$ of the two corresponding thresholds $N'$ of condition ($\beta$), so that $\phi(z_{-n})\in P_R^2\cap P_R^{2'}$.
  Consider the respective repelling Fatou coordinates $\Phi_R^2$ and $\Phi_R^{2'}$, which are two inverse branches of $\Psi_R^{2,\ext}$ on respectively $P_R^2$ and $P_R^{2'}$.
  Let $u_{-n} = \Phi_R^2(\phi(z_{-n}))$ and $u'_{-n} = \Phi_R^{2'}(\phi(z_{-n}))$.
  Then $\forall k\geq 0$, $u_{-N''-k} = u_{-N''}-k$ and $u'_{-N''-k} = u'_{-N''}-k$.
  Note also that $\Psi_R^{2,\ext}(u_{-n}) = \phi(z_{-n}) = \Psi_R^{2,\ext}(u'_{-n})$.
  For all $n$ big enough the two sequences $u_{-n}$ and $u'_{-n}$ belong to domain $D$ where $\Psi_R^{2,\ext}$ is injective (for instance one can take $D$ to be a left half-plane, since there always exists a repelling $\pi/2$-petal).
  It follows that $u_{-n}=u'_{-n}$.
\end{proof}

To deduce from \Cref{lem:indepPsi} that $\psi$ holomorphic, we need to prove that $N'$ can be taken uniform on some neighborhood of $w$.

\begin{lemma}
  If $\wt w\in \C$ satisfies \hyperref[item:alpha]{$(\alpha)$} and \hyperref[item:beta]{$(\beta)$} with some $N$ and $N'$, then all sufficiently nearby points satisfies them with a uniform, but possibly higher, $N$ and $N'$.
\end{lemma}
\begin{proof}
  Consider $\eps>0$ such that $B(\wt w,\eps)$ is contained in the open set $Q_R^1 := \Phi_R^1(P_R^1)$. 
  By $T_{-1}$-stability of $Q_R^1$, the translate $T_{-n}(B(\wt w,\eps))$ is contained in it too.
  Let
  \[W_n :=(\Phi_R^1)^{-1}(T_{-n}(B(\wt w,\eps)))\]
  By choosing $N$ possibly higher, we can ensure that $\forall n\geq N$, $W_n \subset U_1$.
  Note that $z_{-n}\in W_n$ and recall that $z_{-n}\in U_1^0$.
  By continuity, possibly taking $\eps$ smaller, we can ensure that $W_N\subset U_1^0$.
  Then for all $n\geq 0$, $W_n \subset B_{f_1}$ (we have $f_1(W_{n+1})= W_n$).
  Also, for $n\geq N$, the connected set $W_n$ is contained in $U_1\cap B_{f_1}$ and contains the point $z_{-n}\in U_1^0$, so $W_n \subset U_1^0$, which proves \hyperref[item:alpha]{$(\alpha)$}.

  If $N$ was increased, replace $N'$ by $\max(N',N)$.
  Possibly decreasing $\eps$ we can assume that $\phi(W_{N'})$ is contained in $P_R^2$.
  Denote $g_2 = (\Phi_R^2)^{-1} \circ T_{-1} \circ \Phi_R^2 : P_R^2 \to P_R^2$, which is the unique branch of $f_2^{-1}$ sending $P_R^2$ to $P_R^2$.
  In particular, $g_2(\phi(z_{-n})) = \phi(z_{-n-1})$ since both $\phi((z_{-n})$ and $\phi(z_{-n-1})$ belong to $P_R^2$ and $f_2$ maps the second to the first.
  For $n\geq N'$, consider on the set $W_{n}$ the map
  \[ h: z\in W_n\mapsto g_2^{n-N'}\circ \phi\circ f_1^{n-N'}(z)\]
  One computes $h(z_{-n}) = g_2^{n-N'}\circ \phi\circ f_1^{n-N'}(z_{-n}) = g_2^{n-N'}\circ \phi (z_{-N'}) = \phi (z_{-n})$.
  So both holomorphic maps $h$ and $\phi|_{W_n}$: (i) are defined on the connected set $W_n$, (ii) map $z_{-n}$ to the same point $\phi (z_{-n})$ and (iii) give the same map defined on $W_n$ when post-composed with $f_2^{n-N'}$ (since $f_2^{n-N'}\circ h = f_2^{n-N'} \circ g_2^{n-N'}\circ \phi\circ f_1^{n-N'} = \phi\circ f_1^{n-N'} = f_2^{n-N'} \circ \phi$ on $W_n$). 
  It follows that $h$ and $\phi|_{W_n}$ are equal.
  In particular, $\phi(W_n)=g_2^{n-N'}(\phi(W_{N'}))$ is contained in $P_R^2$ for all $n\geq N'$.
  This proves \hyperref[item:beta]{$(\beta)$}.
\end{proof}

So $\Dom \psi$ is open and $\psi$ is holomorphic.
See also \Cref{rk:dom-psi}.

\medskip

We recall that $\h_i$ is the quotient map of $h_i = \Phi_A^{i,\ext}\circ\Psi_R^{i,\ext}$ under $\piZ$, i.e.\ $\h_i \circ\piZ =\piZ\circ h_i$.
\begin{lemma}
  For all $w\in\Dom \psi$, $T_\sigma \circ \h_1(w) = \h_2 \circ \psi(w)$.
\end{lemma}
\begin{proof}
  A presentation of the proof can be found in \cite{t:LM}, page~90 (in Section~4.5), where the focus is on quotient spaces and naturality of the constructions.
  We present here an equivalent approach using diagram chasing in the various semi-conjugacies and definition.
  First note that $\psi(w)$ is indeed in the domain of $\h_2$, because $\psi(w) = \piZ(\wt w'_n)$ with $\wt w'_n:=\Phi_R^2(\phi(z_{-n})))$ for $n\geq N'$ and hence $\Psi_R^{2,\ext}(\wt w'_n) = \phi(z_{-n}) \in B_{f_2} = \Dom  \Phi_A^{2,\ext}$ so $h_2(\wt w'_n)$ is defined.
  Second, $\h_2(\psi(w)) = \piZ (\Phi_A^{2,\ext}\circ \phi(z_{-n})) = \piZ(T_{\wt \sigma}\circ \Phi_A^{1,\ext}(z_{-n}))$ by Diagram~\eqref{eq:cd2:copy}, so $\h_2(\psi(w)) = T_\sigma\circ \piZ\circ\Phi_A^{1,\ext}((\Phi_R^1)^{-1}(\wt w-n)) 
  = T_\sigma\circ \piZ\circ h_1(\wt w-n) 
  = T_\sigma\circ \h_1\circ \piZ(\wt w -n)
  = T_\sigma\circ \h_1(w)$.
\end{proof}

There remains to prove that $\Dd_1^+$ and $\Dd_1^-$ are contained in $\Dom \psi$ and map to $\Dd_2^+$ and $\Dd_2^-$. This will be done in \Cref{subsub:Dd1Dd2}.

In the case of a pseudo-conjugacy $(\phi,\phi')$, the two local semi-conjugacies on immediate basins $\phi$ and $\phi'$ have phases $\sigma$, $\sigma'$ that add up to $0$ modulo $\Z$ by definition, and they create two pairs of maps $\psi^\pm$ and $\psi'^\pm$.
We will prove that they are mutual inverses in \Cref{subsub:psipsiinv}.

\subsubsection{The domain of $\psi$ contains $\Dd_1^+\cup \Dd_1^-$}\label{subsub:Dd1Dd2}

We prove here that $\Dd_1^+$ and $\Dd_1^-$ are contained in $\Dom \psi$ and map by $\psi$ to $\Dd_2^+$ and $\Dd_2^-$.
We treat the case of $\Dd_i^+$ since the argument for $\Dd_i^-$ is completely analogous.

For all initial choices of $w \in \Dd_1^+$,
we already know that \hyperref[item:alpha]{$(\alpha)$} holds (this is \Cref{p:petitPetRep}).
An essential part of the work will be to prove that \hyperref[item:beta]{$(\beta)$} holds too.
This would be obvious if we knew that $\phi$ tends to $0$ at $0$ but this is not part of our hypotheses (and does not hold in some examples, see \Cref{sub:examples}).
We prove it first when $|\Im(w)|$ big in \Cref{c:limPhi}; we then show that it is valid for all $w \in \Dd_1^+$, by exploiting the arc connectedness of this set and by controlling the hyperbolic length of the paths formed at each step of the construction.

\medskip

\underline{\textit{When $\Im(w)$ is big:}} In this case the situation is very favorable due to \Cref{prop:sepals}.

\begin{proposition} \label{c:limPhi}
  Let $\phi$ be a local semi-conjugacy from $f_1$ to $f_2$ on their immediate basins.
  Let $P_R^1\subset U_1$ be a repelling petal of $f_1$ and $P_R^2$ be a repelling petal of $f_2$, with respective Fatou coordinates denoted $\Phi_R^1$, $\Phi_R^2$.
  Then there exists $M > 0$ satisfying the following property.
  Let $w\in\C$ such that $| \Im( w) | > M$.
  Then $\forall n\in\Z$, $w+n\in \Dom \Psi_R^{1,\ext}$ and $z_n := \Psi_R^{1,\ext}(w+n)\in \Dom (\phi)$.
  Moreover, $\phi(z_n) \in P_R^2$ for all $n$ negative enough.
\end{proposition}

Note that in the proposition, we can take $n$ \emph{negative or positive}.

\begin{proof} 
  Consider a neighborhood $U_2$ of $0$ small enough so that $f_2$ is injective and so that the inverse $g$ of $f_2|_{U_2}$ has all its orbits tending to $0$ or exiting $U_2$.
  By \Cref{rk:smallVLPetals} there are very large attracting  petals $P_A^2\subset U_2$ and $P_A^1\subset U_1$.
  By connectedness, $P_A^1 \subset U_1^0$.
  By definition of very large petals, its image by $\Phi_A^1$ contains a union of half planes $|\Im z|>M$ and $\Re z>M$
  for some $M$.
  We can reduce $P_A^1$ to be exactly equal to this union and we can also increase the value of $M$. 
  Being smaller, the new $P_A^1$ is still included in $U_1^0$.
  We choose $M$ big enough so as to ensure that $T_{\tilde\sigma}(\Phi_A^1(P_A^1)) \subset \Phi_A^{2,\ext}(P_A^2)$, which is possible since $P_A^2$ is very large, and this implies by \Cref{p:inclSCPetA} that $\phi(P_A^1)\subset P_A^2$.
  
  By \Cref{prop:sepals} there exists $M'_1>0$ such that $\forall w\in\C$, $|\Im(w)|>M'_1$ implies $w\in\Dom \Psi_R^{1,\ext}(w)$ and $\Psi_R^{1,\ext}(w)\in P_A^1$.
  In particular if $z\in P_R^1$, applying the above to $w=\Phi_R^1(z)$ we get that $|\Im(\Phi_R^1(z))|>M'_1$ implies $z\in P_A^1$, hence $z\in U_1^0=\Dom \phi$.
  Then the sequence $z'_n := \phi(z_{-n})$ for $n\geq 0$ is a backward orbit of $f_2$ contained in $P_A^2$ hence in $U_2$, so it tends to $0$ and hence is eventually in the repelling petal $P_R^2$.
\end{proof}

This in particular proves that Properties~\hyperref[item:alpha]{$(\alpha)$} and~\hyperref[item:beta]{$(\beta)$} of \Cref{subsub:sum} hold when $\Im(w)$ is big enough.

\begin{lemma}\label{lem:psiii}
  $\Im\psi(w)\tends+\infty$ as $\Im(w) \to+\infty$.
\end{lemma}
\begin{proof}
  We proved in \Cref{lem:indepPsi} that $\psi(w)$ is independent of the choice of petals. 
  We put ourselves in the situation of \Cref{c:limPhi} and its proof, where we moreover ask that $U_2$ is small enough so that $f_2$ has an inverse branch fixing $0$ on $U_2$, so as to apply the first part of \Cref{cor:dlCornes} to $f_2$.
  In the proof of \Cref{c:limPhi}, there was no condition on $P_R^2$. Here we choose it to be included in $U_2$.
  
  Now, when $\Im(w)=\Im(w+n)$ is high, then by \Cref{p:hi}, $\Phi_A^{1, \ext}(z_n)$ has high imaginary part, and so has $T_{\tilde\sigma}(\Phi_A^{1, \ext}(z_n)) = \Phi_A^{2,\ext}(\phi(z_n))$.
  We saw in the previous proof that $\forall n\in \Z$, $\phi(z_n) \in P_A^2$ and that for $n$ negative enough, $\phi(z_n)\in P_R^2$.
  Then by the first part of \Cref{cor:dlCornes}, $\Phi_R^2(\phi(z_n))$ is high too.
\end{proof}

This implies, since $\psi^+$ is holomorphic, that its singularity at $+ i \infty$ is removable.

\begin{figure}[!t] 
  \def\svgwidth{\textwidth}
  \vspace{-3cm}
  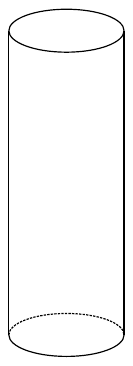
  \vspace{-8cm}
  \caption{Construction of the image of a path $\gamma$ by $\psi^+: \Dd_1^+ \to \Dd_2^+$}
  \label{f:psi+}
\end{figure}

\bigskip

\underline{\textit{General $w\in\Dd_1^+$:}} \Cref{f:psi+} provides an illustration.
Recall that $\Dom \psi$ does not depend on the choice of petals, so we choose $P_R^1$ to be a $\pi/2$-petal contained in $U_1$ (apply \Cref{prop:mil} to a restriction $f$ to a small open neighbourhood of $0$).
The set 
\[Q_R^1:=\Phi_R^1(P_R^1)\] is equal to a left half-plane defined by $\Im(w)<L_0$ for some $L_0\in\R$.
Denote 
\[\Psi_R^1:=(\Phi_R^1|_{P_R^1})^{-1}:Q_1\to P_R^1\]
Let $w \in \Dd_1^+$. Choose any $w_0 \in \Dd_1^+$ with imaginary part $>M$ where $M$ is given by \Cref{c:limPhi}.
There exists a path $\gamma$ included in $\Dd_1^+$ such that $\gamma(0) = w_0$, $\gamma(1) = w$.
Denote $\hat{\gamma}$ a lift of $\gamma$ by $\piZ$, and for $n \in \N$, $\hat{\gamma}_{-n} = T_{-n}\circ \hat{\gamma} = \hat{\gamma} - n$.
By \Cref{p:covPetCompact}, $\hat{\gamma}_{-n}([0,1])\subset Q_R^1$ for all $n$ big enough, say $n\geq n_0$.
It follows that for $n\geq n_0$, $\Psi_R^1$ is defined on the support of $\hat{\gamma}_{-n}$.
Denote for $n\geq n_0$
\[ \alpha_{-n} = \Psi_R^1 \circ \hat\gamma_{-n}\] 
which is defined on $[0,1]$, and
\[ z_{-n} = \alpha_{-n}(1) = \Psi_R^1(\wt w - n) \]
where $\wt w = \hat\gamma(1)$.
Then $\alpha_{-n}$ maps in $\Psi_R^1(Q_R^1)=P_R^1$, hence in $U_1$.
Moreover, from $\gamma([0,1])\subset \Dd$ and by definition of $\Dd$, $\alpha_{-n}([0,1])\subset B_f$.
Since moreover $\alpha_{-n}(0)$ belongs to $\Dom  \phi$ by \Cref{c:limPhi}, hence to $B_f^0$, we have that $\alpha_{-n}([0,1])$ is contained in  $B_f^0$, and hence in $U_1^0$, so $\phi$ is defined on this set.
Denote
\[\beta_{-n} = \phi \circ \alpha_{-n}
\]
for all $n \geq n_0$.
Then $f_2 \circ \beta_{-n-1} = \beta_{-n}$.

The set $\wDd_1^+ = \piZ^{-1}(\Dd_1^+)$ contains the upper half plane $\Im z>M$ so $\wDd_1^+ \cap Q_R^1$ contains the quarter plane $Q$ of equation $\Im z>M$, $\Re z<L_0$.
We denote $X$ the connected component of $\wDd_1^+ \cap Q_R^1$ that contains $Q$.
The set $\Psi_R^1(X)$ is contained in $\Psi_R^1(\wDd_1^+)$ so by  \Cref{p:tildeDpm}, it is contained in $B_{f_1}^0$ thus in $U_1^0$ (since it is also connected, contained in $U_1$ and it contains $\Psi_R^1(Q)\subset U_1^0$).
We claim that
\[T_{-1}(X)\subset X\]
Indeed $T_{-1}(\wDd_1^+) = \wDd_1^+$, $T_{-1}(Q_R^1) \subset Q_R^1$ and $T_{-1}(Q)\subset Q$, so the set $T_{-1}(X)$ is a connected subset of $\wDd_1^+ \cap Q_R^1$ that has non-empty intersection with $Q$.

To state the next claim, let us set a general notation: for all hyperbolic Riemann surface $S$ and $\gamma$ included in $S$, denote by $L_S(\gamma)$ the length of the path $\gamma$ for the hyperbolic metric of $S$.
Open subsets of hyperbolic sets are hyperbolic so the following sets are hyperbolic: $X$ (included in a half plane), $\psi_R^1(X)$ (isomorphic to $X$), $U_1^0$ (contained in $B_{f_1}^0$) and $B_{f_2}^0$.
(By \Cref{lem:hyp}, $B_{f_i}^0$ is hyperbolic.)
We have for all $n\geq n_0$:
\[ L_{B_{f_2}^0}(\beta_{-n})
\underset{i}\leq L_{U_1^0}(\alpha_{-n})
\underset{ii}\leq L_{\Psi_R^1(X)}(\alpha_{-n})
\underset{iii}= L_X(\hat{\gamma}_{-n})
\underset{iv}\leq L_X(\hat{\gamma}_0)
\]
where the successive inequalities and equalities, are justified by the following inclusions and holomorphic maps: 
\begin{align*}
  i.\quad &\phi: U_1^0 \to B_{f_2}^0 \\ 
  ii.\quad &\Psi_R^1(X) \subset U_1^0 \\ 
  iii.\quad &\Psi_R^1: X \to \Psi_R^1(X) \;\;\;\;\;\; \text{ (biholomorphism)} \\ 
  iv.\quad &T_{-n} : X  \to X
\end{align*}
The sequence of hyperbolic lengths $ L_{\phi(U_1^0)}(\beta_{-n}) $ is thus bounded from above.

Since $(\beta_{-n})$ for $n\geq n_0$ is a sequence of paths of $B^0_{f_2}$ with one extremity converging in the Euclidean sense to $0 \in \partial B_{f_2}^0$ and having a bounded hyperbolic length in $B_{f_2}^0$, this implies that $\beta_{-n}([0,1])$ has a Euclidean diameter converging to $0$. Thus $(\beta_{-n})$ converges uniformly to $0$ in the Euclidean sense.
Choose any repelling petal $P_R^2$ such that $s(P_R^2)$ contains a left half-plane where $s(z)=-1/a_2z$ (see \Cref{cor:sec1,prop:petEx} completed by \Cref{rk:smallVLPetals}).
By \Cref{lem:cptUnif} applied to a local branch of $f_2^{-1}$, $\beta_{-n}$ is contained in $P_R^2$ for all $n$ big enough.
%
In particular, $z_{-n}$ is in $P_R^2$ for all $n$ big enough, i.e.\ \hyperref[item:beta]{$(\beta)$} of \Cref{subsub:sum} holds. 

Recall that $\psi$ takes values in $\Dd_2$.
Since $\psi$ is holomorphic, hence continuous, it sends $\Dd_1^+$ to a connected subset of $\Dd_2$.
Since the two sets $\Dd_i^+$ contain upper half cylinders, \Cref{lem:psiii} implies that $\psi(\Dd_1^+)$ contains points in $\Dd_2^+$.
Since $\Dd_2^+$ is a connected component of $\Dd_2$, we get that:
\[ \psi(\Dd_1^+) \subset \Dd_2^+ \]

This completes the proof of \Cref{p:conjRevet}.

\begin{remark}\label{rk:dom-psi}
  \Cref{p:conjRevet} proves that the connected set $\Dd_1^+$ is contained in $\Dom \psi$ from the knowledge that a single point $w_0\in \Dd_1^+$ belongs to $\Dom \psi$, using paths.
  By following to the letter this proof, we can prove that 
  if a connected component of $\Dd_1$ contains a point where $\psi$ is defined, then $\psi$ is defined over the whole component; this shows that the 
  domain of $\psi$ must be a union of connected components of $\Dd_1$.
\end{remark}

\begin{remark}
  If the open set $U_1$ is small enough to apply \Cref{p:petitPetRep2}, condition \hyperref[item:alpha]{$(\alpha)$} holds \emph{if and only if} $w \in \Dd_1^+ \cup \Dd_1^-$. In this setting, $\Dom \psi = \Dd_1^+ \cup \Dd_1^-$.
\end{remark}

\subsubsection{Local pseudo-conjugacies $(\phi,\phi')$ yield inverse maps $\psi,\psi'$.}\label{subsub:psipsiinv}

We complete here the proof of \Cref{cpl:conjRevet}.

Applying the direct sense\footnote{\ldots\ of our main theorem for local-semiconjugacies on parabolic basins, which we recalled as \Cref{p:conjRevet} and just proved.} to $\phi$ and $\phi'$, we obtain two holomorphic maps $\psi^+:\Dd_1^+\to \Dd_2^+$ and $\psi'^+:\Dd_2^+\to \Dd_1^+$ such that
\begin{align*}
  T_{\sigma} \circ \h_1 &= \h_2 \circ \psi^+ 
  \\
  T_{\sigma'} \circ \h_2 &= \h_1 \circ {\psi'}^+ 
\end{align*}
By definition of pseudo-conjugacies, we have $\sigma + \sigma' = \ov 0$.
It follows that 
\[ \h_1 = T_{\sigma'} \circ T_{\sigma} \circ \h_1 = T_{\sigma'} \circ \h_2 \circ \psi^+ = \h_1 \circ {\psi'}^+ \circ \psi^+
\]
on $\Dd_1^+$.
The maps $\h_1$, $\psi$ and $\psi'$ have an erasable singularity at the top end of the cylinder, fixing this end, and $\h_1$ is locally injective near this end. It follows that $\psi'\circ\psi = \id$ near the top end. By analytic continuation:
\[\psi' \circ \psi = \on{id}\text{ on }\Dd_1^+\]

By similar arguments, the same holds on $\Dd_1^-$, and $\psi \circ \psi'=\on{id}$ on $\Dd_2^+$ and on $\Dd_2^-$.

\begin{remark}
  In \cite{t:LM}, a proof is given using instead a property of $\phi$ and $\phi'$ called \emph{pseudo invertibility} see pages~80 to~83 in that reference.
  We plan to develop this notion in a later article. 
\end{remark}

\begin{remark}\label{rk:dom-psi2}
  Contrarily as in \Cref{rk:dom-psi},  
  in the case where the domain of $\psi$ is strictly bigger than $\Dd_1^+\cup\Dd_1^-$, the arguments of the proof does not apply to show that $\psi$ is a biholomorphism over its image with inverse $\psi'$. 
  We do not know if one should expect this to be true.
\end{remark}

\subsection{Converse implication}

We now show the converse implication in \Cref{t:conjRevetBasLoc} and its complement, which we recall below for convenience as \Cref{p:recTh,compl:recTh}.

Let $f_1$, $f_2$ be holomorphic maps from open neighborhoods of $0 \in \C$ to $\C$, with Taylor expansions at $0$: $f_i(z) = z + a_iz^2 + o(z^2)$ where $a_i \neq 0$. 
We recall that we denote $\mathcal D_i^+$ the connected component containing $+i\infty$ of the domain of the non-lifted horn map $\h_i$ of $f_i$.
Similarly we denote $\mathcal D_i^-$ the connected component containing $-i\infty$ and $\h_i^+$ the restriction of $\h_i$ to $\mathcal D_i^+$ and $\h_i^-$ the restriction to $\mathcal D_i^-$.

\begin{proposition} \label{p:recTh}
Suppose that there exists $\sigma \in \CZ$ and a pair of holomorphic maps $\psi = (\psi^+, \psi^-)$, where $\psi^\pm: \Dd_1^\pm \rightarrow \Dd_2^\pm$ such that: 
\[ \h_1^\pm = T_\sigma^{-1} \circ \h_2^\pm \circ \psi^\pm \]
on $\Dd_1^\pm$.
Then there exists a local semi-conjugacy $\phi$ on immediate basins from $f_1$ to $f_2$ (see \Cref{d:semiConjLoc}).
\end{proposition}

\begin{compl} \label{compl:recTh}
Suppose furthermore that  $\psi^+$ and $\psi^-$ are biholomorphisms $\Dd_1^\pm \to \Dd_2^\pm$.
Then there exists a local pseudo-conjugacy $(\phi, \phi')$ of $f_1, f_2$ (see \Cref{d:pseudoConjLoc}).
\end{compl}

As noted in \Cref{p:exPSSync}, we may choose this local pseudo-conjugacy to be synchronous.

\begin{figure}
\def\svgwidth{\textwidth}
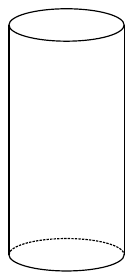

\vspace{-3.5cm}

\caption{Proof of \Cref{p:recTh} and \Cref{compl:recTh}}
\label{f:demRecTheoreme}

\vspace{0.5cm}
\footnotesize{\textit{The sets $P_A^1$, $P_A^2$ are represented by cardioids, $P_R^1$, $P_R^2$ by circles (in the proof we will instead take $\alpha$-petals with $\alpha>2$ for $P_A^1$, $P_R^1$ and $P_A^2$).
$U_1 = P_A^1 \cup P_R^1 \cup \{0\}$.
We denote $(P_R^1)^0 = (\Phi_R^1)^{-1}(\wDd_1^+ \cup \wDd_1^-) = U_1^0\cap P_R^1$.
Centre left, in dark gray: $U_1^0$, in light gray: $U_1 \cap B^0_{f_1} \setminus U_1^0$.
We represent also in dark gray, light gray their diverse images by the maps of the diagram.
The map $\phi$ will be built as the unique function completing this commutative diagram.
Note that a priori, $\phi(U_1^0)$ is not a set of the form $U_2^0$ for some open neighborhood $U_2$ of $0$, and that $Y_R^2:= \Psi_R^{2,\ext}\circ \wt\psi\circ \Phi_R^1((P_R^1)^0)) = $ (the candidate for $\phi(P_R^1)$) is not a set of the form $(P_R^2)^0$ for some repelling petal $P_R^2$ of $f_2$.}}
\end{figure}

\medskip

Before giving the details of the proof of \Cref{p:recTh} in \Cref{subsub:pf:p:recTh} and of \Cref{compl:recTh} in \Cref{subsub:pf:compl:recTh}, we give a brief summary here.
We first lift $\wt\sigma\in\CZ$ into $\sigma\in\C$ and
the maps $\psi^\pm$ and into maps $\wt\psi^\pm$ in a way that
the following diagram commutes:
\[
\xymatrix{
  \wDd_1^\pm \ar[d]_{\wt\psi^\pm} \ar[r]^{\Psi_R^{1,\ext}} & B_{f_1}^0 \ar@{..>}[d]_{\phi} \ar[r]^{\Phi_A^{1,\ext}} & \C \ar[d]^{T_{\tilde \sigma}}
  \\
  \wDd_2^\pm \ar[r]_{\Psi_R^{2,\ext}} & B_{f_2}^0 \ar[r]_{\Phi_A^{2,\ext}} & \C
}
\]
The objective is to use both sides of the diagram to define a map $\phi$ that completes the diagram and then to check that it is a local semi-conjugacy on immediate basins.
It will not be possible to define $\phi$ on the whole immediate basin $B_{f_1}^0$ (see \Cref{sub:ex:3}) but only on a set of the form $U_1^0$, and for a small enough neighborhood $U_1$ of $0$ (so that we can apply \Cref{p:petitPetRep2}).
We then prove that the constructed $\phi$ has for lifted phase $\wt\sigma$.
In the case where the two maps $\psi^\pm$ are biholomorphisms for $\Dd_1^\pm$ to $\Dd_2^\pm$, this allows us to check that the two local semi-conjugacies maps on immediate basins $\phi$ from $f_1$ to $f_2$ (built from $\wt\psi^\pm$) and $\phi'$ from $f_2$ to $f_1$ (built from the inverses $(\wt \psi^\pm)^{-1}$) have opposite phases, and thus form a pseudo-conjugacy.
See \Cref{f:demRecTheoreme}.

\subsubsection{Construction of a local semi-conjugacy on immediate basins}\label{subsub:pf:p:recTh}

We show here \Cref{p:recTh}.
Let $T_\sigma$ and $\psi^\pm$ satisfy the hypothesis of the proposition. 

\medskip

\underline{\textit{Choice of adapted lifts $\wt \sigma$ and $\widetilde{\psi}^\pm$ of $\sigma$ and $\psi^\pm$ by $\piZ$: }} 
Recall we denote by $\piZ: \C \to \CZ$ the canonical projection. The quantity $\sigma$ is only defined modulo $\Z$.
We choose any $\wt\sigma$ such that $\piZ(\wt\sigma) =\sigma$.
Let $\wDd_i = \piZ^{-1}(\Dd_i)$, $\wDd_i^\pm = \piZ^{-1}(\Dd_i^\pm)$ for $i =1, 2$.

Let $\widetilde{\psi}^\pm: \wDd_1^\pm \to \wDd_2^\pm$ be a (continuous) lift of $\psi^\pm$ by $\piZ$.
By connectedness of $\Dom \psi^\pm = \wDd_1^\pm$, choosing another lift of $\wt\psi^\pm$ amounts to adding to it a constant in $\Z$ (not necessarily the same for the top and the bottom domains).
The hypothesis, together with the relation $\h_i \circ \piZ = \piZ \circ h_i$ imply that
\[ \forall w \in \wDd_1^\pm,\ T_{\tilde\sigma} \circ h_1 (w) \equiv h_2 \circ \widetilde{\psi}^\pm(w) \bmod \Z\]
By connectedness of $\wDd_1^\pm$ and continuity, the difference $h_1 - T_{\tilde\sigma}^{-1} \circ h_2 \circ \widetilde{\psi}^\pm$ takes over $\wDd_1^+$ and $\wDd_1^-$ two (possibly equal) constant values, which belong to $\Z$.
We can ensure that this difference is $0$ by subtracting these constant to $\wt\psi^+$ and $\wt\psi^-$: indeed $T_{\tilde\sigma}^{-1}$ and $h_2$ commute with $T_1$.
This also works if $\wDd_1^+ = \wDd_1^-$.

We choose in the sequel these lifts of $\wt\psi^\pm$.
They depend on our choice of $\wt \sigma \in \piZ^{-1}(\{\sigma\})$.
Recalling the definition of $h_i$
we thus have:
\begin{equation}\label{eq:pf1}
  \forall w\in\wDd_1^\pm,\ T_{\tilde \sigma} \circ \Phi_A^{1,\ext} \circ \Psi_R^{1,\ext} (w) =  \Phi_A^{2,\ext} \circ \Psi_R^{2,\ext} \circ \widetilde{\psi}^\pm(w)
\end{equation}

\medskip

\textit{\underline{Choice of the petals and domain $U_1$}}

We choose the petals of $f_1$ and the domain $U_1$ as follows, which is a special case of the conditions of \Cref{ssub:ci}, so that \Cref{p:petitPetRep2} can apply.
More precisely we let $\alpha = \beta \in (\pi/2,\pi)$ (so $\alpha+\beta>\pi$),
$P_A^1$ is an attracting $\alpha$-petal of $f_1$,
$P_R^1$ is a repelling $\alpha$-petal of $f_1$
such that $P_A^1 \cap P_R^1$ has two connected components $C_+, C_-$.
In this case we necessarily\footnote{This does not rely on \Cref{p:petitPetRep} but on a simpler argument explained in \Cref{ssub:ci}.} have
$\piZ\circ\Phi_R^1(C_\pm) \subset \Dd_1^\pm$.
Finally, $U_1$ is initially defined as 
\[U_1 = P_A^1 \cup P_R^1 \cup \{0\}\]
(we will shrink it a couple of times).

Throughout the rest of this section we denote for $i=1$ and $i=2$:
\[Q_A^i=\Phi_A^i(P_A^i),\quad Q_R^i=\Phi_R^i(P_R^i)\]

Choose for the set $P_A^2$ an $\alpha$-petal of $f_2$ such that $Q_A^2\subset T_{\tilde \sigma}(Q_A^1)$.
Then, by decreasing $P_1$ (hence $U_1$), we can ensure that
\[Q_A^2=T_{\tilde \sigma}(Q_A^1)\]
without destroying the previous conditions.

Denote $\Phi_A^i = \Phi_A^{i,\ext}|_{P_A^i}$ and $\Phi_R^i$ the unique repelling Fatou coordinate on $P_R^i$ that is an inverse branch of $\Psi_R^{i,\ext}$ (see \Cref{cor:invPsionPrep}).
\begin{lemma}\label{lem:compatible}
  For $z \in P_A^1 \cap P_R^1$, the two sides of the following equality are defined and the equality holds:
  \begin{equation}\label{eq:e3}
    \forall z \in P_A^1\cap P_R^1,\ (\Phi_A^2)^{-1} \circ T_{\tilde \sigma} \circ \Phi_A^1(z) = \Psi_R^{2,  \ext } \circ \widetilde{\psi} \circ \Phi_R^1(z)
  \end{equation}
\end{lemma} 
\begin{proof}
The left member is defined over $P_A^1$ since $T_{\tilde \sigma}(Q_A^1) = Q_A^2$.
The right member is defined over $P_A^1 \cap P_R^1$ since $\Phi_R^1(P_A^1 \cap P_R^1) = \Phi_R^1(C^+) \cup \Phi_R^1(C^-) \subset \Dd_1^+ \cup \Dd_1^-$.
  
By applying \cref{eq:pf1} to $w=\Phi_R^1(z)$ with $z\in P_A^1\cap P_R^1 = C^+\cup C^-$, we get
\begin{equation}\label{eq:e4}
  \forall z \in P_A^1\cap P_R^1,\ T_{\tilde \sigma} \circ \Phi_A^{1,\ext} (z) = \Phi_A^{2,\ext} \circ \Psi_R^{2,\ext} \circ \widetilde{\psi}^\pm \circ \Phi_R^1(z)
\end{equation}
and we would like to apply $(\Phi_A^2)^{-1}$ to this equation but we do not know if the right hand side belongs to $P_A^2$.

Instead, we are going to use an analytic extension argument.
In this view, it is enough to prove the equality on any subset of $V\subset C^+$ having an accumulation point in $C^+$ (and similarly for $C^-$). A non-empty open set will do.
For $V$, we take a small sector bisected by an axis perpendicular to the attracting and repelling axes, of half-opening angle $<\alpha-\pi/2$ and contained in $C^+$.
By \Cref{cor:sec1} such a sector contained in $C^+$ exists.
By \Cref{lem:i2}, $\widetilde{\psi}^\pm(w) = w + \rho^\pm + o(1)$
for some $\rho^\pm\in\C$.
If $z$ tends to $0$ within $V$ we have by this estimate and by \Cref{lem:asymp1}:
$\wt\psi^+(\Phi_R^1(z))\sim -1/a_1z$ and by \Cref{lem:asymp2}: $\Psi_R^{2,\ext}(\wt\psi^+(\Phi_R^1(z)))\sim \frac{a_1}{a_2}z$.
Since $P_A^2$ is an $\alpha$-petal, \Cref{cor:sec1} implies that it will contain $z':=\Psi_R^{2,\ext}(\wt\psi^+(\Phi_R^1(z))$ provided $z\in V$ is small enough.
So $\Phi_A^{2,\ext}(z') = \Phi_A^2(z')$, so by \cref{eq:e4}, \cref{eq:e3} holds on $V$.
As we explained above, the lemma follows.
\end{proof}

Recall that $U_1 = P_A^1\cup P_R^1 \cup\{0\}$.
Recall that $U_1^0$ denotes the connected component of $U_1\cap B_{f_1}$ that contains an attracting petal.
Let us check that the map $\phi: U_1^0 \to \C$ defined below is well defined, holomorphic, and takes values in $B_{f_2}^0$:
\[
\phi(z) =
\begin{cases}
  (\Phi_A^2)^{-1} \circ T_{\tilde \sigma} \circ \Phi_A^1(z) & \text{ if }z \in P_A^1 \\[4pt]
  \Psi_R^{2,  \ext  } \circ \widetilde{\psi} \circ \Phi_R^1(z) & \text{ if }z \in P_R^1 \cap U_1^0
\end{cases}
\]
The first line is defined as we already saw (it is a bijection from $P_A^1$ to $P_A^2$, though we do not need this fact in this section).
By \Cref{p:petitPetRep2}, $\Phi_R(U_1^0) \subset \wDd^+\cup \wDd^-$ so the expression in the second line is defined.
According to \Cref{lem:compatible}, these two definitions coincide on the intersection of the two domains.
Each line is holomorphic and defined on an open domain, so the conjunction of the two is. 
The first line takes values in $P_A^2\subset B_{f_2}$.
The map $\wt\psi^\pm$ takes values in $\wDd_2^\pm$ hence in $\wDd_2$, so the second line takes values in $B_{f_2}$.
Since the domain of $\phi$ is connected and $\phi$ continuous, its image is connected. This image is contained in $B_{f_2}$ and contains $P_A^2$, so it is contained in $B_{f_2}^0$.

There remains to check the semi-conjugacy relation on $U_1^0\cap f^{-1}(U_1^0)$
\[f_2\circ \phi = \phi\circ f_1\]
for this we will need to shrink the set $U_1$ (see below).

But first, the relation is easily checked on the subset $P_A^1$ of $U_1^0\cap f^{-1}(U_1^0)$, as $z\in P_A^1$ implies that $f_1(z) \in P_A^1$ and by the first line, $\phi(f(z)) =   (\Phi_A^2)^{-1} \circ T_{\tilde \sigma} \circ \Phi_A^1(f(z)) = (\Phi_A^2)^{-1} \circ T_{\tilde \sigma} \circ T_1\circ \Phi_A^1(z) = (\Phi_A^2)^{-1} \circ T_1 \circ T_{\tilde \sigma} \circ \Phi_A^1(z) = f_2\circ (\Phi_A^2)^{-1} \circ T_{\tilde \sigma} \circ \Phi_A^1(z)$ where for the last equality we used that $w:=T_{\tilde \sigma} \circ \Phi_A^1(z)$ belongs to $Q_A^2$ since by construction $Q_A^2=T_{\tilde \sigma}(Q_A^1)$.

Consider the subset ${P'_R}^1$ of $P_R^1$ defined by $\Phi_R^1({P'_R}^1) = T_{-1}(\Phi_R^1(P_R^1))$. 
Like $P_R^1$, the set ${P'_R}^1$ is a repelling $\alpha$-petal.
Let us check the semi-conjugacy relation on $P'_R\cap U_1^0\cap f^{-1}(U_1^0)$.
We have $f_1(z)\in P_R\cap U_1^0$ and by the second line of the definition of $\phi$, $\phi(f_1(z)) = \Psi_R^{2, \ext} \circ \widetilde{\psi} \circ \Phi_R^1\circ f_1(z) = \Psi_R^{2, \ext} \circ \widetilde{\psi} \circ T_1 \circ \Phi_R^1(z) = \Psi_R^{2, \ext} \circ T_1 \circ \widetilde{\psi} \circ \Phi_R^1(z) = f_2\circ \Psi_R^{2, \ext} \circ \widetilde{\psi} \circ \Phi_R^1(z)$.

In the last two paragraphs, we checked the semi-conjugacy relation on the two following subsets of $U_1^0\cap f_1^{-1}(U_1^0)$: $P_A^1$ and $P'_R\cap U_1^0\cap f^{-1}(U_1^0)$.
Now instead of trying to check the relation on the rest of $U_1^0\cap f_1^{-1}(U_1^0)$, it is enough to replace $U_1=P_A^1\cup P_R^1\cup\{0\}$ by a subset $U'_1 = P_A^1\cup {P'_R}^1\cup\{0\}$.
Then ${U'_1}^0\cap f^{-1}({U'_1}^0)\subset P_A^1\cup (P'_R\cap U_1^0\cap f^{-1}(U_1^0))$ and the restriction of $\phi$ to $U'_1$ is a local semi-conjugacy on immediate basins in the sense of \Cref{d:semiConjLoc}.

This completes the proof of \Cref{p:recTh}.

\begin{remark}
  Let $(P_R^1)^0 = (\Phi_R^1)^{-1}(\wDd_1^+ \cup \wDd_1^-) = P_R\cap U_1^0$ by \Cref{p:petitPetRep2}.
  Since $\widetilde{\psi}$ does not admit a continuous extension to the boundary of $\wDd_1^+ \cup \wDd_1^-$, it seems not obvious that the set $Y_R^2:=\Psi_R^{2,\ext}\circ\wt\psi\circ \Phi_R^1((P_R^1)^0)$ is included in a set of the form $(P_R^2)^0 = (\Phi_R^2)^{-1}(\wDd_2^+ \cup \wDd_2^-)$ where $P_R^2$ is a repelling petal of $f_2$. \Cref{f:demRecTheoreme} shows how a counter-example might look like. That is why we used the extended repelling parametrization $\Psi_R^{2,\ext}$ in the definition of $\phi$, since it seems difficult to choose a repelling petal $P_R^2$ for which $(\Phi_R^2)^{-1}$ could be used. 
\end{remark}

Recall that the phase $\sigma(\phi)$ and lifted phase $\wt\sigma(\phi)$ of a local semi-conjugacy on immediate basins were defined in \Cref{def:phase} and satisfy (see \Cref{eq:cd}) $\sigma(\phi)=\piZ(\wt\sigma(\phi))$ where 
\[\Phi_A^{2,\ext} \circ \phi = T_{\tilde\sigma(\phi)} \circ \Phi_A^{1,\ext}\text{ on $P_A^1$.}\]
The following property of the map $\phi$ we constructed above immediately follows from the construction and will be used in the proof of the complement:

\begin{lemma}\label{lem:phase}
  The lifted phase $\wt\sigma(\phi)$ of the local semi-conjugacy on immediate basins $\phi$ constructed above is equal to $\wt\sigma$.
\end{lemma}

\subsubsection{Pseudo-conjugacy}\label{subsub:pf:compl:recTh}

We now show \Cref{compl:recTh}.
We suppose that $\psi^\pm$ are biholomorphisms from $\Dd_1^\pm$ to $\Dd_2^\pm$.

We use $\psi^\pm$ to construct as above a local semi-conjugacy $\phi$ on immediate basins from $f_1$ to $f_2$.
We also use $(\psi^\pm)^{-1}$ to construct a similar map $\phi'$, but from $f_2$ to $f_1$.
By \Cref{prop:pciffssz}, to prove that $(\phi,\phi')$ to be a pseudo-conjugacy, it is enough to check that the respective phases satisfy $\sigma(\phi) + \sigma(\phi') = \ov 0$ in $\CZ$.
For this we use \Cref{lem:phase}.
We get $\sigma(\phi)=\piZ(\sigma)$ while $\sigma(\phi')=\piZ(-\sigma)$ since from
\[ T_\sigma \circ \h_1^\pm = \h_2^\pm \circ \psi^\pm \]
we get
\[  \h_1^\pm \circ (\psi^\pm)^{-1} = T_{-\sigma} \circ \h_2^\pm\]

\begin{remark}
  Another approach is to build simultaneously $\phi$ and $\phi'$, see \cite{t:LM}, pages~92 to~96.
\end{remark}

This completes the proof of \Cref{compl:recTh}.


\section{Examples}\label{sec:ex}

%
%

\subsection{Wild semi-conjugacies on petals}\label{sub:ex:5}

Let $f_1$, $f_2$ be two holomorphic maps defined on open subsets of $\C$ containing $0$ and such that $f_i(z)=z+a_iz^2+\ldots$ with $a_i\in\C^*$.
Recall that, to an attracting petal $P_1$ for $f_1$ and a semi-conjugacy $\phi:P_1\to B_{f_2}$ the \Cref{prop:inter} associates a uniquely defined map $\wt\phi:\C\to\C$.

\begin{lemma}\label{lem:ce5}
  Let $g:\C\to\C$ be a holomorphic map commuting with $T_1$ and $P_2$ a petal of $f_2$.
  Then there exists an attracting petal $P_1$ of $f_1$ and a semi-conjugacy $\phi:P_1\to P_2$ from $f_1$ to $f_2$ such that $\wt\phi=g$.
\end{lemma}

If one takes any $\phi$ that is not a translation (for instance: $\wt\phi(w) = w + \exp(2\pi iw)$) this allows to build a counter-example to \Cref{p:SCQuotientBiholom} if we drop the condition that $P_1$ must contain an $\alpha$-petal.

\medskip

\begin{proof}[Proof of \Cref{lem:ce5}]
The map $g$ descends by $\piZ$ to a self-mapping $\wh g$ of $\CZ$ that is necessarily surjective (by the generalized Liouville theorem), but not necessarily injective.
Consider any petal $P'_1$ for $f$.
Denote $Q'_1=\Phi_A^1(P'_1)$ and $Q_2=\Phi_A^2(P_2)$.
These sets are $T_1$-stable (sent in themselves by $T_1$).
Consider the set $Q''_1 = Q'_1\cap g^{-1}(Q_2)$.
It is $T_1$-stable too.
By the maximum principle $Q''_1$ has all its connected components simply connected.
Also we have $\piZ(T_1(Q''_1))=\CZ$ (use that $\piZ(Q''_1)=\CZ$, $\piZ(Q_2)=\CZ$, $g$ commutes with $T_1$ and the $T_1$-stability of $Q''_1$ and $Q_2$).
Since $\R/\Z$ is compact there must exist $n\in\N$ such that $[n_0,n_0+1]\subset Q''_1$.
The connected component $Q_1$ of $Q''_1$ that contains $[n_0,n_0+1]$ satisfies $T_1(Q_1) \subset Q_1$.
Moreover, $\piZ(Q_1)=\CZ$.
Indeed by the same argument as above, $\forall z$, the segment $[n,n+z]$ is contained in $Q''_1$ for all $n$ big enough.
Then the set $P_A^1 := (\Phi_A^1)^{-1}(Q_1)$ is a petal for $f_1$ and the map $\phi= (\Phi_A^2)^{-1} \circ g|_{Q_1} \circ \Phi_A^1$ is defined on $P_A^1$, takes values in $P_A^2$, is a semi-conjugacy from $f_1$ to $f_2$ and $\wt\phi=g$.
\end{proof}

\subsection{A counterexample to a converse inclusion for \Cref{p:petitPetRep}}\label{sub:ex:1}

\begin{figure}[!t]
  \begin{center}
    \begin{tikzpicture}
      \node at (0,0.2) {\includegraphics[scale=0.5]{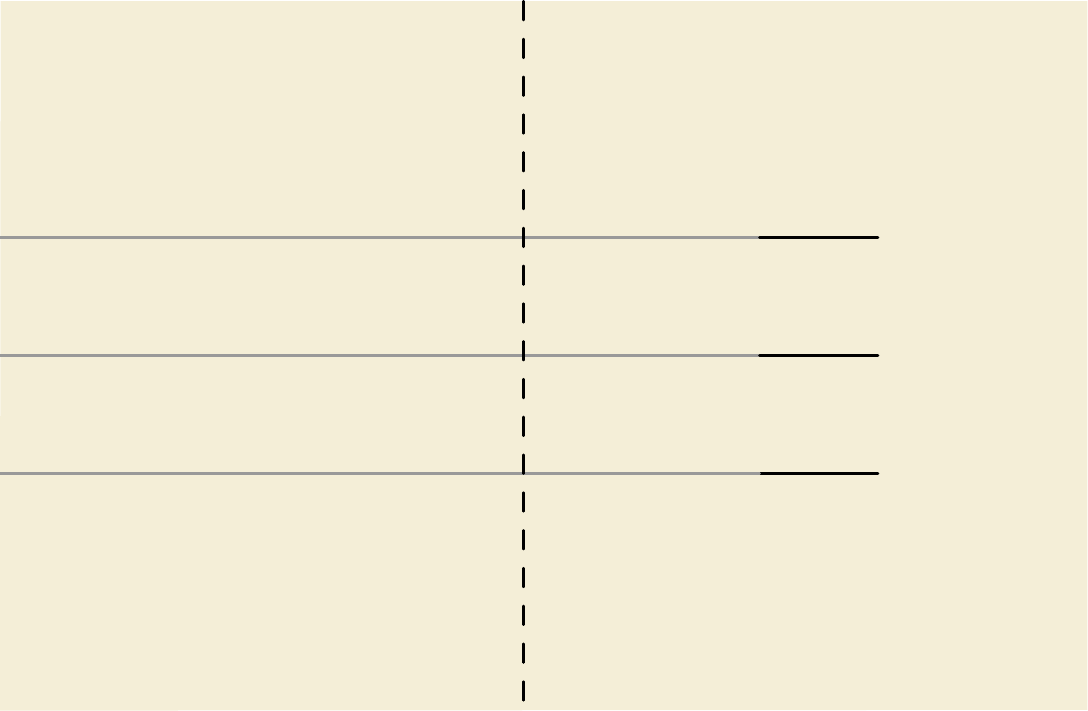}};
      \node at (6.5,0.2) {\includegraphics[scale=0.5]{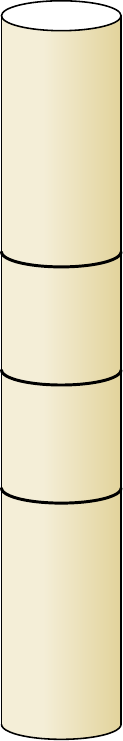}};
      \node at (3.1,0.2) {$1$};
      \node at (1.85,-0.1) {$0$};
      \node at (-0.5,-2.5) {$P_R$};
      \node at (6.6,2.3) {$\Dd^+$};
      \node at (6.6,-2.1) {$\Dd^-$};
    \end{tikzpicture}
  \end{center}
  \caption{The example of \Cref{sub:ex:1} is the restriction of $z\mapsto z+1$ to the complement of the countable union of the three thick black segments.
    The basin of the parabolic point $\infty$ is the complement of the union of black and gray segments (this union is also the union of the translates of the black set by non-positive integer vectors).
    The immediate basin in in yellow.
    The chosen repelling petal is the left half plane delimited by the dashed line.
    On the right: the image by $\Phi_R:P_R\to\CZ$ of the immediate basin is the cylinder minus three horizontal circles.}
  \label{fig:ex:2}
\end{figure}

The converse inclusions in \Cref{p:petitPetRep} are not true in full generality.
Consider the following example, where we work in coordinate $u=s(z):=-1/z$ that puts $0$ at $\infty$ and vice-versa.
The map $z\mapsto f(z)$ corresponds to $u\mapsto g(u)$ with $g=s\circ f\circ s^{-1}$. 
Because $f$ is required to be defined in a subset of $\C$ and to map to $\C$, the point $u=0$ cannot be in the domain nor in the range of the map $g$.
Let $g(u)=u+1$ that we restrict to
$\Dom  g := \hat \C-\Xi$ where
\[\Xi=T_{\bi}([0,1])\cup[0,1]\cup T_{-\bi}([0,1])\]
for some $\eps>0$ and $\bi$ denotes the imaginary unit.
Then $B_g = \C-T_{\bi}((-\infty,1])\cup(-\infty,1]\cup T_{-\bi}((-\infty,1])$ is connected so $B_g = B_g^0$.
A repelling petal for $g$ is for instance given by the left half plane $\Re(u)<-2$ and $\Phi_R(u)=u$.
So $\Phi_R(B^0_{g}) = B^0_{g} \cap P_R$ so
$\piZ \circ \Phi_R(B^0_{g})$ has 4 connected components: it is the cylinder cut by three horizontal lines.
The sets $\Dd^+$ and $\Dd^+$ are the unbounded ones.
If we take $U=\hat \C-\{0\}$ we get that $U^0=B^0_{g}$ and $\Wt{U}$ is $\C$ minus three horizontal lines, hence
\[ \piZ \circ \Phi_R(U^0) = \W{U}  \not\subset \Dd^+ \cup \Dd^-\]

\subsection{Immediate basin with a trace in repelling Fatou coordinate that is not stable by $T_{-1}$}\label{sub:ex:2}

\begin{figure}[!t]
  \begin{center}
    \begin{tikzpicture}
      \node at (0,0.2) {\includegraphics[scale=0.5]{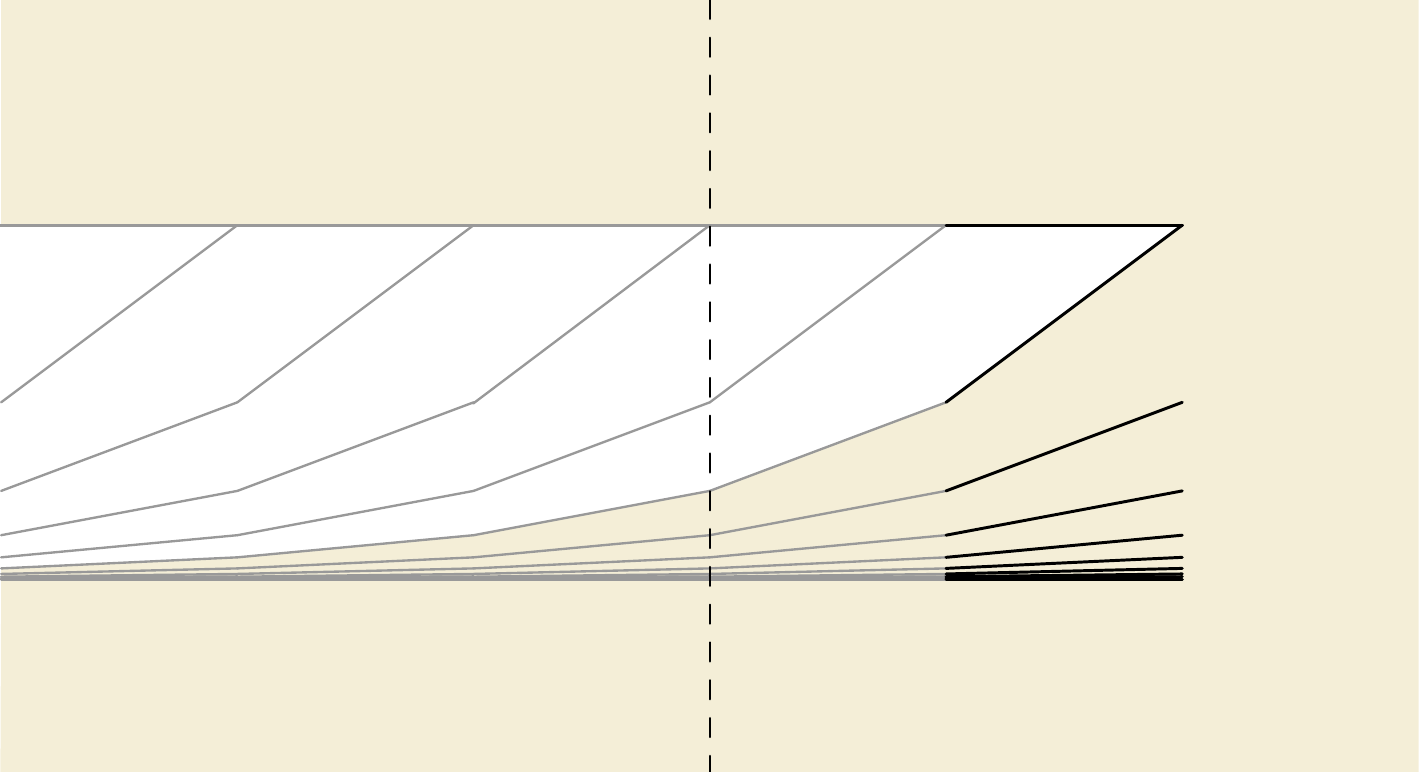}};
      \node at (1.95,-1.8) {$-1$};
      \node at (4,-1.8) {$0$};
      \node at (4.3,1.7) {$i t_1$};
      \node at (4.3,0.1) {$i t_2$};
      \node at (2.1,1.85) {$i t_1-1$};
      \node at (-0.35,-2.75) {$P_R$};
    \end{tikzpicture}
  \end{center}
  \caption{The example of \Cref{sub:ex:2} is the restriction of $z\mapsto z+1$ to the complement of the countable union of solid black segments.
    The basin of the parabolic point $\infty$ is the complement of the union of black and gray segments (this union is also the union of the translates of the black set by non-positive integer vectors).
    The immediate basin in in yellow and the rest of the basin in white.
    The chosen repelling petal is the left half plane delimited by the dashed line.}
  \label{fig:ex:2}
\end{figure}

The image of $B_f^0$ in by the repelling Fatou coordinate $\Phi_R$ of a repelling petal may fail to be stable (sent in itself) by $T_{-1}$.
To give an example, illustrated in \Cref{fig:ex:2}, we take the same construction as in \Cref{sub:ex:1} where we replace the compact set $\Xi$ by a compact set $K$ that we define as follows.
Consider the sequence of imaginary numbers $i t_n$ for $n>0$, where $t_n>0$ is a real sequence decreasing to $0$ and $i$ is the complex unit.
We let $K$ be the union of the real interval $[-1,0]$, the horizontal segment $[it_1-1,it_1]$ above it, and of the segments $[it_{n+1}-1,it_n]$ for all $n\geq 1$.
As previously we let $u=-1/z$, $g(u)=u+1$ restricted to $\{\infty\}\cup\C-K$, $P_R = ``\Re(u)<-2"$ and $\Phi_R(u)=u$.
The basin $B_g$ is the union of the three half-planes $\Re(u)>0$, $\Im(u)>1$, $\Im(u)<0$ and of countably many infinite tongues contained in the complement of the above planes, and of the form ($\Re(u)\leq 0$ and $h_{n+1}(\Re(u))<\Im(u)<h_n(\Re(u))$) where $h_n:(-\infty,0]\to\R$ is affine on each interval $[-k-1,-k]$, $k\in\N$, and $h_n(-k) = -k + it_{n+k}$.
In particular $B_g$ is connected, i.e.\ $B_g^0 = B_g$.
The set $\Phi_R(B_g^0) = P_R \cap B_g^0$ is not connected.
If we take $U=\wh\C -\{0\}$, then $W(U) = \Dd^+\cup\Dd^- \not\supset \piZ\circ \Phi_R(B_g^0)$.

\subsection{On extending local semi-conjugacies}\label{sub:ex:3}

\subsubsection{}\label{sub:ex:3:sub:1}

We still work in coordinates $u=-1/z$.
Let us take $g_1=g_2$ to be the second example above.
Then $\wDd_1^+=\wDd_2^+$ is the half plane of equation $\Re(u)>t_1$, and $\wDd_1^- = \wDd_2^-$ is the one of equation $\Re(u)<0$.
Let $\wt\sigma \in\R-\Z$ and $\wt\psi^\pm(u) = u + \wt\sigma$.
Then we can define the local semi-conjugacy $\phi(z)=z+\wt\sigma$ on a quite large set: we can take for instance $U=\wh\C-\{0\}$ and $\Dom \phi = $ the union of the half planes $\Re(u)>\max(0,-\wt\sigma)$, $\Im(u)>t_1$, $\Im(u)<0$.
However, $\phi$ cannot be extended to a map $B_{g_1}^0\to B_{g_2}^0$.

\subsubsection{}\label{sub:ex:3:sub:2}

Consider the cauliflower map $f(z) = z+z^2$.
It is known that $B_f = B_0^f$ and that $B_f$ is a Jordan domain with fractal boundary and that $\Dd=\Dd^+\cup\Dd^-$ is bounded by two Jordan curves (see \cite{b:LY}).
A local branch of $f^{-1}$ restricts to an injective local semi-conjugacy $\phi$ on immediate basins from $f$ to itself, with lifted phase $\wt\sigma(\phi)=-1$ and phase $\sigma(\phi)=\ov 0$.
It is impossible to extend $\phi$ it to the whole immediate basin.
The associated map $\psi$ is particularly simple: it is the identity on $\Dd^+\cup \Dd^-$.
The map $\wt\psi$ associated to $\wt\sigma=-1$ is the translation by $-1$ on $\wDd^+\cup \wDd^-$.
But things can be fixed:
taking another lift of $\sigma$, namely $\wt\sigma =0$, yields $\wt\psi$ being the identity on $\wDd^+\cup \wDd^-$ and is associated to the map $\phi=\on{id}$ on $U_1^0$, which obviouslty extends to a conjugacy from $f$ to itself (see also \cite{b:DavidMorris}). 

\subsubsection{}\label{sub:ex:3:sub:3}

In the example of \Cref{sub:ex:3:sub:1}, $f_1=f_2$ was an artificial restriction of a very simple map.
In the example of \Cref{sub:ex:3:sub:2}, $f_1=f_2$ is a global map but the non-extensible local semi-conjugacy $\phi$ on immediate basins can be post-composed by $f_2$ to get an extensible one.
A more interesting example, where at the same time the $f_i$ not be an artificial restrictions, but where a local semi-conjugacy on immediate basins would not extend even if we post-compose it by a fixed iterate of $f_2$, would necessarily be complicated according to the main theorem (and methods) of \cite{b:DavidMorris}.
We believe this can be done but this is out of scope of the present article.

\begin{figure}[!t]
  \begin{center}
    \begin{tikzpicture}
    \node at (0,0) {\includegraphics[width=10cm]{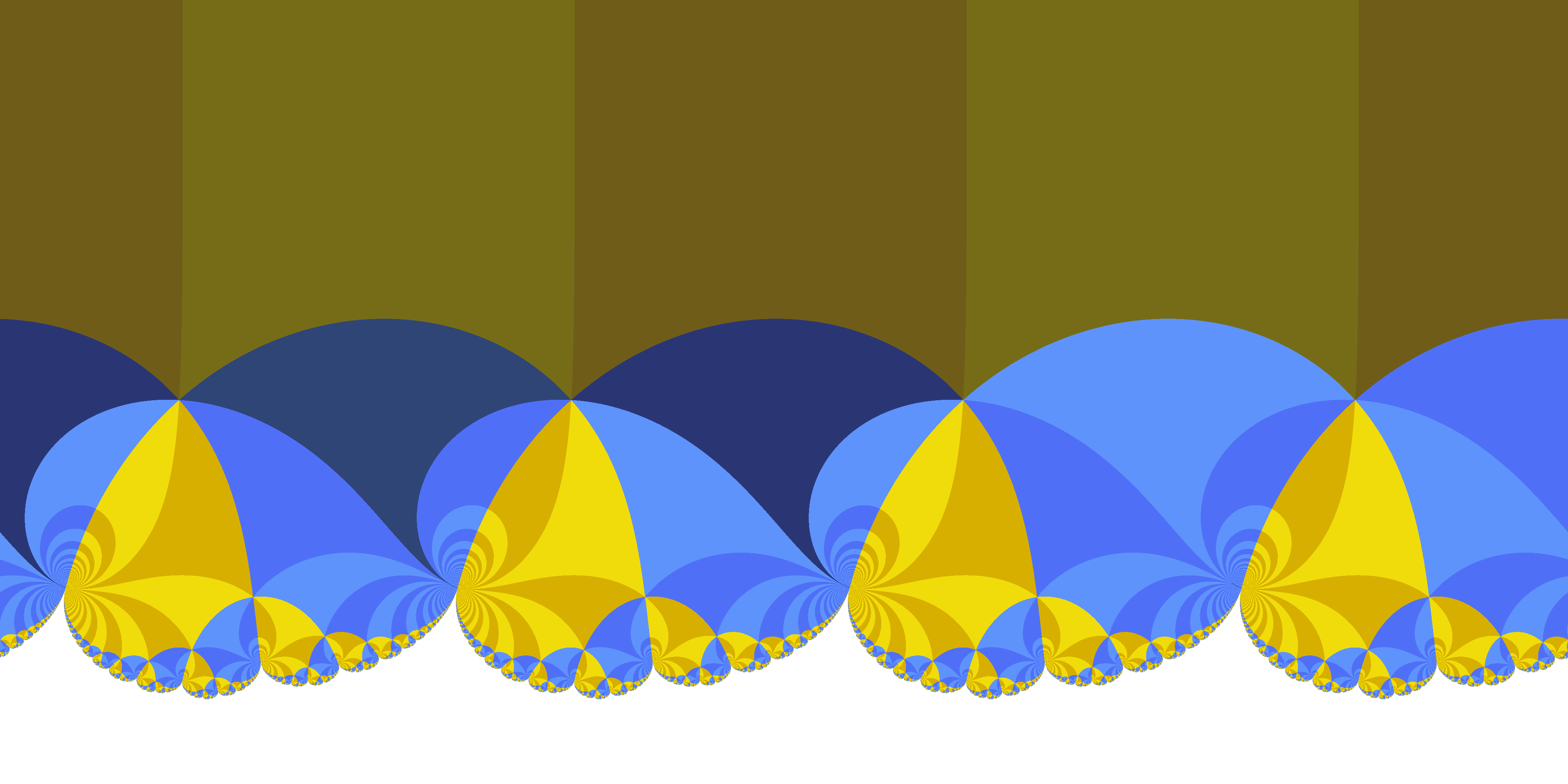}};
    \node at (0,-7) {\includegraphics[width=6cm]{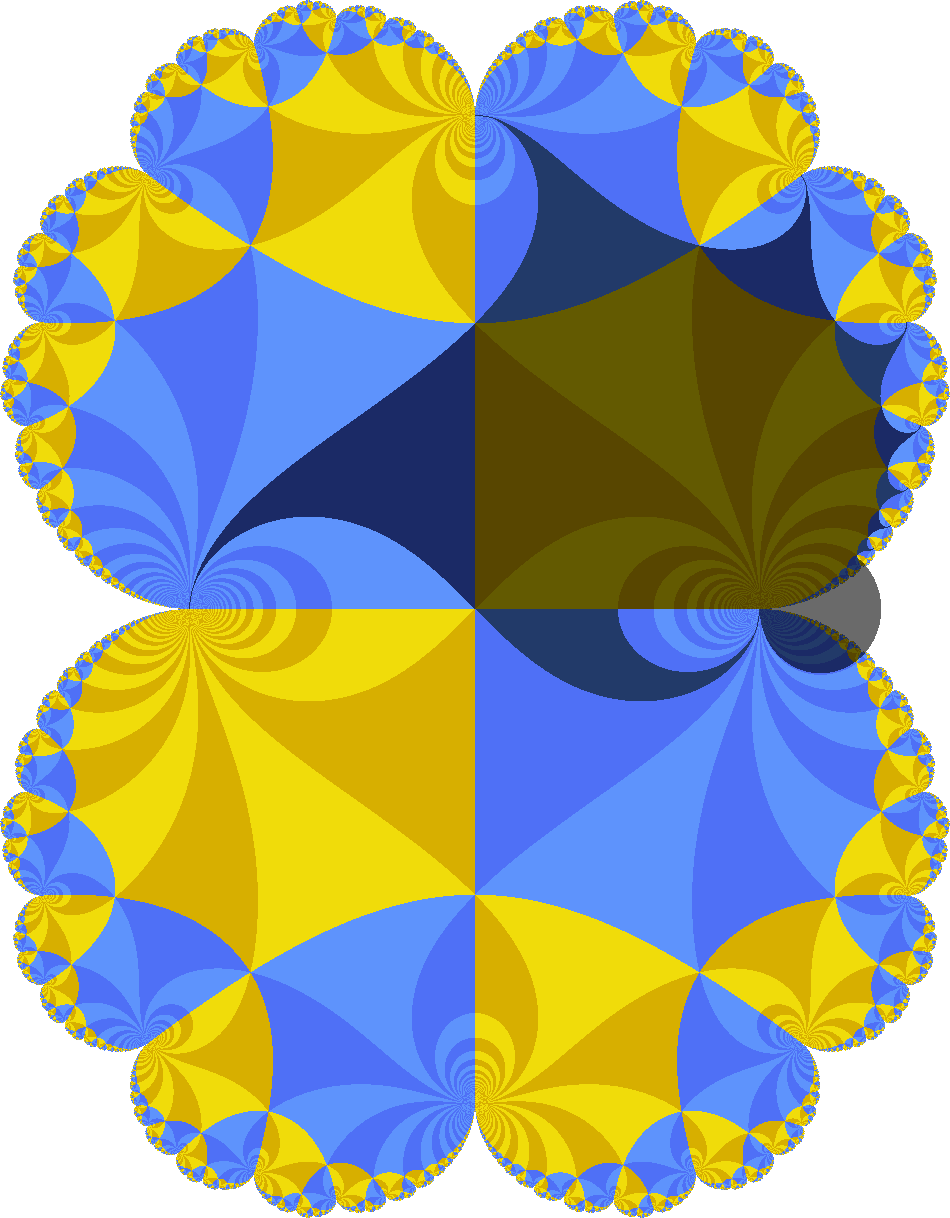}};
    \draw[->] (1.75,-7.55)[out=-65, in=-55] to (-1.55,-6.95);
    \node[circle, fill=white, draw=black] at (0,1.5) {$V$};
    \node[rectangle, fill=white, draw=black] at (1.3,-6) {$\Psi_R^\ext(V)$};
    \draw[->] (0,-2) -- node[right]{$\Psi_R^\ext$} (0,-3);
    \end{tikzpicture}
  \end{center}
  \caption{
    Illustration of the example of \Cref{ss:ex:4:1}.
    The attracting Fatou coordinate is normalized by $\Psi_A^\ext(c)=0$ where $c$ is the critical point.
    \\
    \textit{Bottom:} each point $z$ in the parabolic basin of the cauliflower map $f:z\mapsto z+z^2$ is colored in blue or yellow, according to whether $w=\Psi_A^\ext(z)$ has positive or negative imaginary part. There are two shades of blue and of yellow, to indicate the parity of the entire part of the real part of $w$ (i.e.\ alternate strips of width $1$).
    \\
    \textit{Top:} preimage (ignoring the distinction between bright and dark shades) of the bottom picture by the extended inverse repelling Fatou coordinate $\Psi_R^\ext$.
    \\
    \textit{In dark shades:} 
    the set $V$ on top and $\Psi_R^\ext(V)$ at bottom.
    The shown set $V$ is completed by taking the union with a left half plane far on the left (not visible here), whose image is the circular-like shape.
  }
  \label{fig:ex:4}
\end{figure}

\begin{figure}
  \begin{center}
  \begin{tikzpicture}
    \node at (0,0) {\includegraphics[width=10cm]{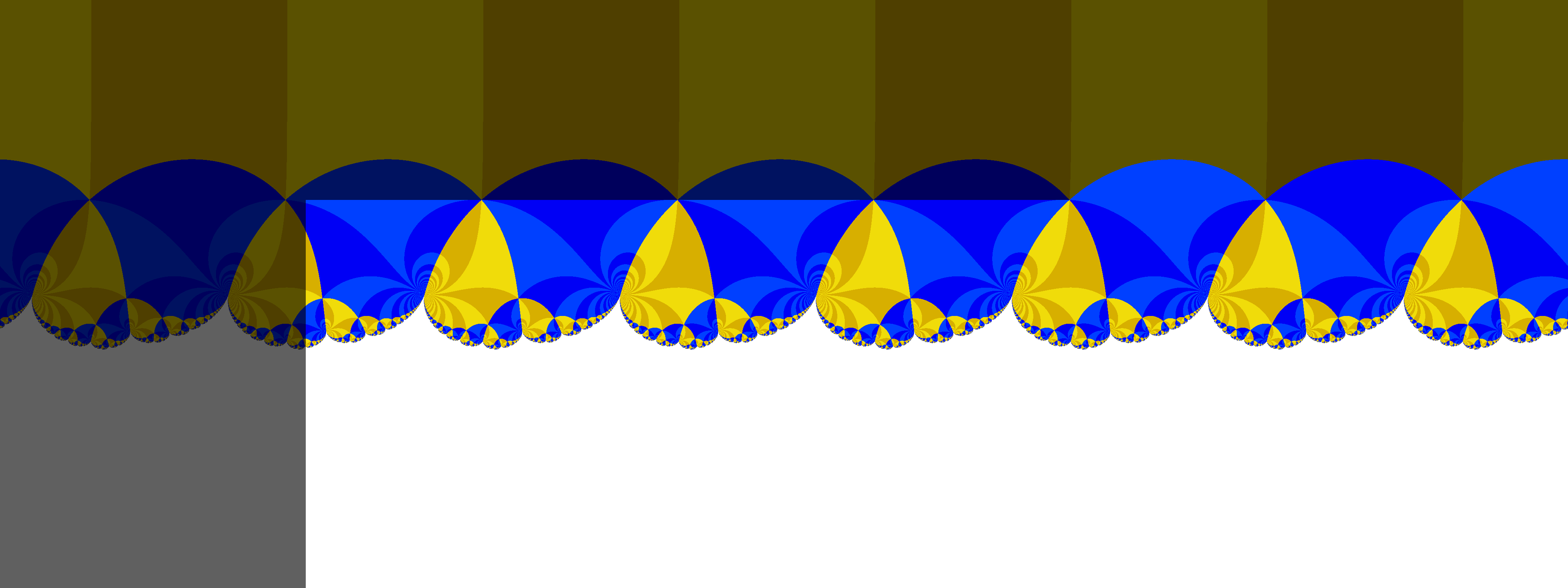}};
    \node at (-3.3,-7) {\includegraphics[width=6cm]{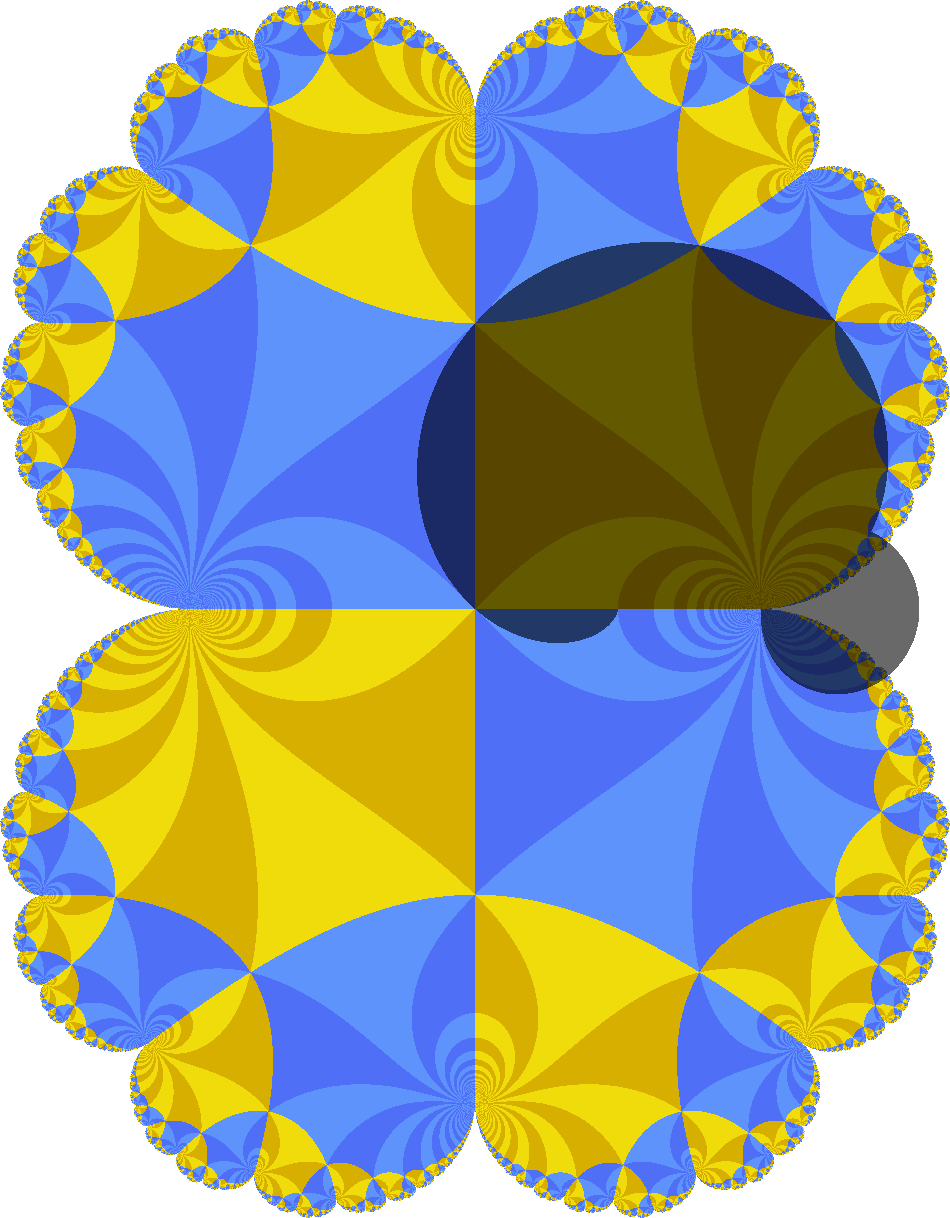}};
    \node at (3.3,-7) {\includegraphics[width=6cm]{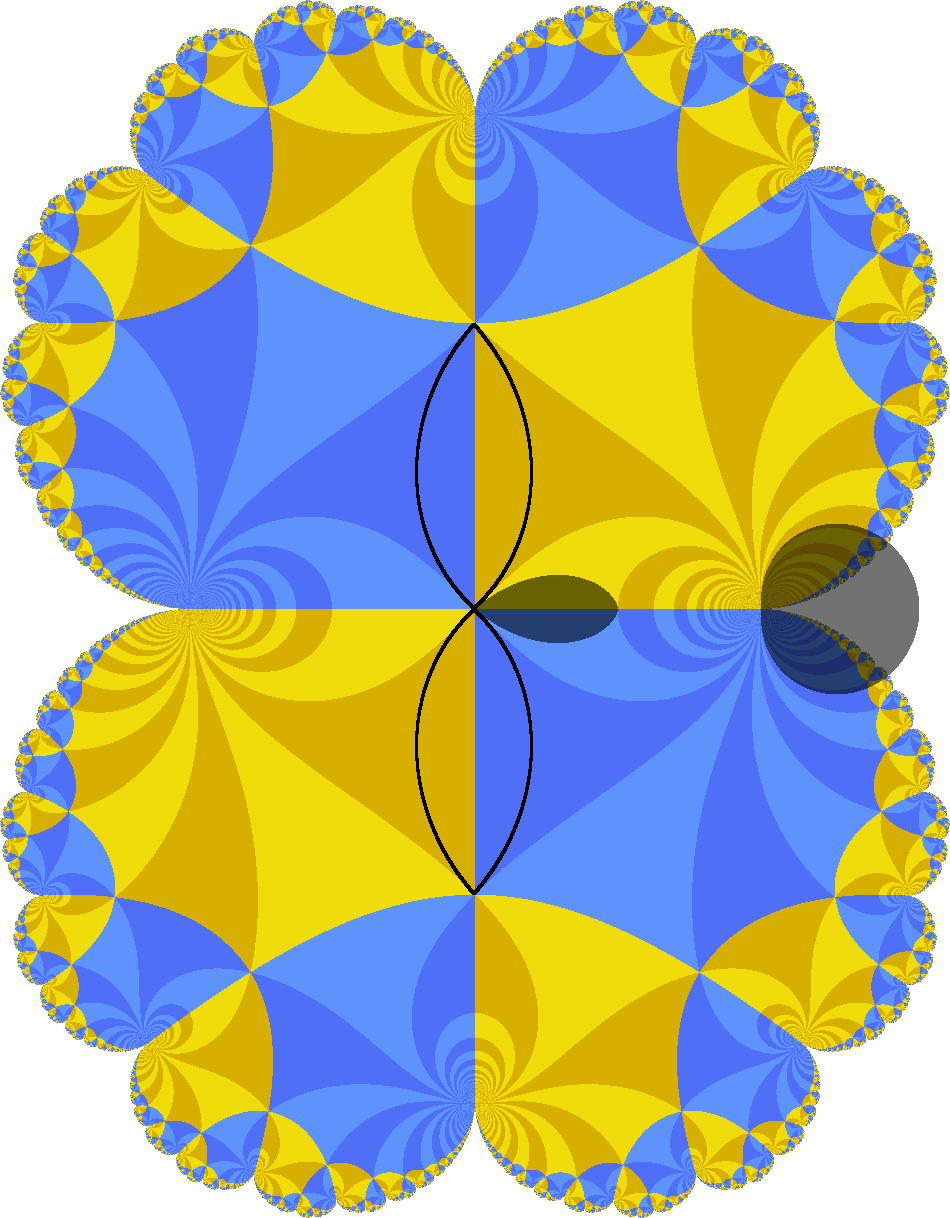}};
    \node[circle, fill=white, draw=black] at (0,1.3) {$V$};
    \node[rectangle, fill=white, draw=black] at (-2,-6) {$\Psi_R^\ext(V)$};
    \draw[->] (-3.3,-2) -- node[right]{$\Psi_R^\ext$} (-3.3,-3);
  \end{tikzpicture}
\end{center}
  \caption{Illustration of the example of \Cref{ss:ex:4:2}. Same conventions as on \Cref{fig:ex:4}. The set $P = \Psi_R^\ext(V)$ is a repelling petal. Its reflection $\ov P$ along the real axis is also a repelling petal. However, $P\cap \ov P$, which is the dark set in the lower right image has two connected components, so is not a petal. Actually, it is not even a repelling quasi-petal: the preimage by $f$ of the small left connected component is disjoint from $P\cap \ov P$.} 
  \label{fig:ex:4:2}
\end{figure}

\subsection{On repelling petals}\label{sub:ex:4}

\subsubsection{A set failing to be a repelling petal}\label{ss:ex:4:1}

We define here, with pictures, an example of set $V$ for which the converse of \Cref{lem:V} does not hold.
Let $f(z)=z+z^2$ be the Cauliflower map and $V$ be the set defined on \Cref{fig:ex:4}.
The set $P=\Psi_R^\ext(V)$ almost qualifies as a repelling petal.
It fails because when conjugating $T_{-1}:V\to V$ by $\Psi_R^\ext|_V$, we obtain a map which is discontinuous at $0$, so there is no branch of $f^{-1}$ defined on an open set containing $0$ and for which $P$ could be an attracting petal (see \Cref{def:rp}).

\subsubsection{An intersection of two repelling petals}\label{ss:ex:4:2}

\Cref{fig:ex:4:2} illustrates an example, for $f(z)=z^2+z$, of two repelling petals whose intersection is not a repelling petal, nor even a repelling quasi-petal.

\clearpage

\section{Parabolic points with several petals}\label{sec:grl}

We generalize here the definitions of semi-conjugacy, pseudo-conjugacy and the main results to parabolic fixed points with (one or) several attracting axes.

Let $f_1, f_2$ be two holomorphic functions from an open neighborhood of $0\in \C$ to $\C$ with a general parabolic point at $0$: \[ f_i(z) = e^{2i\pi p_i/q_i} z + o(z) \] 
where $p_i \in \Z$, $q_i \in \N^*$ and $p_i \wedge q_i = 1$.
We recall that the number $r_i$ of attracting axes of $f_i$ at $0$ is a multiple of $q_i$:
\[r_i=m_i q_i\]
for some positive integer $m_i$ and that
\[f^{q_i}(z) = z + a_i z^{r_i+1}+\cdots\]
We recall that an appropriate change of coordinate is $w = -1/r_i a_i z^{r_i}$, for which the dynamics is close to the translation $w\mapsto w+1$ for $w$ close to $\infty$, with the subtlety that the change of variable is now non-injective.

Consider attracting axes $A_1, A_2$ of $f_1$, $f_2$.
Let $R_i^-, R_i^+$ be the adjacent repelling axes of $A_i$ such that $R_i^-, A_i, R_i^+$ are in the trigonometrical order (if $q_i = 1$ then $R_i^- = R_i^+$). Let $B_{f_1}, B_{f_2}$ be the parabolic basins of $f_1, f_2$ and $C_1$, $C_2$ the connected components of $B_{f_1}, B_{f_2}$ containing a germ of $A_1, A_2$. Notice that $q_i$ is the least positive integer such that $f_i^{q_i}(C_i) \subset C_i$.

For $U_1, U_2$ open neighborhoods of $0 \in \C$, denote $U_i^0$ the connected component of $U_i \cap B_{f_i}^0$ containing a germ of $A_i$.
Equivalently, $U_i^0$ is the connected component of $U_i \cap C_i$ containing a germ of $A_i$.
For all attracting petal $P_A^i$ of $f_i$ associated to the axis $A_i$ and included in $U_i$, we have $P_A^i \subset U_i^0$.

\begin{figure}[!t] 
\def\svgwidth{\textwidth}
\vspace{-2cm}
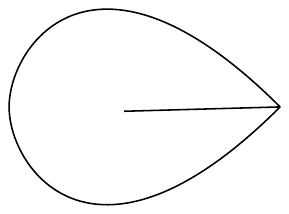
\vspace{-15cm}
\caption{General definition of pseudo-conjugacy}
\label{f:pseudoConjGen}
\vspace{0.5cm}
\footnotesize{\textit{The zones represented in gray are respectively $U_1^0$ at the left and $U_2^0$ at the right.
Note these sets are not images one of the other by $\phi$, $\phi'$.
Notice also that the angles between the attracting axis and the adjacent repelling axis are not necessarily the same, which is equivalent to say that $f_1$ and $f_2$ do not necessarily have the same number of petals.
Last, we represent here basins which are Jordan domains but the situation may by much less regular.}}
\end{figure}

See \Cref{f:pseudoConjGen} for the next two definitions.

\begin{definition}\label{d:semiConjGen}%
A holomorphic map $\phi$ is said to be a local semi-conjugacy from $(f_1^{q_1}, A_1)$ to $(f_2^{q_2}, A_2)$ on their immediate basins if there exists an open neighborhood of $0 \in \C$ such that $
\Dom (\phi) = U_1^0$, $\phi$ maps to $C_2$, and $\phi$ satisfies the semi-conjugacy relation: $\phi \circ f_1^{q_1}(z) = f_2^{q_2} \circ \phi(z)$ for all $z \in U_1^0$ such that $f_1^{q_1}(z) \in U_1^0$.
\end{definition}

As in the case of simple parabolic points, the sets $C_i / f_i^{q_i}$ are isomorphic to $\CZ$, and the map $\ov\phi: B^0_{f_1}/f_1 \rightarrow B^0_{f_2}/f_2$ is a translation too, i.e.\ the analogue of \Cref{p:SCQuotientBiholom} holds too.
For this, one sees that \Cref{lem:bz} still holds in the case of several attracting axes (actually the proof is easier in that case: $C_i$ is disjoint from any attracting petal of any axis different from $A_i$).
And in the first part of \Cref{p:SCQuotientBiholom} one replaces the equivalent by $\phi(z)\sim c z^{q_1/q_2}$ for some $c>0$ and some branch of the function $z\mapsto z^{q_1/q_2}$ on the sector $S[r]$: in its proof, apply \Cref{p:equivSect} to $z\mapsto\phi(z^{2/q_1})^{q_2/2}$ instead of $\phi$.

\begin{definition}\label{d:pseudoConjGen}%
A pair of local semi-conjugacies on immediate basins $(\phi, \phi')$, with $\phi$ from $(f_1^{q_1}, A_1)$ to $(f_2^{q_2}, A_2)$ and $\phi'$ from $(f_2^{q_2}, A_2)$ to $(f_1^{q_1}, A_1)$, is called a local pseudo-conjugacy between $(f_1^{q_1}, A_1)$ and $(f_2^{q_2}, A_2)$ if $\ov\phi : C_1/f_1^{q_1} \to C_2/f_2^{q_2}$ has for inverse $\overline{\phi'}: C_2/f_2^{q_2} \to C_1/f_1^{q_1}$.
\end{definition}

The notion of synchronous pseudo-conjugacies also extends.

\begin{remark}
Suppose that $(f_1^{q_1}, A_1)$ and $(f_2^{q_2}, A_2)$ are pseudo-conjugate, with sets $U_1, U_2$ small enough such that $f_i$ admits a local inverse at $0$: 
\[ f_i^{-1}: f_i(U_i) \to U_i\] 
Write $V_i = f_i(U_i)$, and $V_i^0$ the connected component of $V_i \cap B_{f_i}$ containing a germ of the attracting axis $Tf_i(A_i)$, where $Tf_i$ denotes the tangent map of $f_i$ at $0$. 
Saying that $(f_1^{q_1}, A_1)$ and $(f_2^{q_2}, A_2)$ are pseudo-conjugate by $(\phi, \phi')$ is equivalent to saying that $(f_1^{q_1}, Tf_1(A_1))$ and $(f_2^{q_2}, Tf_2(A_2))$ are pseudo-conjugate by $(\phi \circ f_1^{-1}, \phi' \circ f_2^{-1})$, where $\phi \circ f_1^{-1}$ and $\phi' \circ f_2^{-1}$ are defined over $V_1^0, V_2^0$. (One can check that $f_i(U_i^0)=V_i^0$)
\end{remark}

In this paragraph, we omit the index/exponent $i\in\{1,2\}$.
Note $\Phi_{A}^{ \ext  }: C \to \C$ an extended Fatou coordinate associated to the attracting axis $A$, and $\Psi_{R^\pm}^{ \ext  }$ an extended repelling Fatou parametrization associated to the repelling axis $R^\pm$.
Let
\[h^u= \Phi_{A}^{\ext} \circ \Psi_{R^+}^{\ext}\] 
\[h^d= \Phi_{A}^{\ext} \circ \Psi_{R^-}^{\ext}\] 
where $u$, $d$ stand for \emph{up} and \emph{down}, and call them the \emph{horn maps}\footnote{This definition is classical.
We specify here the domain where we want to consider horn maps to be defined.} of $(f, A)$.
These maps have respective domains $\wDd^u=\Phi_{R^{-}}^{-1}(C)$ and $\wDd^d=\Phi_{R^{+}}^{-1}(C)$ that are invariant by $T_1$, they commute with $T_1$, and we denote $\h^u$ and $\h^d$ their quotient at the domain and range by $\piZ:\C\to\CZ$, defined on $\Dd^u=\piZ(\wDd^u)$ and $\Dd^d=\piZ(\wDd^d)$ (if $q=1$ then $\h^u=\h^d$).
Denote $\Dd^+$ the connected component containing a punctured neighborhood of $+ i \infty$ of $\Dd^u$; denote $\h^+$ the restriction of $\h^u$ to $\Dd^+$.
Similarly, denote $\Dd^-$ the connected component containing a punctured neighborhood of$- i \infty$ of $\Dd^d$ and denote $\h^-$ the restriction of $\h^d$ to $\Dd^-$.
Notice that $\h^+$ and $\h^-$ both correspond to the same attracting axis $A$, but to distinct repelling axes $R^-$ and $R^+$. 
Notice that if $q>1$, since we can take arbitrary and independent normalizations for $R^+$ and $R^-$, there is no guarantee that the domains of $\h^+$ and $\h^-$ are disjoint.

The proof of \Cref{t:conjRevetBasLoc} and the proof of \Cref{prop:complement} may be followed in this more general case, Here are the analogue statements we obtain.

\begin{theorem} \label{en:semiConjLoc}
There is a local semi-conjugacy $\phi$ from $(f_1^{q_1}, A_1)$ to $(f_2^{q_2}, A_2)$ on their immediate basins if and only if
there exists $\sigma \in \CZ$ and a pair of holomorphic maps $\psi = (\psi^+, \psi^-)$, where $\psi^\pm: \Dd_1^\pm \rightarrow \Dd_2^\pm$, such that: 
\[ \h_1^\pm = T_\sigma^{-1} \circ \h_2^\pm \circ \psi^\pm
\]
\end{theorem}

\begin{theorem} \label{en:conjLoc}
The maps, together with selected attracting axes, $(f_1^{q_1}, A_1)$ and $(f_2^{q_2}, A_2)$, are locally pseudo-conjugate if and only if
there exists $\sigma \in \CZ$, and a biholomorphism pair $\psi = (\psi^+, \psi^-)$, where $\psi^\pm: \Dd_1^\pm \to \Dd_2^\pm$, such that: 
\[ \h_1^\pm = T_{\sigma}^{-1} \circ \h_2^\pm \circ \psi^\pm \]
\end{theorem}

Furthermore, for both statements, the maps $\psi^\pm$ admit at $\pm i \infty$ a removable singularity and an expansion of the form $\psi^\pm(w) = w + \rho^\pm + o(1)$ as $\Im(z)\to\pm\infty$.

The commuting diagram relating $\psi^\pm$ and $\phi$ is adapted as follows:
\[
\xymatrix@R=10pt@C=35pt{
  \wDd_1^\pm \ar[dd]_{\wt\psi^\pm} \ar[r]^{\Psi_{R^\pm}^{1,\ext}} & B_{f_1}^0  \ar[r]^{\Phi_A^{1,\ext}} & \C \ar[dd]^{T_{\tilde \sigma}}
  \\
  & U_1^0 \ar[d]_{\phi} \ar[u]^{\iota} &
  \\
  \wDd_2^\pm \ar[r]_{\Psi_{R^\pm}^{2,\ext}} & B_{f_2}^0 \ar[r]_{\Phi_A^{2,\ext}} & \C
}
\]
the only difference being that $R$ has been replaced by $R^\pm$.

Concerning the proof of the two statements above:


As already mentioned, the set $\wDd$ has to be replaced by the pair $\wDd^u$ and $\wDd^d$.

In the adaptation of \Cref{sub:dhu0}, the set $W(U)$ must be replaced by the pair of sets $W^\pm(U)$ defined as the set of $w\in\C$ such that for all $n$ big enough, $w-n\in \Phi_{R^\pm}(U^0)$.
In \Cref{p:petitPetRep}, the statement becomes $\wDd^\pm \subset \wt W^\pm(U)$.

In the reverse inclusion proved in \Cref{ssub:ci}, the condition $U_1 = P_A^1 \cup P_R^1 \cup \{0\}$ with $P_A^1 \cup P_R^1$ has two connected components is changed to $U_1$ being the union of $\alpha$-petals ($\alpha>\pi/2$) in a Leau flower, three of which being $P_A^1$, $P_{R^-}^1$ and $P_{R^+}^1$, such that any two consecutive petals have connected intersection.

\bibliographystyle{alpha} 
\bibliography{These} 

\end{document}